\documentclass[oneside,english]{amsart}
\usepackage[T1]{fontenc}
\usepackage[latin9]{inputenc}
\usepackage{geometry}
\usepackage{verbatim}
\usepackage{amstext}
\usepackage{amsthm}
\usepackage{amssymb}
\usepackage[foot]{amsaddr}
\usepackage{graphicx}
\usepackage{esint}
\usepackage{mathabx}
\usepackage{color}

\makeatletter
\numberwithin{equation}{section}
\numberwithin{figure}{section}
\theoremstyle{plain}
\newtheorem{thm}{\protect\theoremname}
\theoremstyle{plain}
\newtheorem{lem}[thm]{\protect\lemmaname}
\theoremstyle{plain}
\newtheorem{prop}[thm]{\protect\propositionname}
\theoremstyle{remark}
\newtheorem{rem}[thm]{\protect\remarkname}
\theoremstyle{remark}
\newtheorem{notation}[thm]{\protect\notationname}
\theoremstyle{definition}
\newtheorem{defn}[thm]{\protect\definitionname}

\global\long\def\Re{\operatorname{Re}}

\global\long\def\Im{\operatorname{Im}}

\global\long\def\Arg{\operatorname{Arg}}

\global\long\def\Log{\operatorname{Log}}

\global\long\def\Res{\operatorname{Res}}

\makeatother

\usepackage{babel}
\providecommand{\definitionname}{Definition}
\providecommand{\lemmaname}{Lemma}
\providecommand{\remarkname}{Remark}
\providecommand{\theoremname}{Theorem}
\providecommand{\propositionname}{Proposition}
\providecommand{\notationname}{Notation}

\begin{document}
\title[On Sheffer sequences]{On combinatorial properties and the zero distribution of certain Sheffer sequences}
\author{Gi-Sang Cheon${}^1$}
\address{${}^1$Department of Mathematics/ Applied Algebra and Optimization Research Center, Sungkyunkwan University, Suwon 16419, Rep. of Korea}
\email{gscheon@skku.edu}
\author{Tam\'{a}s Forg\'{a}cs${}^2$}
\address{${}^2$Department of Mathematics, California State University, Fresno, Fresno, CA 93740-8001, USA}
\email{tforgacs@csufresno.edu}
\author{Hana Kim${}^1$}
\author{Khang Tran${}^2$}
\email{khangt@csufresno.edu}
\begin{abstract} We present combinatorial and analytical results concerning a Sheffer sequence with a generating function of the form $G(x,z)=Q(z)^{x}Q(-z)^{1-x}$, where $Q$ is a quadratic polynomial with real zeros. By using the properties of Riordan matrices we address combinatorial properties and interpretations of our Sheffer sequence of polynomials and their coefficients. We also show that apart from two exceptional zeros, the zeros of polynomials with large enough degree in such a Sheffer sequence lie on the line $x=1/2+it$. \\

\end{abstract}
\maketitle
\setcounter{tocdepth}{3}
\tableofcontents
\section{Introduction}

A systematic study of the Sheffer sequence of polynomials was carried out by G.-C. Rota and his
collaborators as part of their development of Umbral Calculus (\cite{RR,RK}). Recall that a {\it Sheffer sequence} $(s_{n}(x))_{n\ge0}$ is uniquely associated to a pair $(g,f)$ of formal power series in $z$ generating the polynomials $s_n(x)$ of degree $n$ via the relation
\begin{eqnarray}
g(z)e^{xf(z)}=\sum_{n=0}^\infty s_{n}(x)\frac{z^{n}}{n!},\label{e:sheffer}
\end{eqnarray}
where $g$ and $f$ are invertible with respect to the product and the composition of series, respectively. 
This notion led to the development of Riordan group theory by means
of infinite lower triangular matrices called Riordan matrices or Riordan arrays generated by a pair of formal power series. It is shown in \cite{BShap} that (exponential) Riordan matrices constitute
a natural way of describing Sheffer sequences and various combinatorial situations such
as ordered trees, generating trees, and lattice paths, etc. (\cite{BMS, DFR}).

As is well-known, a number of classical polynomial sequences are Sheffer sequences, e.g., Laguerre polynomials, Bernoulli
polynomials, Hermite polynomials, Poisson-Charlier polynomials and Stirling polynomials. In addition to their combinatiorial importance, many of these sequences have also been extensively studied from an analytical perspective, and in particular, from the perspective of their zero distribution. While there is a wealth of knowledge assembled about the zeros of the classical orthogonal polynomials and other special functions (see \cite{rainville} for example), there are many Sheffer sequences whose zero distribution is not known.  On one hand, the requirement that a polynomial sequence be a Sheffer sequence restricts the type of generating function the sequence may have. On the other hand, the sequences under consideration in the current paper have generating functions that are quite different from those that have been studied in some recent works (see \cite{tk1}, \cite{tk2}, \cite{tk3}, \cite{innoc}, \cite{tran}, and \cite{tranzumba}). Broadly speaking, these works study the zero distribution of sequences $(s_n)_{n \geq 0}$ generated by certain 'rational'-type bivariate generating functions. Studying the zero distribution of Sheffer sequences $(H_n(x))_{n \geq0}$ with generating functions of the form $G(x,z)=Q(z)^{x}Q(-z)^{1-x}$, where $Q$ is a quadratic polynomial is a new contribution to the growing body of work in this area. While the basic ideas involved are standard, their implementation requires some careful asymptotic analysis and is at times tedious. The reader will be rewarded with an appealing result, namely, that Sheffer polynomials $H_n(x)$ with generating functions like $G(x,z)$ and large enough degree have all their zeros (apart from the ones at $x=0$ and $x=1$) on the 'critical line' $x=1/2+it$, $t \in \mathbb{R}$.   
\subsection{Organization} The paper consists of two main sections. Part I addresses the combinatorial properties and interpretations of our sequence of polynomials and their coefficients. After  introducing the notion of a Riordan matrix, we proceed to translate the coefficient matrix of a Sheffer sequence $(H_n(x))_{n\ge0}$ into an exponential Riordan matrix. This allows us to obtain two different possible combinatorial interpretations for the Sheffer sequence of polynomials. We arrive at the first interpretation by introducing generating trees associated to the production matrix of the exponential Riordan matrix (Theorem \ref{thm1}, Lemma \ref{lem_Q}, Theorem \ref{generating}).
Since the best known applications of Sheffer
sequences occur in enumeration problems of lattice paths, we offer a second interpretation for our Sheffer sequence in relation to weighted lattice paths (Theorem \ref{thm4}). We obtain the results of this approach by employing the Stirling transform of the sequence $(x^{n})_{n\ge0}$.
Part II provides the analysis of the zeros of the generated Sheffer sequence.  
It begins with the standard representation of the polynomial $H_n(x)$ as a line integral on a small circle around the origin using the Cauchy integral formula. In subsection \ref{sec:gammadeform} we show that this line integral is essentially the imaginary part (or $i$ times the real part) of the integral $\int_{\Gamma_2} f(z,t)dz$ for an appropriately defined function $f(z,t)$, and $\Gamma_2$ the boundary of a small, tubular neighborhood of the ray $[z_1,\infty)$. 

 In order to estimate $\int_{\Gamma_2} f(z,t)dz$, we use the saddle point method, which requires that after we write 
\[
f(z,t)=e^{-n\phi(z,t)}\psi(z,t),
\]
we identify the critical points of $\phi(z,t)$ (in $z$), as a function of $t$. These critical points will trace two curves as $0<t<T$, of which we select the one more suitable (Lemmas \ref{lem:moduluszeta}, \ref{lem:asympzetaT}, \ref{lem:asympzetaT1}, Proposition \ref{lem:dominantcrit}, Lemmas \ref{lem:imphiinc}, \ref{lem:rephizetacurve} and \ref{lem:rephiimaxis}) for the saddle point method, which we call $\zeta(t)$. 
 Sections \ref{sec:maintermapprox} and \ref{sec:tailapprox} (Lemmas \ref{lem:rephiposcomp}, \ref{lem:distinctcomp}, \ref{lem:gammae-int} and Proposition \ref{lem:existenceGamma}) are dedicated to proving that through each point of $\zeta(t)$, there is a deformation $\Gamma$ (see Figure \ref{fig:Gammacurve}) of $\Gamma_2$ with certain desirable properties vis-\'a-vis the saddle point method. In these sections we also establish that for each $t$, $\Gamma$ can be divided three segments -- two tails, and a segment centered on $\zeta(t)$ -- so that the integral on the central segment dominates those over the tails. Sections \ref{sec:smallT-t}, \ref{sec:SmallT1-t} and \ref{sec:smallt}  address ranges of $t$ for which the estimates in the preceeding sections are not good enough to ascertain that that the integral over the central segment dominates, and provide asymptotic expressions for $\int_{\Gamma_2} f(z,t)dz$ for values of $t$ in these ranges. With these expression in hand, in Section \ref{sec:zerodist} we proceed to compute the change in the argument of the integral representing $H_n(1/2+int)$ along a slightly deformed version of the closed loop $\zeta \cup -\bar \zeta \cup -\zeta \cup \bar \zeta$ which avoids the singularities of $\zeta(t)$ (Lemmas 29-35) in order to get a lower bound on the number of zeros of $H_n$ on the critical line. After accounting for a few exceptional zeros (Lemma \ref{lem:trivialzeros}) a quick computation of the degree of $H_n$ (Lemma \ref{lem:degree}) and the Fundamental Theorem of Algebra complete the proof of the main result.

\subsection*{Akcnowledgements} G.-S. Cheon was partially supported by the National Research Foundation of Korea (NRF) grant funded by the Korea government (MSIP) (2016R1A5A1008055 and 2019R1A2C1007518). T. Forg\'acs and K. Tran would like to acknowledge the research support of the California State University, Fresno.

\section{Part I - Combinatorial results concerning Sheffer sequences}

We begin by briefly describing the notion of a Riordan array by focusing on exponential Riordan matrices, as these are ones we mostly use in this paper.

An infinite lower triangular matrix $A=[a_{n,k}]_{n,k\ge0}$ is called
a \textit{Riordan matrix} if its $k$th column has generating function
$gf^{k}$ for some $g,f\in{\mathbb{C}}[[z]]$, where $f(0)=0$. If,
in addition, $g(0)\ne0$ and $f'(0)\ne0$ then $A$ is invertible
and it is said to be a proper Riordan matrix. We may write $A=(g,f)$.
If the $k$th column of $A$ has an exponential generating function $g(z)\frac{f(z)^{k}}{k!}$,
then $A$ is called an \textit{exponential} Riordan matrix and we denote it by $A=[g,f]$. By definition, if $A=(g,f)=[a_{n,k}]_{n,k\ge0}$
and $B=[g,f]=[b_{n,k}]_{n,k\ge0}$, then it is obvious that 
\begin{eqnarray}
b_{n,k}=\frac{n!}{k!}a_{n,k}\;\;{\rm or}\;\;B=EAE^{-1}\label{e:DAD}
\end{eqnarray}
where $E={\rm diag}(0!,1!,2!,\ldots)$ is the diagonal matrix. A well-known
\textit{fundamental property} of a Riordan matrix is that if $A=[g,f]$
and $b^{T}=(b_{n})_{n\ge0}$ is generated by exponential function
$b(z)$, then the sequence $Ab$ has exponential generating function
given by $gb(f)$; we simply write this property as $[g,f]b=gb(f)$.
In particular, if $b=e^{xz}$ then
\[
[g,f]e^{xz}=ge^{xf}.
\]
 From the definition it follows
at once that if $s_{n}(x)=\sum_{k=0}^{n}s_{n,k}x^{k}$, then the Sheffer sequence $(s_{n}(x))_{n\ge0}$
can be rewritten as a matrix product:
\begin{eqnarray}
(s_{0}(x),s_{1}(x),\ldots)^{T}=[s_{n,k}](1,x,x^{2},\ldots)^{T},\label{e:matrix}
\end{eqnarray}
where $[s_{n,k}]_{n,k\ge0}$ is a lower triangular matrix as the coefficient
matrix of $(s_{n}(x))_{n\ge0}$. Since the sequence $(1,x,x^{2},\ldots)$ has exponential generating
function $e^{xz}$, by the fundamental property it follows from \eqref{e:sheffer} and \eqref{e:matrix} that
$(s_{n}(x))_{n\ge0}$ is a Sheffer sequence for $(g,f)$
if and only if its coefficient matrix $[s_{n,k}]_{n,k\ge0}$ is an
exponential Riordan matrix given by $[g,f]$. It may be also shown
that if $A=[g,f]$ and $B=[h,\ell]$ then the product is given by
$AB=[gh(f),\ell(f)]$. Moreover, $I=[1,z]$ is the usual identity
matrix and $[g,f]^{-1}=[1/g(\overline{f}),\overline{f}]$ where $\overline{f}$
is the compositional inverse of $f$, i.e., $\overline{f}(f)=f(\overline{f})=z$.

We now turn to the polynomial sequence $(H_{n}(x))_{n\ge0}$ generated
by a bivariate function 
\[
G(x,z):=Q(z)^{x}Q(-z)^{1-x},
\]
where $Q(z)$ is a quadratic polynomial whose zeros are real, $Q(0)\ne0$
and $Q'(0)\ne0$. Since 
\[
Q(z)^{x}Q(-z)^{1-x}=Q(-z)\exp\left({x\ln\frac{Q(z)}{Q(-z)}}\right),
\]
the sequence $(H_{n}(x))_{n\ge0}$ can be considered as a Sheffer
sequence with the coefficient matrix: 
\[
A:=[a_{i,j}]_{i,j\ge0}=\left[Q(-z),\ln\frac{Q(z)}{Q(-z)}\right].
\]
Using the fundamental property, we thus have 
\begin{eqnarray}
\sum_{n=0}^{\infty}H_{n}(x)\frac{z^{n}}{n!}=Q(z)^{x}Q(-z)^{1-x}=\left[Q(-z),\ln\frac{Q(z)}{Q(-z)}\right]e^{xz}.\label{e:H}
\end{eqnarray}

In this section, we are interested in finding combinatorial interpretations
for the Sheffer polynomials $H_{n}(x)=\sum_{k=0}^{n}a_{n,k}x^{k}$
of degree $n$ where $A=[a_{i,j}]_{i,j\ge0}$. For combinatorial counting
purposes, our interest is the polynomials with non-negative integer
coefficients $a_{n,k}$. Throughout this section, we assume that $Q(z)=(1+az)(1+bz)$
and $a,b>0$ are integers. We first note that the coefficient matrix
$A$ might have negative entries, but $\widehat{A}:=\left[Q(z),\ln\frac{Q(z)}{Q(-z)}\right]$
has no negative entries whenever $a,b>0$. Since 
\[
Q(z){\exp\left({x\ln\frac{Q(z)}{Q(-z)}}\right)}=Q(z)^{1+x}Q(-z)^{-x},
\]
the Sheffer sequence associated to $\widehat{A}$, say $(\widehat{H}_{n}(x))_{n\ge0}$
has the bivariate generating function given by $G(-x,-z)$.
In addition, using Riordan multiplication we obtain $LD$-decomposition
$\widehat{A}=L_{Q}D$, where 
\begin{eqnarray}
L_{Q}=[q_{n,k}]_{n,k\ge0}=\left[Q(z),\frac{1}{2(a+b)}\ln\frac{Q(z)}{Q(-z)}\right]\label{e:LQ}
\end{eqnarray}
is a unit lower triangular matrix with ones on the main diagonal,
and $D=[1,2(a+b)z]$ is a diagonal matrix of the form ${\rm diag}(1,(2a+2b),(2a+2b)^{2},\ldots)$.
If we use the notation $[x^{k}]$ for the coefficient extraction operator,
we have for $k=0,1,\ldots,n$, 
\begin{eqnarray}
[x^{k}]H_{n}(x)=(-1)^{n+k}[x^{k}]\widehat{H}_{n}(x)=(-1)^{n+k}(2a+2b)^{k}q_{n,k},\;\;q_{0,0}=1.\label{e:Qnk}
\end{eqnarray}
A few rows of the matrix $L_{Q}$ are displayed by 
\begin{eqnarray}
L_{Q}=\left(\begin{array}{cccccc}
1 & 0 & 0 & 0 & 0\\
a+b & 1 & 0 & 0 & 0\\
2ab & 2(a+b) & 1 & 0 & 0 & \cdots\\
0 & 2(a+b)^{2} & 3(a+b) & 1 & 0\\
0 & 8(a^{3}+b^{3}) & 4(2a^{2}+2b^{2}+ab) & 4(a+b) & 1\\
\vdots &  & \cdots &  &  & \ddots
\end{array}\right).\label{e:LQ-1}
\end{eqnarray}
Thus from \eqref{e:Qnk} and \eqref{e:LQ-1} we obtain: 
\begin{eqnarray*}
H_{0}(x) & = & 1,\\
H_{1}(x) & = & -(a+b)+2(a+b)x,\\
H_{2}(x) & = & 2ab-4(a+b)^{2}x+4(a+b)^{2}x^{2},\\
H_{3}(x) & = & 4(a+b)^{3}\left(x-3x^{2}+2x^{3}\right),\\
 & \cdots
\end{eqnarray*}

The following theorem is useful for finding the combinatorial interpretation
for $L_{Q}$.
\begin{thm}
\label{thm1}(\cite{CJB,DFR}) Let $A=[a_{n,k}]$ be an infinite lower
triangular matrix with $a_{0,0}=1$. Then $A$ is an exponential Riordan
matrix given by $A=[g,f]$ if and only if there exists a horizontal
pair $\{c_{n},r_{n}\}_{n\ge0}$ of the sequences with $c=\sum_{n\ge0}c_{n}z^{n}$
and $r=\sum_{n\ge0}r_{n}z^{n}$ such that 
\begin{equation}
c(f(z))=\frac{g^{\prime}(z)}{g(z)}\quad\text{and}\quad r(f(z))=f^{\prime}(z),\label{eq:cb}
\end{equation}
or for all $n\ge k\ge0$, 
\begin{eqnarray}
a_{n+1,k}=\frac{1}{k!}\sum_{i=k-1}^{n}i!(c_{i-k}+kr_{i-k+1})a_{n,i},\quad(c_{-1}:=0).\label{e:rec}
\end{eqnarray}
\end{thm}

At times the exponential generating function approach provides explicit
forms for various (increasing) tree counting problems. Under various
guises, such trees have surfaced as tree representations of permutations,
as data structures in computer science, and as probabilistic models
in diverse applications. There is a unified generating function approach
to the enumeration of parameters on such trees (\cite{BFS,MM}). Indeed,
it was shown in \cite{BFS} that the counting generating functions for
several basic parameters, e.g. root degree, number of leaves, path
length, and level of nodes, are related to a simple ordinary differential
equation: 
\begin{eqnarray}
Y'(z)=\varphi(Y(z)),\quad Y(0)=0.\label{e:ode}
\end{eqnarray}

Comparing this differential equation with the second equation in (\ref{eq:cb}),
we see that an exponential Riodan matrix with non-negative integer
entries is closely related to tree counting problems. For instance,
if we consider $\varphi(z)$ to be the generating function for the degree-weight
sequence under a `nonnegativity' condition, then $Y(z)=\ln\frac{Q(z)}{Q(-z)}$
satisfying \eqref{e:ode} can be regarded as the generating function
for total weights of certain simple family of increasing trees.

The formula in equation \eqref{e:rec} can be rewritten as a matrix equality,
$AP_{A}=UA$ where $U$ is the upper shift matrix with ones only on
the superdiagonal and zeros elsewhere, and $P_{A}:=[p_{i,j}]_{i,j\ge0}$
where 
\begin{eqnarray}
p_{i,j}=\frac{i!}{j!}(c_{i-j}+jr_{i-j+1}),\;\;c_{i},r_{i}=0\;{\rm for}\;i<0.\label{e:prod}
\end{eqnarray}
We call $P_{A}$ the production matrix of $A$, or sometimes the Stieltjes
transform of $A$.
\begin{lem}
\label{lem_Q} Let $L_{Q}=[g,f]$ denote the exponential Riordan matrix
given by \eqref{e:LQ}. Then the horizontal pair $\{c_{n},r_{n}\}_{n\ge0}$
of $L_{Q}$ is given by 
\begin{eqnarray}
c(z) & = & \frac{a}{1+a\overline{f}}+\frac{b}{1+b\overline{f}},\label{e:cr}\\
r(z) & = & \frac{1}{a+b}\left(\frac{a}{1-(a\overline{f})^{2}}+\frac{b}{1-(b\overline{f})^{2}}\right),
\end{eqnarray}
where 
\begin{eqnarray*}
\overline{f}(z)=\frac{a+b}{2ab}\left({\rm coth}\;(a+b)z-\sqrt{{\rm coth}^{2}(a+b)z-\frac{4ab}{(a+b)^{2}}}\right).
\end{eqnarray*}
\end{lem}

\begin{proof} Note that $f(z)=\frac{1}{2(a+b)}\ln\frac{1+(a+b)z+abz^{2}}{1-(a+b)z+abz^{2}}$.
Since $f(\overline{f})=z=\overline{f}(f)$, we have 
\[
1+(a+b){\overline{f}}+ab{\overline{f}}^{2}=e^{2(a+b)z}\left(1-(a+b){\overline{f}}+ab{\overline{f}}^{2}\right).
\]
Solving the above quadratic equation we obtain 
\begin{eqnarray*}
\overline{f}(z) & = & \frac{a+b}{2ab}\left({\rm coth}\;(a+b)z-\sqrt{{\rm coth}^{2}(a+b)z-\frac{4ab}{(a+b)^{2}}}\right)\\
 & = & z-\left(\frac{2(a^{3}+b^{3})}{a+b}\right)\frac{z^{3}}{3!}+8\left(\frac{2(a^{5}+b^{5})}{a+b}-5ab(a-b)^{2}\right)\frac{z^{5}}{5!}-+\cdots.
\end{eqnarray*}
Thus it follows from (\ref{eq:cb}) that 
\begin{eqnarray*}
c(z) & = & \frac{g^{\prime}(\overline{f})}{g(\overline{f})}=\frac{(a+b)+2abz}{1+(a+b)z+abz^{2}}\mid_{z=\overline{f}}=\frac{a}{1+a\overline{f}}+\frac{b}{1+b\overline{f}}\\
 & = & (a+b)-(a^{2}+b^{2})z+2(a^{3}+b^{3})\frac{z^{2}}{2!}-2(a+b)^{2}(ab+2(a-b)^{2})\frac{z^{3}}{3!}+-\cdots,
\end{eqnarray*}
and 
\begin{eqnarray*}
r(z) & = & f^{\prime}(\overline{f})=\frac{1}{a+b}\left(\frac{a}{1-(a\overline{f})^{2}}+\frac{b}{1-(b\overline{f})^{2}}\right)\\
 & = & 1+2((a-b)^{2}+ab)\frac{z^{2}}{2!}+8\left(a^{4}+b^{4}+ab((a-b)^{2}-ab)\right)\frac{z^{4}}{4!}+\cdots.
\end{eqnarray*}
\end{proof}
In particular, if $a=b$ then we obtain: 
\begin{eqnarray*}
\overline{f}(z) & = & \frac{1}{a}\left(\frac{e^{2az}-1}{e^{2az}+1}\right)=\frac{1}{a}\sum_{n=1}^{\infty}\frac{2^{2n}(2^{2n}-1)B_{2n}(az)^{2n-1}}{(2n)!},\\
c(z) & = & a\left(1+e^{-2az}\right)=a\left(1+\sum_{n=0}^{\infty}\frac{(-2az)^{n}}{n!}\right),\\
r(z) & = & \left(\frac{e^{az}+e^{-az}}{2}\right)^{2}=1+\sum_{n=1}^{\infty}2^{2n-1}\frac{(az)^{2n}}{(2n)!},
\end{eqnarray*}
where $B_{n}$ is the $n$th Bernoulli number. \bigskip{}

Let $P_{L_{Q}}=[p_{i,j}]_{i,j\ge0}$ denote the production matrix
of $L_{Q}$ in \eqref{e:LQ}. Then 
\begin{eqnarray}
P_{L_{Q}}=\left(\begin{array}{ccccc}
a+b & 1 & 0 & 0\\
-(a^{2}+b^{2}) & a+b & 1 & 0\\
2(a^{3}+b^{3}) & -2ab & a+b & 1 & \cdots\\
-2(a+b)^{2}(ab+2(a-b)^{2}) & 6(a^{3}+b^{3}) & 3(a-b)^{2} & a+b\\
\vdots & \vdots & \vdots &  & \ddots
\end{array}\right).\label{e:production}
\end{eqnarray}

It is well-known \cite{BMS,DFR} that if a production matrix $P_{A}$
is an integer matrix then every element of the Riordan matrix $A$
has a combinatorial interpretation of counting marked or non-marked
nodes in the associated generating tree. A \textit{marked generating
tree} \cite{BMS} is a rooted labeled tree with the property that
if $v_{1}$ and $v_{2}$ are any of two nodes with the same label
$k$ then for each label $k=0,1,2,\ldots$, $v_{1}$ and $v_{2}$
have the same number of children. The nodes $(k)$
with label $k$ may or may not be marked depending on whether an element
in $P_{A}$ is negative or positive. To specify a generating tree
it therefore suffices to specify: 
\begin{itemize}
\item[(a)] the label of the root; 
\item[(b)] a set of production rules explaining how to derive the quantity of
children and their labels, from the label of a parent. 
\end{itemize}
Noticing that the $k$th row of the production matrix of a Riordan
matrix defines a production rule for the node $(k)$,
we can similarly associate an exponential Riordan matrix with its
production matrix to a marked generating tree specification using the notation
\[
(k)^{p}=\underbrace{(k)\cdots(k)}_{p\;{\rm times}}\;(p>0)\;\;{\rm and}\;\;(\overline{k})^{p}=\underbrace{(\overline{k})\cdots(\overline{k})}_{-p\;{\rm times}}\;(p<0)
\]
where ${(k)}^{0}$ is the empty sequence and $(\overline{\overline{k}})=(k)$.
We are now ready to give a combinatorial interpretation for coefficients
of the polynomials $H_{n}(x)$ by means of a marked generating tree.
Note that it follows from \eqref{e:Qnk} that $\frac{(-1)^{n+k}}{(2a+2b)^{k}}[x^{k}]H_{n}(x)=q_{n,k}$
where $L_{Q}=[q_{i,j}]_{i,j\ge0}$.
\begin{thm}\label{generating}
For $k=0,1,\ldots,n$, let $\mu_{n}(k)$ denote the number of nodes
$(k)$ with label $k$ at level $n$ in the marked
generating tree specification where the root is at level $0$:
\begin{align}
\left\{ \begin{tabular}{ll}
 root :  & \ensuremath{(0)}\\
 rule :  & \ensuremath{(k)\rightarrow (0)^{p_{k,0}}(1)^{p_{k,1}}\cdots(k+1)^{p_{k,k+1}}}\\
  &  \ensuremath{(\overline{k})\ensuremath{\rightarrow}(\overline{0})^{p_{k,0}}(\overline{1})^{p_{k,1}}\cdots(\overline{k+1})^{p_{k,k+1}}}
\end{tabular}\right.\label{eq:gts}
\end{align}

where $p_{k,j}={\frac{k!}{j!}(c_{k-j}+jr_{k-j+1})}$ for the horizontal
pair $\{c_{n},r_{n}\}_{n\ge0}$ in Lemma \ref{lem_Q}. Then
\begin{eqnarray}
[x^{k}]H_{n}(x)=(-1)^{n+k}(2a+2b)^{k}\left(\mu_{n}(k)-\mu_{n}(\overline{k})\right).\label{e:interpretation-1}
\end{eqnarray}
\end{thm}

\begin{proof} Let $P_{L_{Q}}=[p_{i,j}]_{i,j\ge0}$ be
the production matrix of $L_{Q}$. Since $L_{Q}$ is a unit lower
triangular matrix of nonnegative integers for $a,b>0$, it is obvious
from $P_{L_{Q}}=L_{Q}^{-1}UL_{Q}$ that all entries of $P_{L_{Q}}$
are integers. Thus it follows from \eqref{e:prod} that $p_{k,j}={\frac{k!}{j!}(c_{k-j}+jr_{k-j+1})} \in \mathbb{Z}$, where $(c_{n})_{n\ge0}$ and $(r_{n})_{n\ge0}$
are the sequences obtained from Theorem \ref{thm1}. In addition,
for $n,k\ge0$ we have 
\begin{eqnarray*}
q_{0,0}=1,\;\;q_{n+1,k}=\sum_{i=k-1}^{n}\frac{i!}{k!}(c_{i-k}+kr_{i-k+1})q_{n,i}=\sum_{i=k-1}^{n}q_{n,i}p_{i,k},\quad(c_{-1}:=0).
\end{eqnarray*}

By using the succession rule (\ref{eq:gts}), we show that $q_{n,k}=\mu_{n}(k)-\mu_{n}(\overline{k})$.
In the marked generating tree at level zero we have only one node
$(0)$ with label 0. This is represented by the row
vector $R_{0}=(1,0,\ldots)$. At the next levels of the generating
tree, the distribution of the nodes $(1),(2),\ldots$
is given by the row vectors $R_{i}$, $i\ge1$, defined by the recurrence
relation $R_{i}=R_{i-1}P_{L_{Q}}$. Stacking these row vectors we
obtain the matrix $L_{Q}=[R_{0},R_{1},\ldots]^{T}$ satisfying $L_{Q}P_{L_{Q}}=UL_{Q}$.
Since the marked nodes kill or annihilate the non-marked nodes with
the same number in this process, it follows that $q_{n,k}$ counts
the difference between the number of non-marked nodes $(k)$
and the number of marked nodes $(\overline{k})$
at level $n$. Hence \eqref{e:interpretation-1} immediately follows
from \eqref{e:Qnk}, as required.
\end{proof}
The best known applications of Sheffer
sequences occur in enumeration problem of lattice paths (see \cite{HN}). We propose
now a second combinatorial interpretation for coefficients of the
polynomials $H_{n}(x)$ by weighted lattice paths. This approach can
be obtained from the Stirling transform of the sequence $(x^{k})_{k\ge0}$
given by 
\begin{eqnarray}
x^{n}=\sum_{k=0}^{n}S(n,k)(x)_{k},\label{e:stirling}
\end{eqnarray}
where $S(n,k)$ is the Stirling number of the second kind and $(x)_{k}:=\prod_{i=0}^{k-1}(x-i)$
is the $k$th falling factorial with $(x)_{0}=1$. Using the exponential
generating functions $e^{xz}$ for $\left((x^{n})\right)_{n\ge0}$,
and $(1+z)^{x}$ for $\left((x)_{n}\right)_{n\ge0}$, we obtain from
\eqref{e:stirling} that $e^{xz}=[1,e^{z}-1](1+z)^{x}$ where $[1,e^{z}-1]$
is the Stirling matrix of the second kind whose $(i,j)$-entry is
$S(i,j)$. Hence we obtain 
\begin{eqnarray}\label{e:QS}
\left[Q(-z),\ln\frac{Q(z)}{Q(-z)}\right]e^{xz} & = & \left[Q(-z),\ln\frac{Q(z)}{Q(-z)}\right]\left[1,e^{z}-1\right](1+z)^{x}\nonumber \\
 & = & \left[Q(-z),\frac{Q(z)}{Q(-z)}-1\right](1+z)^{x}.
\end{eqnarray}
It folllows from \eqref{e:H}, \eqref{e:stirling} and \eqref{e:QS}
that 
\begin{eqnarray}
H_{n}(x)=\sum_{k=0}^{n}a_{n,k}x^{k}=\sum_{k=0}^{n}\left(\sum_{i=0}^{n}a_{n,i}S(i,k)\right)(x)_{k}=\sum_{k=0}^{n}c_{n,k}(x)_{k},\label{e:IC}
\end{eqnarray}
where $A=[a_{n,k}]_{n,k\ge0}=\left[Q(-z),\ln\frac{Q(z)}{Q(-z)}\right]$
and $C:=[c_{n,k}]_{n,k\ge0}=\left[Q(-z),\frac{Q(z)}{Q(-z)}-1\right]$.
A few rows of $C$ shown below:
\begin{eqnarray}
\qquad \qquad C=\left[\begin{array}{cccccc}
1 & 0 & 0 & 0 & 0\\
-(a+b) & 2(a+b) & 0 & 0 & 0\\
2ab & 0 & 4(a+b)^{2} & 0 & 0 & \cdots\\
0 & 0 & 12(a+b)^{3} & 8(a+b)^{3} & 0\\
0 & 0 & 48(a+b)^{2}(a^{2}+ab+b^{2}) & 64(a+b)^{4} & 16(a+b)^{4}\\
\vdots &  & \cdots &  &  & \ddots
\end{array}\right].\label{e:C}
\end{eqnarray}

We show that every element of the matrix $C$ can be represented in
terms of the weights of some lattice path. For this purpose, consider
a weighted lattice path in the plane from $(0,0)$ to $(n,k)$ with
up-steps $U=(1,1)$, and level-steps $H=(1,0)$ or double level-steps
$H^{2}=(2,0)$. We denote by $\omega(U)$ the weight of $U$, and
by $\omega_{ij}(H)$ and $\omega_{ij}(H^{2})$ the weights of $H$
for $(i-1,j)\rightarrow(i,j)$ and $H^{2}$ for $(i-2,j)\rightarrow(i,j)$,
respectively. As usual, the weight of a lattice path is defined by
the product of all weights assigned to the steps along the
path.
\begin{thm}\label{thm4}
For $k=0,1,\ldots,n$, let $\sigma(n,k)$ denote the sum of the weights
of lattice paths in the plane from $(0,0)$ to $(n,k)$ with the steps
in $\{U,H,H^{2}\}$ where $w(D)=2(a+b)$, $w_{ij}(H)=(i+j-2)(a+b)$,
and $w_{ij}(H^{2})=-(i-1)(i+2j-4)ab$. Then 
\[
\sigma(n,k)=c_{n,k}=[z^{n}]Q(-z)\left(\frac{Q(z)}{Q(-z)}-1\right)^{k}.
\]
In particular, 
\begin{eqnarray}
[x^{k}]H_{n}(x)=\sum_{i=0}^{n}(-1)^{i+k}\sigma(n,i)s(i,k)\label{e:interpretation-2}
\end{eqnarray}
where $s(i,k)$ is the Stirling number of the first kind. 
\end{thm}

\begin{proof} Let $H_{n}(x)=\sum_{k=0}^{n}c_{n,k}(x)_{k}$. Consider
a weighted lattice path from $(0,0)$ to $(n,k)$ with the step set
$\{U,H,H^{2}\}$. Since a lattice point $(n,k)$ may be approached
from any of the lattice points $(n-1,k-1)$, $(n-1,k)$, or $(n-2,k)$,
it is immediate that $\sigma(n,k)=c_{n,k}$ if and only if $c_{n,k}$
satisfies the following recurrence relation for $n\ge2$ and $k\ge1$:
\begin{eqnarray}
c_{n,k}=2(a+b)c_{n-1,k-1}+(n+k-2)(a+b)c_{n-1,k}-(n-1)(n+2k-4)abc_{n-2,k},\label{e:recurrence}
\end{eqnarray}
with the initial conditions $c_{0,0}=1,c_{1,0}=-(a+b),c_{1,1}=2(a+b),c_{2,0}=2ab$,
and $c_{n,0}=0$ for $n\ge3$. It is clear that \eqref{e:recurrence}
together with the initial conditions determines $c_{n,k}$.

For simplicity, we substitute by $d_{n,k}=\frac{k!}{n!}c_{n,k}$.
Then the recurrence \eqref{e:recurrence} is equivalent to 
\begin{align}
nd_{n,k}=2(a+b)kd_{n-1,k-1}+(n+k-2)(a+b)d_{n-1,k}-(n+2k-4)abd_{n-2,k},\label{eq:rec1}
\end{align}
where $d_{0,0}=1,d_{1,0}=-(a+b),d_{1,1}=2(a+b),d_{2,0}=ab$, and $d_{n,0}=0$
for $n\ge3$. If we put $\varphi_{k}(z)=\sum_{n\ge0}d_{n,k}z^{n}$,
it follows from (\ref{eq:rec1}) that 
\begin{eqnarray*}
\varphi_{k}^{\prime} & = & 2(a+b)k\varphi_{k-1}+(a+b)(z\varphi_{k})^{\prime}+(k-2)(a+b)\varphi_{k}-ab(z^{2}\varphi_{k})^{\prime}-ab(2k-4)z\varphi_{k}\\
 & = & 2(a+b)k\varphi_{k-1}+((a+b)z-abz^{2})\varphi_{k}^{\prime}+(k-1)(a+b-2abz)\varphi_{k}.
\end{eqnarray*}
Using $Q(z)=1+(a+b)z+abz^{2}$ we obtain the following differential equation for $\varphi_{k}$: 
\begin{eqnarray}
Q(-z)\varphi_{k}^{\prime}-(k-1)Q^{\prime}(-z)\varphi_{k}=2(a+b)k\varphi_{k-1},\;\;\varphi_{0}=Q(-z).\label{e:rec2}
\end{eqnarray}
If $k=1$ then $\varphi_{1}^{\prime}=2(a+b)$. Since $\varphi_{1}(0)=0$,
we get that 
\[
\varphi_{1}=2(a+b)z=Q(z)-Q(-z)=Q(-z)\left(\frac{Q(z)}{Q(-z)}-1\right).
\]
If $k=2$ then it follows from \eqref{e:rec2} that 
\[
Q(-z)\varphi_{2}^{\prime}-Q^{\prime}(-z)\varphi_{2}=\left(Q(-z)\varphi_{2}\right)^{\prime}=4(a+b)\varphi_{1}=8(a+b)^{2}z.
\]
Since $\varphi_{2}(0)=0$, we find in this case that $Q(-z)\varphi_{2}=4(a+b)^{2}z^{2}=(Q(z)-Q(-z))^{2}.$
Thus 
\[
\varphi_{2}=\frac{(Q(z)-Q(-z))^{2}}{Q(-z)}=Q(-z)\left(\frac{Q(z)}{Q(-z)}-1\right)^{2}.
\]
Repeated application of recurrence in the differential equation \eqref{e:rec2}
gives 
\begin{eqnarray}
\varphi_{k}=Q(-z)\left(\frac{Q(z)}{Q(-z)}-1\right)^{k}.\label{e:dgen}
\end{eqnarray}
Since $\varphi_{k}(z)=\sum_{n\ge0}d_{n,k}z^{n}$ and $d_{n,k}=0$
for $n<k$, the array $D:=[d_{n,k}]_{n,k\ge0}$ is a lower triangular
matrix with $\varphi_{k}$ in \eqref{e:dgen} as the $k$th column
generating function for $k\ge0$. By definition, $D$ is a Riordan
matrix given by $\left(Q(-z),\frac{Q(z)}{Q(-z)}-1\right)$. Since
$c_{n,k}=\frac{n!}{k!}d_{n,k}$ it immediately follows from \eqref{e:DAD}
that 
\[
[c_{n,k}]_{n,k\ge0}=\left[Q(-z),\frac{Q(z)}{Q(-z)}-1\right]
\]
is the exponential Riordan matrix, which agrees with the coefficient
matrix of $H_{n}(x)$ with respect to the basis $\{1,(x)_{1},(x)_{2},\ldots\}$.
Hence $\sigma(n,k)=c_{n,k}$, as required. In particular, since $(x)_{n}=\sum_{k=0}^{n}(-1)^{n+k}s(n,k)x^{n}$,
using the inverse Stirling transform we obtain from \eqref{e:IC}
that 
\[
H_{n}(x)=\sum_{k=0}^{n}c_{n,k}(x)_{k}=\sum_{k=0}^{n}\left(\sum_{i=0}^{n}(-1)^{i+k}c_{n,i}s(i,k)\right)x^{k},
\]
which completes the proof.
\end{proof}

In concluding this section, we remark that if the polynomials in
a Sheffer sequence have non-negative integer coefficients, combinatorially the sequence may be  interpreted as corresponding to a (pair of) generating functions counting increasing
trees, generating trees, or weighted lattice paths, etc. In terms
of enumeration, one may consider a Sheffer sequence as allowing several
different kind of trees or lattices paths.

\section*{Part II - the zeros of the Sheffer sequence}
We now turn our attention to the zero distribution of the Sheffer sequence $\left(H_{n}(x)\right) _{n \geq 0}$ with generating function $Q(z)^{x}Q(-z)^{1-x}$,  where $Q(z)$ is a quadratic polynomial whose zeros are real, $Q(0)\ne0$ and $Q'(0)\ne0$The main result of this part of the paper is the following theorem.
\begin{thm}
\label{thm:maintheorem}For $z_{2}>z_{2}>0$, let $Q(z)=(z_{1}-z)(z_{2}-z)$
and $\left(H_{n}(x)\right) _{n \geq 0}$ be the sequence of
polynomials generated by 
\begin{equation}
\sum_{n=0}^{\infty}H_{n}(x)\frac{z^{n}}{n!}=Q(z)^{x}Q(-z)^{1-x} \qquad (z,x \in \mathbb{C}).\label{eq:genfunc}
\end{equation}
Then for all large $n$, other than the two trivial zeros at $x=0,1$, all the
zeros of $H_{n}(x)$ lie on the critical line $\Re x=1/2$. 
\end{thm}

\subsection{Deforming the path of integration}\label{sec:gammadeform} The substitution $x=1/2+int$ and the Cauchy integral formula give
\begin{align*}
H_{n}\left(\frac{1}{2}+int\right) & =\frac{n!}{2\pi i}\ointctrclockwise_{|z|=\epsilon}\frac{Q(z)^{1/2+int}Q(-z)^{1/2-int}}{z^{n+1}}dz, \qquad (n \in \mathbb{N} \cup \{0\}). 
\end{align*}
The integrand (as a function of $z$) has an analytic continuation to the complement of  $\{0\}\cup[z_{1},\infty)\cup(-\infty,-z_{1}]$ defined by 
\begin{eqnarray*}
f(z,t)&:=&\frac{1}{z^{n+1}}\exp\left(\left(\frac{1}{2}+int\right)\Log(z_{1}-z)+\left(\frac{1}{2}+int\right)\Log(z_{2}-z) \right.\\
&+&\left.\left(\frac{1}{2}-int\right)\Log(z_{1}+z)+\left(\frac{1}{2}-i n t\right)\Log(z_{2}+z)\right),
\end{eqnarray*}
where $\Log z$ denotes the principal logarithm\footnote{In the remaining of the paper, we always use the principle cut for complex power functions without explicitly stating so. If we need a different cut, it will be made clear and explicit in the text.}. On any circular arc $\mathcal{C}_{R}$ in the complement of $\{0\}\cup[z_{1},\infty)\cup(-\infty,-z_{1}]$ centered at the origin with large radius $R$, the expression
\[
\left|\int_{\mathcal{C}_{R}}f(z,t)dz\right|
\]
is at most 
\[
\ointctrclockwise_{|z|=R}\frac{|Q(z)Q(-z)|^{1/2}\exp^{-nt}\left(\Arg(z_{1}-z)+\Arg(z_{2}-z)+\Arg(z_{1}+z)+\Arg(z_{2}+z)\right)}{R^{n+1}}d|z|.
\]
Using the estimates
\[
\left|\Arg(z-z_{1})+\Arg(z-z_{2})+\Arg(z+z_{1})+\Arg(z+z_{2})\right| \le2\pi,
\]
and
\[
|Q(z)Q(-z)|^{1/2}  =\mathcal{O}(R^{2})
\]
for $|z|=R$ we conclude that 
\[
\lim_{R\rightarrow\infty}\int_{\mathcal{C}_{R}}f(z,t)dz=0.
\]
As a consequence, 
\[
H_{n}\left(\frac{1}{2}+int\right)=\frac{n!}{2\pi i}\int_{\Gamma_{1}\cup\Gamma_{2}}f(z,t)dz,
\]
where $\Gamma_{1}$ and $\Gamma_{2}$ are two loops around two cuts
$(-\infty,-z_{1}]$ and $[z_{1},\infty)$ with counter clockwise
orientation. Using the substitution $z \mapsto -z$ we see that 
\[
\frac{1}{2\pi i}\int_{\Gamma_{1}}f(z,t)dz=\frac{(-1)^{n+1}}{2\pi i}\int_{\Gamma_{2}}f(z,-t)dz.
\]
On the other hand, the substitution $z \mapsto \overline{z}$ leads to the identity
\[
\frac{(-1)^{n+1}}{2\pi i}\int_{\Gamma_{2}}f(z,-t)dz=\frac{(-1)^{n+1}}{2\pi i}\overline{\int_{\Gamma_{2}}f(z,t)dz}.
\]
We deduce that $\pi H_{n}(1/2+int)$ is either the imaginary part, or $-i$ times
the real part of the integral
\begin{equation}
\int_{\Gamma_{2}}f(z,t)dz, \label{eq:intloopcut}
\end{equation}
depending on the parity of $n$.
\subsection{Approximating $\displaystyle{\int_{\Gamma_{2}}f(z,t)dz}$ - the saddle point method} \label{sec:approxsaddlepoint}
In order to approximate \eqref{eq:intloopcut} using the saddle point method (see for example \cite[Ch.\,4]{temme}),
we write 
\begin{equation} \label{eq:fdefn}
f(z,t)=e^{-n\phi(z,t)}\psi(z),
\end{equation}
where 
\begin{equation}
\phi(z,t)=\Log z-it\left(\Log(z_{1}-z)+\Log(z_{2}-z)-\Log(z_{1}+z)-\Log(z_{2}+z)\right),\label{eq:phidef}
\end{equation}
and 
\begin{equation}
\psi(z)=\frac{1}{z}\exp^{1/2}\left(\Log(z_{1}-z)+\Log(z_{2}-z)+\Log(z_{1}+z)+\Log(z_{2}+z)\right).\label{eq:psidef}
\end{equation}
As a function in $z$, the critical points of $\phi(z,t)$ are the
solutions of the equation 
\begin{align*}
0 & =\frac{1}{z}-it\left(\frac{1}{z-z_{1}}+\frac{1}{z-z_{2}}-\frac{1}{z+z_{1}}-\frac{1}{z+z_{2}}\right),
\end{align*}
which (after clearing denominators) is equivalent to 
\begin{equation}
z^{4}-2it(z_{1}+z_{2})z^{3}-(z_{1}^{2}+z_{2}^{2})z^{2}+2itz_{1}z_{2}(z_{1}+z_{2})z+z_{1}^{2}z_{2}^{2}=0.\label{eq:critpointseq}
\end{equation}
The discriminant in $z$ of \eqref{eq:critpointseq} is a polynomial
in $t$, whose positive zeros (in $t$) are 
\[
T_{1}:=\frac{z_{2}-z_{1}}{z_{1}+z_{2}}\qquad\text{and }\qquad T_{2}:=\frac{z_{1}+z_{2}}{4\sqrt{z_{1}z_{2}}}.
\]
We note that $T_{2}\ge T_{1}$ since 
\[
T_{2}^{2}-T_{1}^{2}=\frac{(z_{1}^{2}-6z_{1}z_{2}+z_{2}^{2})^{2}}{16z_{1}z_{2}(z_{1}+z_{2})^{2}} \geq 0.
\]
Set 
\begin{align} \label{eq:Tdefn}
T & =\begin{cases}
T_{1} & \text{if }z_{1}^{2}-6z_{1}z_{2}+z_{2}^{2}\ge0\\
T_{2} & \text{if }z_{1}^{2}-6z_{1}z_{2}+z_{2}^{2}<0.
\end{cases}
\end{align}
For each $t\in(0,T)$, the four solutions
of \eqref{eq:critpointseq} are $\zeta_{1}(t)$, $\zeta_{2}(t)$,
$-\overline{\zeta_{1}(t)}$, $-\overline{\zeta_{2}(t)}$ where 
\begin{align}
\zeta_{1}(t) & =\frac{z_{1}+z_{2}}{2}\left(it-\sqrt{T_{1}^{2}-t^{2}}+\sqrt{1-2t^{2}-2it\sqrt{T_{1}^{2}-t^{2}}}\right) \quad \text{for} \quad 0\le t<T_{1},\label{eq:zeta1defnT1}\\
\zeta_{1}(t) & =\frac{z_{1}+z_{2}}{2}\left(it+i\sqrt{t^{2}-T_{1}^{2}}+\sqrt{1-2t^{2}-2t\sqrt{t^{2}-T_{1}^{2}}}\right) \quad \text{for} \quad T_{1}\le t\le T_{2},\label{eq:zeta1defnT2}
\end{align}
and 
\begin{align}
\zeta_{2}(t) & =\frac{z_{1}+z_{2}}{2}\left(it+\sqrt{T_{1}^{2}-t^{2}}+\sqrt{1-2t^{2}+2it\sqrt{T_{1}^{2}-t^{2}}}\right) \quad \text{for} \quad 0\le t<T_{1},\label{eq:zeta2defnT1}\\
\zeta_{2}(t) & =\frac{z_{1}+z_{2}}{2}\left(it-i\sqrt{t^{2}-T_{1}^{2}}+\sqrt{1-2t^{2}+2t\sqrt{t^{2}-T_{1}^{2}}}\right) \quad \text{for} \quad T_{1}\le t\le T_{2}.\label{eq:zeta2defnT2}
\end{align}
We next establish some properties of the two curves $\zeta_{1}(t)$
and $\zeta_{2}(t)$ and their images under the map $\phi(z,t)$. These properties provide the justification for the the choice of the loop we use in applying the argument principle in Section \ref{sec:zerodist}.

\subsubsection{Properties of $\zeta_{1}(t)$ and $\zeta_{2}(t)$}

Using the definitions of $\zeta_{1}(t)$ and $\zeta_{2}(t)$ it is straightforward to
verify that in the case $T=T_{2}$ and $t\in[T_{1},T_{2})$,
the two curves $\zeta_{1}(t)$ and $\zeta_{2}(t)$ are parts of the
circle centered at the origin with radius $\sqrt{z_{1}z_{2}}$. 

\begin{figure}
\begin{centering}
\includegraphics[scale=0.3]{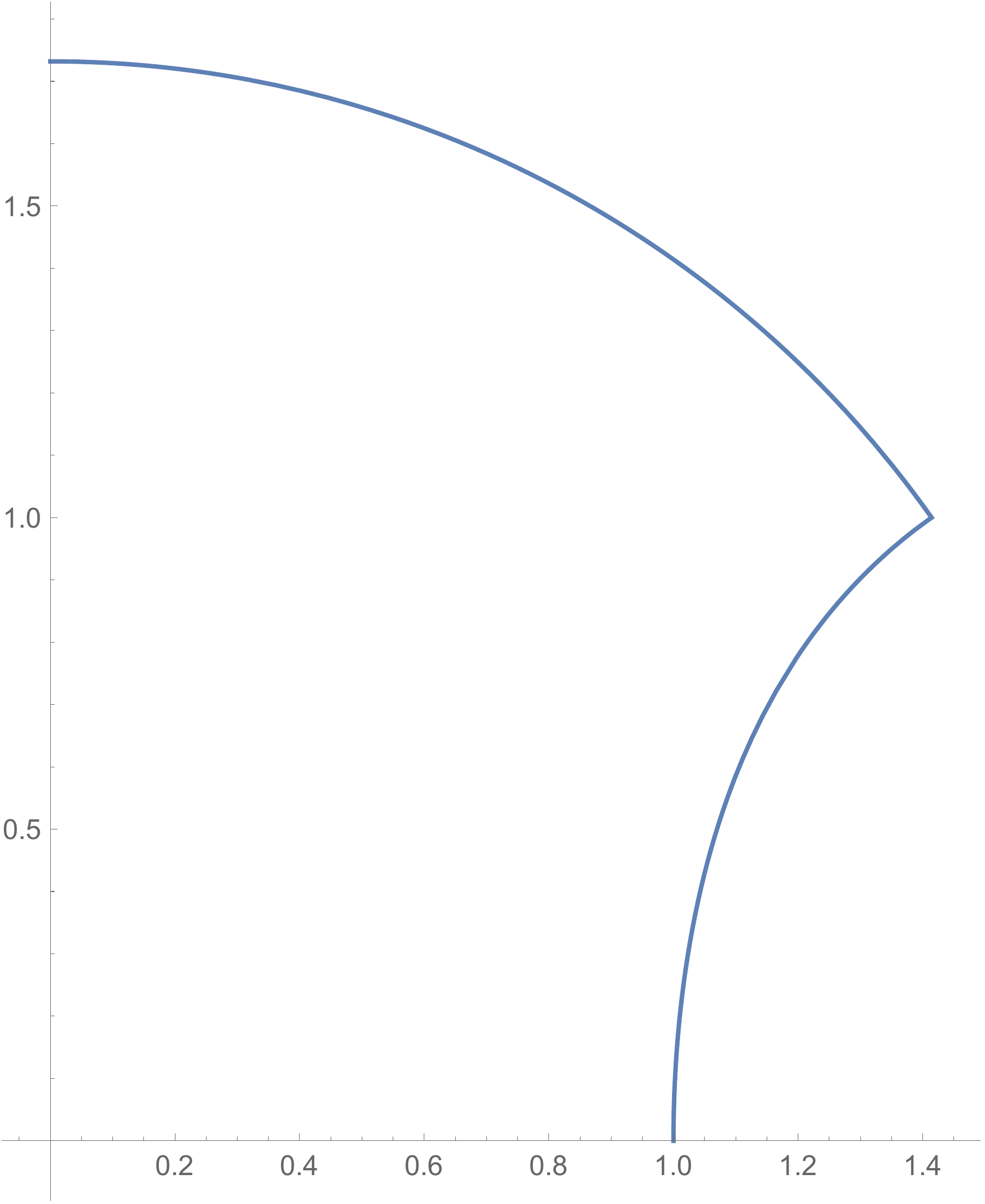} $\qquad$\includegraphics[scale=0.3]{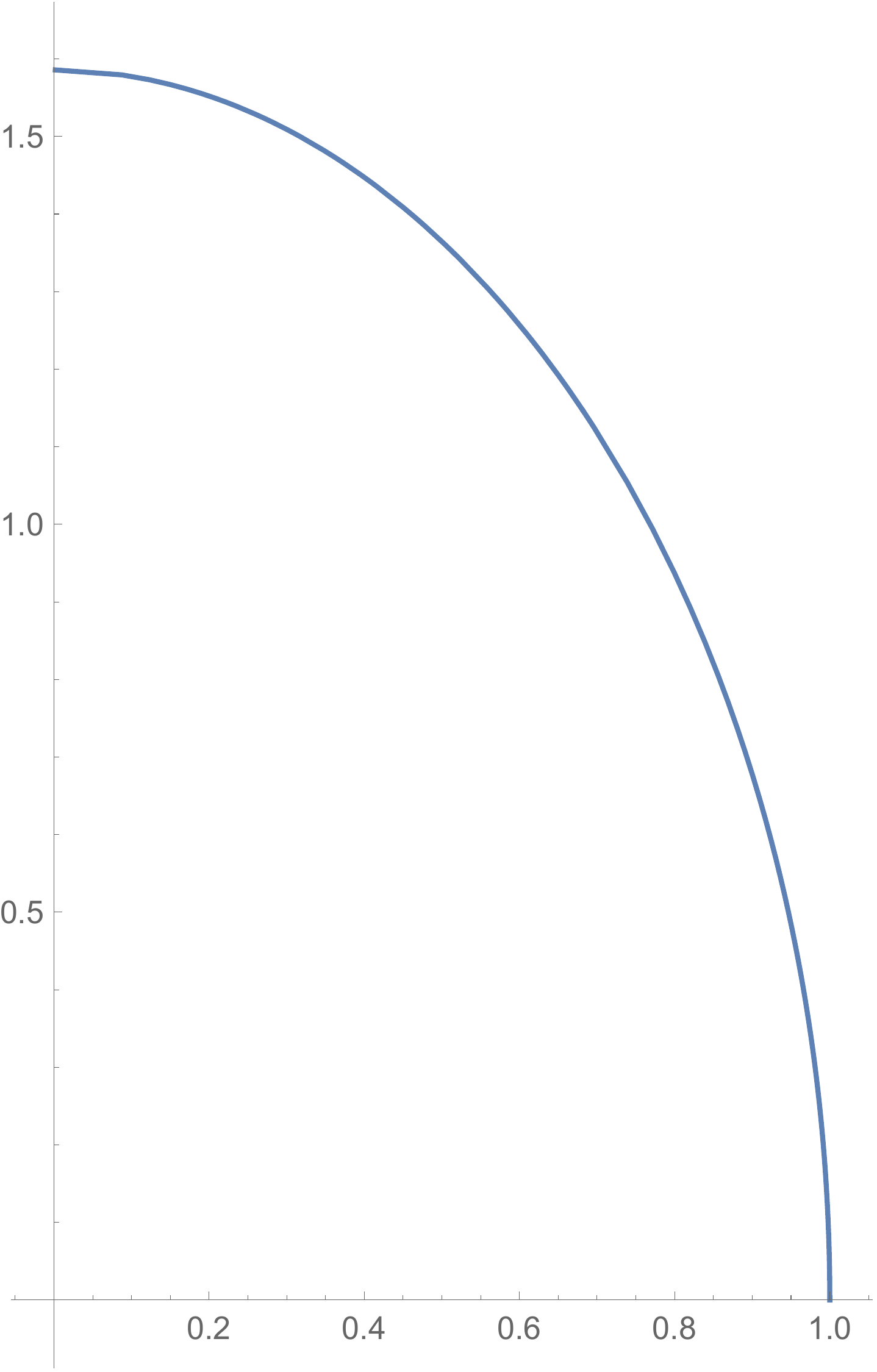}
\par\end{centering}
\caption{\label{fig:zetamap}The map $\zeta_{1}(t)$ on $[0,T]$ for $(z_{1},z_{2})=(1,3)$
(left) and $(1,7)$ (right)}
\end{figure}

\begin{lem}
\label{lem:zetalocation}If $\zeta_{k}(t)$, $1\le k\le2$, is defined
as in equations \eqref{eq:zeta1defnT1}-\eqref{eq:zeta2defnT2}, then
for $0<t<T$ these curves lie in the first open quadrant. 
\end{lem}

\begin{proof}
We need to study the solutions to equation \eqref{eq:critpointseq}.
For the purposes of this proof, we set 
\[
f(z)=z^{4}-2it(z_{1}+z_{2})z^{3}-(z_{1}^{2}+z_{2}^{2})z^{2}+2itz_{1}z_{2}(z_{1}+z_{2})z+z_{1}^{2}z_{2}^{2},
\]
\[
h(z)=4t^{2}z_{1}z_{2}z^{2},
\]
and 
\[
g(z)=f(z)-h(z)=(z^{2}-2itz_{1}z-z_{1}^{2})(z^{2}-2itz_{2}z-z_{2}^{2}).
\]
We note that the zeros of $h$ are at the origin, while the zeros
of $g$ are given by 
\[
r_{1,2,3,4}=iz_{k}(t\pm\sqrt{t^{2}-1}),\qquad k=1,2.
\]
By assumption, $0<t<T$, and the reader will recall that $0<z_{1}<z_{2}$.
If $T=T_{1}$, then trivially $0<t<1$, since $T_{1}<1$. On the other
hand, if $T=T_{2}$, then we are in the case $z_{1}^{2}-6z_{1}z_{2}+z_{2}^{2}<0$.
This implies that 
\[
z_{1}^{2}+2z_{1}z_{2}+z_{2}^{2}<8z_{1}z_{2}<16z_{1}z_{2},
\]
and consequently, $z_{1}+z_{2}<4\sqrt{z_{1}z_{2}}$. We conclude that
$T_{2}<1$ and hence $t<1$. This means that we may rewrite the roots
of $g$ as 
\[
r_{1,2,3,4}=\pm z_{k}\sqrt{1-t^{2}}+iz_{k}t,\qquad k=1,2.
\]
It is clear that the roots $r_{1,2,3,4}$ lie on the two circular
arcs $|z|=z_{1}$ and $|z|=z_{2}$. We wish to employ Rouch\'{e}'s theorem
to show that $g$, and $g+h=f$ have the same number of zeros inside
the curve $\Gamma$, which consists of the arcs $z=\epsilon e^{i\theta}$
and $z=Re^{i\theta}$ for $0\leq\theta\leq\pi/2$ and the two line
segments $\epsilon\leq z\leq R$, and $iz$, with $\epsilon\leq z\leq R$.
Here $0<\epsilon\ll1$, and $0<R-z_{2}\ll1$. Since $g$ has two
zeros there for any $0<t<T$, so does $f$, thereby establishing the claim. We now demonstrate that $|h|<|g|$ on $\Gamma$, which in turn will imply that $g$ and $f$ have the same number of zeros inside $\Gamma$.  \\
Note first that 
\begin{eqnarray}
|h(z)|^{2} & = & 16t^{4}z_{1}^{2}z_{2}^{2}|z|^{4},\quad\text{and}\nonumber \\
|g(z)|^{2} & = & 16t^{4}z_{1}^{2}z_{2}^{2}|z|^{4}\nonumber \\
 & + & 4t^{2}|z|^{2}\left(z_{1}^{2}(|z^{2}-z_{2}^{2}|^{2}-\Re(4itz_{2}z(z_{2}^{2}-\bar{z}^{2})))+z_{2}^{2}(|z^{2}-z_{1}^{2}|^{2}-\Re(4itz_{1}z(z_{1}^{2}-\bar{z}^{2})))\right)\nonumber \\
 & + & \left(|z^{2}-z_{1}^{2}|^{2}-\Re(4itz_{1}z(z_{1}^{2}-\bar{z}^{2}))\right)\left(|z^{2}-z_{2}^{2}|^{2}-\Re(4itz_{2}z(z_{2}^{2}-\bar{z}^{2})))\right) \label{eq:gmod}
\end{eqnarray}
Consider first the piece $\epsilon\leq z\leq R$ of $\Gamma$, which
lies on the real axis. In this case every expression in equation \eqref{eq:gmod}
that involves taking the real parts is zero. Consequently, 
\begin{eqnarray*}
|g(z)|^{2} & = & 16t^{4}z_{1}^{2}z_{2}^{2}|z|^{4}+4t^{2}|z|^{2}\left(z_{1}^{2}(|z^{2}-z_{2}^{2}|^{2})+z_{2}^{2}(|z^{2}-z_{1}^{2}|^{2})\right)+|z^{2}-z_{1}^{2}|^{2}|z^{2}-z_{2}^{2}|^{2}\\
 & > & 16t^{4}z_{1}^{2}2z_{2}^{2}|z|^{4}=|h(z)|^{2}.
\end{eqnarray*}
Next we consider the segment of $\Gamma$ lying on the imaginary axis:
$z=is$, with $\epsilon\leq s\leq R$. On this segment, 
\[
\Re(4itz_{k}z(z_{k}^{2}-\bar{z}^{2}))=4tz_{k}z(z_{k}^{2}+s^{2})),\qquad k=1,2.
\]
It follows that 
\begin{eqnarray*}
|g(z)|^{2} & = & 16t^{4}z_{1}^{2}z_{2}^{2}|z|^{4}\\
 & + & 4t^{2}|z|^{2}\left(z_{1}^{2}(|s^{2}+z_{2}^{2}|^{2}+4tz_{2}s(z_{2}^{2}+s^{2})))+z_{2}^{2}(|s^{2}+z_{1}^{2}|^{2}+4tz_{1}s(z_{1}^{2}+s^{2})))\right)\\
 & + & \left(|s^{2}+z_{1}^{2}|^{2}+4tz_{1}s(z_{1}^{2}+s^{2}))\right)\left(|s^{2}+z_{2}^{2}|^{2}+4tz_{2}s(z_{2}^{2}+s^{2})))\right)\\
 & > & 16t^{4}z_{1}^{2}z_{2}^{2}|z|^{4}=|h(z)|^{2}.
\end{eqnarray*}
We now turn our attention to the arc $z=\varepsilon e^{i\theta},0\leq\theta\leq\pi/2$.
On this segment of $\Gamma$, we have 
\begin{eqnarray*}
|z^{2}-z_{k}^{2}|^{2} & > & (z_{k}^{2}-\varepsilon^{2})^{2},\qquad k=1,2\\
\Re(4itz_{k}z(z_{k}^{2}-\bar{z}^{2})) & = & -4tz_{k}\varepsilon(\varepsilon^{2}(\cos\theta\sin(2\theta)-\sin\theta\cos(2\theta))+\sin\theta z_{2}^{2})\\
 & = & -4tz_{k}\varepsilon\sin\theta(\varepsilon^{2}+z_{2}^{2})\\
 & \leq & 0,\qquad k=1,2.
\end{eqnarray*}
Therefore, 
\begin{eqnarray*}
|g(z)|^{2} & > & 16t^{4}z_{1}^{2}z_{2}^{2}|z|^{4}\\
 & + & 4t^{2}\varepsilon^{2}\left(z_{1}^{2}((z_{2}^{2}-\varepsilon^{2})^{2})+z_{2}^{2}((z_{1}^{2}-\varepsilon^{2})^{2})\right)\\
 & + & (z_{1}^{2}-\varepsilon^{2}))^{2}(z_{2}^{2}-\varepsilon^{2})^{2}\\
 & > & 16t^{4}z_{1}^{2}z_{2}^{2}|z|^{4}=|h(z)|^{2}.
\end{eqnarray*}
Finally, on the arc $z=Re^{i\theta},0\leq\theta\leq\pi/2$, we have
\begin{eqnarray*}
|z^{2}-z_{k}^{2}|^{2} & > & (R-z_{k})^{2}(R^{2}+z_{k}^{2}),\qquad k=1,2\\
\Re(4itz_{k}z(z_{k}^{2}-\bar{z}^{2})) & = & -4tz_{k}R(R^{2}(\cos\theta\sin(2\theta)-\sin\theta\cos(2\theta))+\sin\theta z_{2}^{2})\\
 & = & -4tz_{k}R\sin\theta(R^{2}+z_{2}^{2})\\
 & \leq & 0,\qquad k=1,2.
\end{eqnarray*}
We conclude that 
\begin{eqnarray*}
|g(z)|^{2} & > & 16t^{4}z_{1}^{2}z_{2}^{2}|z|^{4}\\
 & + & 4t^{2}R^{2}\left(z_{1}^{2}((R-z_{2})^{2}(R^{2}+z_{2}^{2})+z_{2}^{2}(R-z_{1})^{2}(R^{2}+z_{1}^{2})\right)\\
 & + & (R-z_{1})^{2}(R^{2}+z_{1}^{2})(R-z_{2})^{2}(R^{2}+z_{2}^{2})\\
 & > & 16t^{4}z_{1}^{2}z_{2}^{2}|z|^{4}=|h(z)|^{2}
\end{eqnarray*}
on this piece as well. In summary, $|h(z)|<|g(z)|$ on $\Gamma$,
and by Rouch\'{e}'s Theorem, $f$ has two zeros in the first open
quadrant for all $0<t<T$ as desired. 
\end{proof}
\begin{lem}
\label{lem:moduluszeta} Let $T$ be as defined in \eqref{eq:Tdefn}. Then for any $t\in[0,T]$, 
\[
|\zeta_{1}(t)|\le\sqrt{z_{1}z_{2}}\le|\zeta_{2}(t)|
\]
with equality only when $T=T_{2}$ and $t\in[T_{1},T_{2}]$. 
\end{lem}

\begin{proof}
One can easily verify from equations \eqref{eq:zeta1defnT2} and \eqref{eq:zeta2defnT2}
that
\[
|\zeta_{1}(t)|=\sqrt{z_{1}z_{2}}=|\zeta_{2}(t)|
\]
 in case $T=T_{2}$ and $t\in[T_{1},T_{2}]$. Suppose now that $t\in(0,T_{1})$. Then 
 \[ 
 \Im \left(1-2t^2-2it\sqrt{T_1^2-t^2} \right)<0,
 \]
 and hence 
 \[
\Im \left(\sqrt{1-2t^2-2it\sqrt{T_1^2-t^2}} \right) <0
 \]
 as well. It follows from equations \eqref{eq:zeta1defnT1} and \eqref{eq:zeta2defnT1} that $|\Re(\zeta_1(t)|< |\Re(\zeta_2(t)|$ and $|\Im(\zeta_1(t)|< |\Im(\zeta_2(t)|$, and consequently, 
\[
|\zeta_{1}(t)|<|\zeta_{2}(t)|.
\]
Since the product of the four roots of the polynomial in equation \eqref{eq:critpointseq} is equal to its constant term, we see that 
\[
|\zeta_1(t)|^2|\zeta_2(t)|^2=z_1^2z_2^2.
\]
Taking square roots and applying the preceding inequality finishes the proof.
\end{proof}
\begin{lem}
\label{lem:asympzetaT} Let $T$ be as defined in \eqref{eq:Tdefn}. Then as $t\rightarrow T$, 

\[
\zeta_1(t)-\zeta_1(T)=\begin{cases}
\displaystyle{\frac{z_{1}+z_{2}}{2}\left(-\sqrt{T_{1}^{2}-t^{2}}+\frac{T_{1}\sqrt{T_{1}^{2}-t^{2}}}{\sqrt{2T_{1}^{2}-1}}+\mathcal{O}(T_{1}-t)\right)} & \text{if }T=T_{1}\\
\displaystyle{\frac{\sqrt{32}(z_{1}z_{2})^{3/4}}{\sqrt{(z_{1}+z_{2})(-z_{1}^{2}+6z_{1}z_{2}-z_{2}^{2})}}\sqrt{T_{2}-t}+\mathcal{O}(T_{2}-t)} & \text{if }T=T_{2}
\end{cases}.
\]
\end{lem}

\begin{proof}
Suppose first that $T=T_{1}$. We note that 
\[
1-2T_{1}^{2}=-\frac{z_{1}^{2}-6z_{1}z_{2}+z_{2}^{2}}{(z_{1}+z_{2})^{2}}<0,
\]
and conclude that as $t\rightarrow T_{1}$,
\begin{align*}
\sqrt{1-2t^{2}-2it\sqrt{T_{1}^{2}-t^{2}}} & =\sqrt{(1-2T_{1}^{2})\left(1-2iT_{1}\frac{\sqrt{T_{1}^{2}-t^{2}}}{1-2T_{1}^{2}}+\mathcal{O}(T_{1}-t)\right)}\\
 & =-i\sqrt{2T_{1}^{2}-1}\left(1-iT_{1}\frac{\sqrt{T_{1}^{2}-t^{2}}}{1-2T_{1}^{2}}+\mathcal{O}(T_{1}-t)\right).
\end{align*}
Thus equation \eqref{eq:zeta1defnT1} implies that as $t\rightarrow T_{1}$
\[
\zeta_1(t)-\zeta_1(T_{1})=\frac{z_{1}+z_{2}}{2}\left(-\sqrt{T_{1}^{2}-t^{2}}+\frac{T_{1}\sqrt{T_{1}^{2}-t^{2}}}{\sqrt{2T_{1}^{2}-1}}+\mathcal{O}(T_{1}-t)\right).
\]
In the case $T=T_{2}$, using equation \eqref{eq:zeta1defnT2} and
a CAS we obtain that as $t\rightarrow T_{2}$,
\[
\zeta_1(t)-\zeta_1(T_{2})=\frac{\sqrt{32}(z_{1}z_{2})^{3/4}}{\sqrt{(z_{1}+z_{2})(-z_{1}^{2}+6z_{1}z_{2}-z_{2}^{2})}}\sqrt{T_{2}-t}+\mathcal{O}(T_{2}-t).
\]
\end{proof}
\begin{lem}
\label{lem:asympzetaT1}Suppose $T=T_{2}$, and let 
\begin{equation}
d=-\frac{(z_{1}+z_{2})\sqrt{T_{1}}}{\sqrt{2}}\left(1+\frac{iT_{1}}{\sqrt{1-2T_{1}^{2}}}\right).\label{eq:dform}
\end{equation}
Then for $k=1,2$ the following hold:
\[
\zeta_{k}(t)-\zeta_{k}(T_{1})=\begin{cases}
(-1)^{k+1}d\sqrt{T_{1}-t}+\mathcal{O}(T_{1}-t) & \text{ if }t\rightarrow T_{1}^{-}\\
(-1)^k id\sqrt{t-T_{1}}+\mathcal{O}(t-T_{1}) & \text{ if }t\rightarrow T_{1}^{+}
\end{cases}.
\]
\end{lem}

\begin{proof}
We treat the case $t\rightarrow T_{1}^{-}$ and $k=2$. The
remaining cases follow from similar computations. Using
\[
1-2T_{1}^{2}=-\frac{z_{1}^{2}-6z_{1}z_{2}+z_{2}^{2}}{(z_{1}+z_{2})^{2}}>0
\]
we conclude as $t\rightarrow T_{1}^{-}$,
\begin{align*}
\sqrt{1-2t^{2}+2it\sqrt{T_{1}^{2}-t^{2}}} & =\sqrt{(1-2T_{1}^{2})\left(1+2iT_{1}\frac{\sqrt{T_{1}^{2}-t^{2}}}{1-2T_{1}^{2}}+\mathcal{O}(T_{1}-t)\right)}\\
 & =\sqrt{1-2T_{1}^{2}}\left(1+iT_{1}\frac{\sqrt{T_{1}^{2}-t^{2}}}{1-2T_{1}^{2}}+\mathcal{O}(T_{1}-t)\right).
\end{align*}
Thus \eqref{eq:zeta2defnT1} implies that as $t\rightarrow T_{1}^{-}$
\begin{align*}
\zeta_{2}(t)-\zeta_{2}(T_{1})= & \frac{z_{1}+z_{2}}{2}\left(\sqrt{T_{1}^{2}-t^{2}}+\frac{iT_{1}\sqrt{T_{1}^{2}-t^{2}}}{\sqrt{1-2T_{1}^{2}}}+\mathcal{O}(T_{1}-t)\right)\\
= & -d\sqrt{T_{1}-t}+\mathcal{O}(T_{1}-t).
\end{align*}
\end{proof}
\begin{rem}
\label{rem:dprop} Suppose that $\phi$ is as defined in equation \eqref{eq:phidef} and $d$ as in equation \eqref{eq:dform}. Using a CAS one can verify that 
\[
d^3 \cdot \left(\frac{\partial^3 \phi(z,t)}{\partial z^3} \Bigg|_{(\zeta_2(T_1),T_1)} \right)=d^3 \cdot \left(\frac{\partial^3 \phi(z,t)}{\partial z^3} \Bigg|_{(\zeta_1(T_1),T_1)} \right)=\frac{2\sqrt{2}(z_{1}+z_{2})}{z_{2}-z_{1}}\sqrt{\frac{z_{1}^{2}-z_{2}^{2}}{z_{1}^{2}-6z_{1}z_{2}+z_{2}^{2}}}\in\mathbb{R}^{+}
\]
and
\[
d \cdot \phi_{z,t}(\zeta_{2}(T_{1}),T_{1})=\frac{\sqrt{2}(z_{1}+z_{2})}{z_{2}-z_{1}}\sqrt{\frac{z_{1}^{2}-z_{2}^{2}}{z_{1}^{2}-6z_{1}z_{2}+z_{2}^{2}}}\in\mathbb{R^{+}}.
\]
The fact that these quantities are positive and real will play important roles as we develop asymptotic expressions for the integral in \eqref{eq:intloopcut} in the coming sections. 
\end{rem}
Our next result provides a key guiding component of the implementation of the saddle point method in our asymptotic approximation of the integral in \eqref{eq:intloopcut}.
\begin{prop}
\label{lem:dominantcrit} Let $T$ be as defined in \eqref{eq:Tdefn}. Then for any $t\in(0,T)$, 
\[
\Re\phi(\zeta_{1}(t),t)\le\Re\phi(\zeta_{2}(t),t),
\]
with equality only when $T=T_{2}$ and $T_1 \leq t < T_2$.
\end{prop} 

\begin{proof}
We first consider the case $t\in(0,T_{1})$ regardless of whether $T=T_{1}$
or $T=T_{2}$. Since $\zeta_{1}(t)$ and $\zeta_{2}(t)$ approach
$z_{1}$ and $z_{2}$ respectively as $t\rightarrow0$, the definition of $\phi(z,t)$ (c.f. equation \eqref{eq:phidef}) implies that
\[
\Log(z_1)=\lim_{t\rightarrow0}\Re\phi(\zeta_{1}(t),t)<\lim_{t\rightarrow0}\Re\phi(\zeta_2(t),t)=\Log(z_2).
\]
Suppose by way of contradiction that $\Re\phi(\zeta_{1}(t),t)\ge\Re\phi(\zeta_{2}(t),t)$
for some $t\in(0,T_{1})$. Then the equality
\[
\Re\phi(\zeta_{1}(T_{1}),T_{1})=\Re\phi(\zeta_{2}(T_{1}),T_{1})
\]
 along with the Mean Value Theorem implies the existence of a $t^* \in(0,T_{1})$ for which
\begin{equation} \label{eq:RePhider}
\frac{d}{dt}\left(\Re\phi(\zeta_{1}(t),t)-\Re\phi(\zeta_{2}(t),t)\right)\Big|_{t=t^*}=0.
\end{equation}
 We use the chain rule to compute 
 \[
 \frac{d\phi(\zeta_k(t),t)}{dt}=\frac{\partial \phi}{\partial z}\Big|_{(\zeta_k(t),t)}\frac{d\zeta_k(t)}{dt}+\frac{\partial \phi}{\partial t}=\frac{\partial \phi}{\partial t},
 \] since $\partial\phi/\partial z|_{(\zeta_{k}(t),t)}=0$ for $k=1,2$ by virtue of $\zeta_k(t)$ being a critical point of $\phi$. Using this relation we rewrite equation \eqref{eq:RePhider} to obtain
 \begin{align*}
0 & =\Re\left(\frac{\partial \phi(z,t)}{\partial t} \Big|_{(\zeta_1(t^*),t^*)}- \frac{\partial \phi(z,t)}{\partial t} \Big|_{(\zeta_2(t^*),t^*)}\right)\\
 & =\Im\left(F(\zeta_{1}(t^*))-F(\zeta_{2}(t^*))\right),
\end{align*}
where 
\begin{equation} \label{eq:Fdefn}
F(z)=\Log(z_1-z)+\Log(z_2-z)-\Log(z_1+z)-\Log(z_2+z).
\end{equation}
 It is straightforward to check that
\[
\Im\left(F(\zeta_1(T_1))-F(\zeta_2(T_1))\right)=0,
\]
hence applying the Mean Value theorem to the function $\Im\left(F(\zeta_1(t))-F(\zeta_2(t))\right)$  provides a $t^{**}\in(t^*,T_{1})$ so that
\begin{equation} \label{eq:ImdF}
\Im\left(\frac{d}{dt}F(\zeta_{1}(t^{**}))-\frac{d}{dt}F(\zeta_{2}(t^{**}))\right)=\Im\left(\frac{dF}{dz}\Big|_{\zeta_{1}(t^{**})}\zeta'_{1}(t^{**})-\frac{dF}{dz}\Big|_{\zeta_{2}(t^{**})}\zeta'_{2}(t^{**})\right)=0.
\end{equation}
Since $\phi(z,t)=\Log z-i t F(z)$, and $\displaystyle{\frac{\partial \phi(z,t)}{\partial z} \Big|_{ (\zeta_k(t),t)}=0}$ for all $t \in (0,T)$, we see that
\[
\frac{1}{\zeta_k(t)}-it\frac{dF}{dz} \Big|_{\zeta_k(t)}=0,
\] 
or equivalently, 
\[
\frac{dF}{dz}\Big|_{\zeta_{k}(t)}=\frac{1}{it\zeta_{k}(t)}.
\]
This, together with the second equation in \eqref{eq:ImdF}, implies that
\begin{equation} \label{eq:Relogderzeta}
0 =\Re\left(\frac{\zeta_{1}'(t^{**})}{\zeta_{1}(t^{**})}-\frac{\zeta_{2}'(t^{**})}{\zeta_2(t^{**})}\right).
\end{equation}
Using the definition of $\zeta_{1}(t)$ and $\zeta_{2}(t)$
in equations \eqref{eq:zeta1defnT1} and \eqref{eq:zeta2defnT1} we find the explicit experssions
\begin{align*}
\frac{\zeta_{1}'(t)}{\zeta_{1}(t)} & =\frac{t+i\sqrt{T_{1}^{2}-t^{2}}}{\sqrt{T_{1}^{2}-t^{2}}\sqrt{-2it\sqrt{T_{1}^{2}-t^{2}}-2t^{2}+1}}, \qquad \textrm{and}\\
\frac{\zeta_{2}'(t)}{\zeta_{2}(t)} & =\frac{i\left(\sqrt{T_{1}^{2}-t^{2}}+it\right)}{\sqrt{T_{1}^{2}-t^{2}}\sqrt{2it\sqrt{T_{1}^{2}-t^{2}}-2t^{2}+1}},
\end{align*}
from which we readily deduce that
\[
\frac{\zeta_{1}'(t)}{\zeta_{1}(t)}=-\overline{\frac{\zeta'_{2}(t)}{\zeta_2(t)}}.
\]
It follows that 
\[
\frac{\zeta_1'(t)}{\zeta_1(t)}-\frac{\zeta_2'(t)}{\zeta_2(t)}=\frac{\zeta_1'(t)}{\zeta_1(t)}+\overline{\frac{\zeta_1'(t)}{\zeta_1(t)}}=2 \Re \left(\frac{\zeta_1'(t)}{\zeta_1(t)} \right),
\]
and hence equation \eqref{eq:Relogderzeta} can be reformulated as 
\[
\Re\frac{\zeta_{1}'(t^{**})}{\zeta_{1}(t^{**})}=0.
\]
But then $\Im\left(\frac{\zeta_{1}'(t^{**})}{\zeta_{1}(t^{**})}\right)^{2}=0$, which implies that $1-T_1^2=0$, contradicting the fact that $0<T_1<1$. We point out here that in light of the arguments above, $\lim_{t\rightarrow T_{1}} \zeta_{1}'(t)/\zeta_{1}(t)>0$ implies that
\[
\Re\frac{\zeta_{1}'(t)}{\zeta_{1}(t)}>0,\qquad\forall t\in(0,T_{1}),
\]
a fact we shall use shortly (see the proof of Lemma \ref{lem:rephizetacurve}).\\
We now turn our attention to the second statement in the proposition.
If $T=T_{2}$, then $\forall t\in[T_{1},T_{2})$
\begin{equation} \label{eq:eqreparts}
\Re\frac{\zeta_{1}'(t)}{\zeta_{1}(t)}=\Re\frac{\zeta_{2}'(t)}{\zeta_2(t)}=0,
\end{equation}
since on this range of the parameter $t$ we have the explicit formulas
\begin{align*}
\frac{\zeta_{1}'(t)}{\zeta_{1}(t)} & =\frac{i\left(\sqrt{t^{2}-T_{1}^{2}}+t\right)}{\sqrt{t^{2}-T_{1}^{2}}\sqrt{1-2t^{2}-2t\sqrt{t^{2}-T_{1}^{2}}}}, \qquad \textrm{and}\\
\frac{\zeta_{2}'(t)}{\zeta_{2}(t)} & =\frac{i\left(\sqrt{t^{2}-T_{1}^{2}}-t\right)}{\sqrt{t^{2}-T_{1}^{2}}\sqrt{1-2t^{2}+2t\sqrt{t^{2}-T_{1}^{2}}}}.
\end{align*}
Using equation \eqref{eq:eqreparts} and reversing the direction in the argument leading to equation \eqref{eq:Relogderzeta} we conclude that $d/dt (\Im F(\zeta_1(t))-F(\zeta_2(t))=0$. That is, $\Im F(\zeta_{1}(t))=\Im F(\zeta_{2}(t))=C$ for some constant $C$, and all $t \in [T_1,T_2)$. Recalling (c.f. equations \eqref{eq:zeta1defnT2} and \eqref{eq:zeta2defnT2}) that
\[
\left|\zeta_{1}(t)\right|=|\zeta_{2}(t)|=\sqrt{z_{1}z_{2}}, \qquad t \in [T_1,T_2),
\]
and that $\Re \phi(\zeta_k(t),t)=\ln|\zeta_k(t)|+t\Im F(\zeta_k(t))$ (c.f. equation \eqref{eq:phidef}), the conclusion $\Re\phi(\zeta_{1}(t),t)=\Re\phi(\zeta_{2}(t),t)$ readily follows.
\end{proof}
\begin{rem}
\label{rem:rephicirc} Using the definition of $F$ in equation \eqref{eq:Fdefn} and a CAS we find that  for any $(z,t)$ in the domain of $\phi(z,t)$ with $z=\sqrt{z_{1}z_{2}}e^{i\theta}$, $0 \leq \theta \leq \pi/2$,
\[
\frac{d}{d\theta}\left(\Im F(z)\right)=0.
\]
Consequently, $\Re\phi(z,t)=\lim_{\theta\rightarrow0^{+}}\Re\phi(\sqrt{z_{1}z_{2}}e^{i\theta},t)=\ln(\sqrt{z_{1}z_{2}})-2\pi t$. 
\end{rem}

\begin{lem} \label{lem:imphiinc} Let $\phi$ be as defined in \eqref{eq:phidef}, $T$ be as defined in \eqref{eq:Tdefn}, and $\zeta_k$ be as defined in \eqref{eq:zeta1defnT1}-\eqref{eq:zeta2defnT2}. Then for $k=1,2$, 
\begin{itemize}
\item[(i)] the function $\Im\phi(\zeta_{k}(\cdot),\cdot): (0,T) \to (0, \pi/2)$ is strictly monotone increasing and onto, and
\item[(ii)] the function $\Re\phi(\zeta_{k}(\cdot),\cdot)$ is strictly decreasing on $(0,T)$. 
\end{itemize}
\end{lem}

\begin{proof}
(i) Since $\zeta_{k}(t)$ is a critical point of $\phi$, we have 
\begin{equation}
\frac{d\phi}{dt}=\frac{\partial\phi}{\partial z}\Big|_{(\zeta_k(t),t)}\frac{d\zeta_{k}(t)}{dt}+\frac{\partial\phi}{\partial t}=\frac{\partial\phi}{\partial t}.\label{eq:dphidt}
\end{equation}
We take the imaginary part of both sides to obtain 
\begin{equation}
\dfrac{d}{dt}\left(\Im\phi(\zeta_{k}(t),t)\right)=-\ln\frac{|z_{1}-\zeta_{k}(t)||z_{2}-\zeta_{k}(t)|}{|z_{1}+\zeta_{k}(t)||z_{2}+\zeta_{k}(t)|}, \label{eq:derivImphi}
\end{equation}
which is positive for $t\in(0,T)$ since $\zeta_{k}(t)$ lies in
the first open quadrant by Lemma \ref{lem:zetalocation}. Finally,
since $\zeta_{k}(t)$ approaches $z_{k}$ (resp. a purely imaginary number)
as $t\rightarrow 0$ (resp. $t\rightarrow T$), the definition of $\phi$ implies that
\begin{align*}
\lim_{t\rightarrow0}\Im\phi(\zeta_{k}(t),t) & =0,\\
\lim_{t\rightarrow T}\Im\phi(\zeta_{k}(t),t) & =\pi/2.
\end{align*}
For part (ii), we compute 
\[
\frac{d}{dt}\left(\Re\phi(\zeta_{k}(t),t)\right)=\Arg(z_{1}-\zeta_{k}(t))+\Arg(z_{2}-\zeta_{k}(t))-\Arg(z_{1}+\zeta_{k}(t))-\Arg(z_{1}+\zeta_{k}(t))<0,
\]
since $\zeta_{k}(t)$ lies in the first open quadrant. The proof is complete.
\end{proof}
\begin{lem}
\label{lem:rephizetacurve} Let $\phi$ be as defined in \eqref{eq:phidef}, $T$ be as defined in \eqref{eq:Tdefn}, and $\zeta_k$ be as defined in \eqref{eq:zeta1defnT1}-\eqref{eq:zeta2defnT2}. Then for any $t,t_{0}\in(0,T)$,
\begin{align*}
\Re\phi(\zeta_{1}(t_{0}),t_{0}) & \le\Re\phi(\zeta_{1}(t),t_{0}) \qquad \textrm{and}\\
\Re\phi(\zeta_{2}(t_{0}),t_{0}) & \ge\Re\phi(\zeta_{2}(t),t_{0}),
\end{align*}
with equality only when (i) $t=t_{0}$ or (ii) $T=T_{2}$
and $t,t_{0}\in[T_{1},T_{2})$.
\end{lem}

\begin{proof}
We begin with the first inequality. Let $F$ be as defined in equation \eqref{eq:Fdefn}. Then
\begin{align*}
\frac{d}{dt}\Re\phi(\zeta_{1}(t),t_{0}) & =\Re\frac{d}{dt}\phi(\zeta_{1}(t),t_{0})\\
 & =\Re\frac{\zeta'_{1}(t)}{\zeta(t)}+t_{0}\Im\frac{d}{dt}F(\zeta_{1}(t))\\
 & =\Re\frac{\zeta'_{1}(t)}{\zeta_{1}(t)}-\frac{t_{0}}{t}\Re\frac{\zeta'_{1}(t)}{\zeta_{1}(t)}\\
 & =\Re\frac{\zeta'_{1}(t)}{\zeta_{1}(t)}\left(1-\frac{t_0}{t} \right).
\end{align*}
If $t\in(0,T_{1})\backslash\{t_{0}\}$, the claim $\Re(\phi(\zeta_{1}(t_{0}),t_{0})<\Re\phi(\zeta_{1}(t),t_{0})$
follows from the fact that $\Re\left(\zeta'_{1}(t)/\zeta_{1}(t)\right)>0$.
If $T=T_{2}$ and $t\in[T_{1},T_{2})$, then by equation \eqref{eq:eqreparts} we have $\Re\left(\zeta'_{1}(t)/\zeta_{1}(t)\right)=0$, which implies
\[
\Re\phi(\zeta_{1}(t),t_{0})=\Re\phi(\zeta_{1}(t_{0}),t_{0})
\]
if $t_{0}\in[T_{1},T_{2})$. If on the other hand $t_{0}\in(0,T_{1})$, then
\[
\Re\phi(\zeta_{1}(t_{0}),t_{0})>\Re\phi(\zeta_{1}(T_{1}),T_{1})=\Re\phi(\zeta_{1}(t),t_{0}).
\]
 The second inequality in the lemma follows from analogous arguments, utilizing now that $\Re\left(\zeta_{2}'(t)/\zeta_{2}(t)\right)<0$ for $t\in(0,T_{1})$. 
\end{proof}
\begin{lem}
\label{lem:rephiimaxis} Let $\phi$ be as defined in \eqref{eq:phidef}, $T$ be as defined in \eqref{eq:Tdefn}. For any $t\in(0,T)$, the function $\Re\phi( \cdot ,t)$
is increasing on the positive imaginary axis, and decreasing on the negative imaginary axis. 
\end{lem}

\begin{proof}
We first show that $\frac{\partial}{\partial y}\Re\phi(iy,t)>0$ for
all $y>0$. By the chain rule, this is equivalent to showing that for all $y>0$,
\[
\Re\left(i\phi_{z}(iy,t)\right)=\frac{y^{4}-2t(z_{1}+z_{2})y^{3}+(z_{1}^{2}+z_{2}^{2})y^{2}-2tz_{1}z_{2}(z_{1}+z_{2})y+z_{1}^{2}z_{2}^{2}}{y(y^{2}+z_{1}^{2})(y^{2}+z_{2}^{2})}>0.
\]
Since equation \eqref{eq:critpointseq} has no purely imaginary solutions when
$t\in(0,T)$, $\Re(i\phi_{z}(iy,t))\ne0$ for any $y \in \mathbb{R}$. By evaluating the numerator of the fraction above at $y=0$, we conclude $\Re\phi(iy,t)$ is increasing
at every $y>0$. Applying the same reasoning mutatis mutandis, we obtain that $\Re\phi(iy,t)$ decreasing at every $y<0$. 
\end{proof}
Heuristically, for each $t\in(0,T)\backslash T_{1}$, the saddle
point method gives the (non-uniform in $t$) approximation 
\[
\int_{\Gamma_{2}}e^{-n\phi(z,t)}\psi(z,t)dz \sim \pm\frac{\sqrt{2\pi}\psi(\zeta_{k})e^{-n\phi(\zeta_{k},t)}}{\sqrt{n\phi_{z^{2}}(\zeta_{k},t)}}\qquad (k=1,2, \ n \to \infty).
\]
When estimating the integral, we must therefore consider the quantity
\[
\left| e^{-n \phi(\zeta_k(t),t)}\right|=e^{-n \Re \phi(\zeta_k(t),t)},
\]
which, together with Proposition \ref{lem:dominantcrit}, suggests that $\zeta_{1}(t)$ plays a more important role that $\zeta_{2}(t)$.Thus for the remainder of
the paper we denote $\zeta:=\zeta(t):=\zeta_{1}(t)$. The goal of
the ensuing sections is to provide rigorous arguments for
this approximation and to provide a condition under which the approximation
is uniform in $t$. 

\subsubsection{The main term of the approximation } \label{sec:maintermapprox}
The aim of this section, given $\zeta$, is to find a curve $\Gamma_I: I \to \mathbb{C}$ on some real interval containing the origin so that $z(0)=\zeta(t)$, and so that the integral over this curve becomes the dominant term in estimating the integral $\displaystyle{\int_{\Gamma_{2}}e^{-n\phi(z,t)}\psi(z,t)dz}$. We begin our quest by studying the behavior of $\phi$ near the curve $\zeta(t)$.\\
For each $t\in(0,T)$, the function $\phi(z,t)$ is analytic as a function
in $z$ on the open ball with center $\zeta$ and radius $\min(\Re\zeta,\Im\zeta)$,
with the power series representation 
\begin{equation}
\phi(z,t)=\phi(\zeta,t)+\frac{\phi_{z^{2}}(\zeta,t)}{2}(z-\zeta)^{2}+\sum_{k=3}^{\infty}\frac{\phi_{z^{k}}(\zeta,t)}{k!}(z-\zeta)^{k}.\label{eq:Taylorseriesphi}
\end{equation}
 Letting $t\rightarrow0$ in \eqref{eq:zeta1defnT1} we obtain
\begin{align}
z_1-\zeta & =z_1-\frac{z_1+z_2}{2}\left(1-T_{1}+it-itT_{1}+\mathcal{O}(t^2)\right)\nonumber \\
 & =-iz_1t+\mathcal{O}(t^{2})\label{eq:z1-zeta}.
\end{align}
This means that the curve $\zeta(t)$ approaches the real axis at an angle of $\pi/2$. Consequently, 
\[
z_1-\zeta \asymp\Im(\zeta-z_{1})=\Im\zeta.
\]
In addition, as $t \to 0$, we also have the relation 
\[
\phi_{z^{2}}(\zeta,t)\asymp\frac{t}{(\zeta-z_{1})^{2}}\asymp\frac{1}{t},
\]
and hence $t^2 \phi_{z^2}(\zeta,t) \asymp \Im \zeta$.
Since $\phi_{z^3}(\zeta(T),T) \neq 0$, expanding $\phi_{z^2}$ in a Taylor series and using Lemma \ref{lem:asympzetaT} we see that as $t\rightarrow T$,
\[
\phi_{z^{2}}(\zeta,t)\asymp\zeta-\zeta(T)\asymp\Re\zeta.
\]
Putting all this together we conclude that there exists small $\xi$ (independent of $t$) such that given any $t \in (0,T)$ the expansion in  \eqref{eq:Taylorseriesphi} is valid for all $z$ satsifying
\begin{align}
|z-\zeta| & <\xi|\phi_{z^{2}}(\zeta,t)|t^{2},\label{eq:convradius}
\end{align}
as the right side the above inequality is less than $\min(\Re\zeta,\Im\zeta)$. 
For $k \geq 3$, the definition of $\phi$ and our preceding discussion also yield the estimate
\begin{align} \label{eq:phikderest}
\phi_{z^{k}}(\zeta,t) & \ll\frac{(k-1)!}{|\zeta|^{k}}+\frac{t(k-1)!}{|\zeta-z_{1}|^{k}}+\frac{t(k-1)!}{|\zeta+z_{1}|^{k}}+\frac{t(k-1)!}{|\zeta-z_{2}|^{k}}+\frac{t(k-1)!}{|\zeta+z_{2}|^{k}}\\
 & \ll\frac{tk!}{|\zeta-z_{1}|^{k}}+\frac{k!}{\eta^{k}}\ll\frac{tk!A^{k}}{|\zeta-z_{1}|^{k}} \nonumber
\end{align}
for some small $\eta>0$ and large $A>0$ independent of $t$, $\zeta$
and $k$. Combining \eqref{eq:convradius} and \eqref{eq:phikderest} we conclude that for sufficiently small $\xi$, 
\begin{align}
\sum_{k=3}^{\infty}\frac{\phi_{z^{k}}(\zeta,t)}{k!}(z-\zeta)^{k} & =(z-\zeta)^{3}\sum_{k=3}^{\infty}\frac{\phi_{z^{k}}(\zeta,t)}{k!}(z-\zeta)^{k-3}\nonumber \\
 & \ll t\left|z-\zeta\right|^{3}\sum_{k=0}^{\infty}\frac{A^{3}}{|\zeta-z_{1}|^{3}}\left|\frac{A\xi|\phi_{z^{2}}(\zeta,t)|t^{2}}{\zeta-z_{1}}\right|^{k}\nonumber \\
 & \ll\frac{|z-\zeta|^{3}}{t^{2}}\sum_{k=0}^{\infty}\left|A\xi \phi_{z^{2}}(\zeta,t)t\right|^{k}\nonumber \\
 & \ll\frac{|z-\zeta|^{3}}{t^{2}}.\label{eq:Taylorphitail}
\end{align}
Using the above estimate for the tail of the series we write 
\begin{equation}
\phi(z,t)=\phi(\zeta,t)+\frac{\phi_{z^{2}}(\zeta,t)}{2}(z-\zeta)^{2}(1+h(z,t)),\label{eq:phizaroundzeta}
\end{equation}
where
\begin{equation}
h(z,t)=\mathcal{O}\left(\frac{z-\zeta}{t^{2}\phi_{z^{2}}(\zeta,t)}\right) \qquad (t \neq T_1).\label{eq:boundhzt}
\end{equation}
Consequently, if $\xi\ll1$, $t \in (0,T) \setminus \{T_1\}$ and $z$ satisfies
\eqref{eq:convradius}, then $|h(z,t)|<1/2$, and in turn 
\[
|\phi(z,t)-\phi(\zeta,t)|<\frac{3\xi^{2}}{4}|\phi_{z^{2}}(\zeta,t)|^{3}t^{4}.
\]
We now establish the existence of the desired curve. For $\xi \ll 1$ and for any fixed $y\in\mathbb{C}$ satisfying 
\begin{equation}
|y|<\frac{\xi}{2}\left|\phi_{z^{2}}(\zeta,t)\right|^{3/2}t^{2},\label{eq:yball}
\end{equation}
we apply Rouch\'{e}'s theorem to the functions\footnote{By virtue of $\xi$ being small enough, the square roots are all defined using the principle cut of the logarithm, and are analytic on the domain under consideration. Although it is possible that $\sqrt{\phi_{z^{2}}(\zeta,t)}$ is not continuous in $t$, we will later remove potential discontinutities by squaring this quantity.} 
\[
\frac{\sqrt{\phi_{z^{2}}(\zeta,t)}}{\sqrt{2}}(z-\zeta)\sqrt{1+h(z,t)}\qquad\textrm{and}\qquad\frac{\sqrt{\phi_{z^{2}}(\zeta,t)}}{\sqrt{2}}(z-\zeta)\sqrt{1+h(z,t)}-y
\]
to demonstrate that the equation 
\[
y=\frac{\sqrt{\phi_{z^{2}}(\zeta,t)}}{\sqrt{2}}(z-\zeta)\sqrt{1+h(z,t)}
\]
has exactly one solution in $z$ satisfying \eqref{eq:convradius}.
Indeed, the equation 
\[
0=\frac{\sqrt{\phi_{z^{2}}(\zeta,t)}}{\sqrt{2}}(z-\zeta)\sqrt{1+h(z,t)}
\]
has exactly one solution $z=\zeta$ satisfying \eqref{eq:convradius}.
In addition, if 
\[
|z-\zeta|=\xi|\phi_{z^{2}}(\zeta,t)|t^{2},
\]
then the inequality $\sqrt{1+h(z,t)}>1/\sqrt{2}$ implies that
\[
\frac{\sqrt{|\phi_{z^{2}}(\zeta,t)|}}{\sqrt{2}}|z-\zeta|\sqrt{1+h(z,t)}\ge\xi\frac{|\phi_{z^{2}}(\zeta,t)|^{3/2}}{2}t^{2}>|y|.
\]
It follows then that the equation 
\[
y=\frac{\sqrt{\phi_{z^{2}}(\zeta,t)}}{\sqrt{2}}(z-\zeta)\sqrt{1+h(z,t)}
\]
also has exactly one solution in $z$ satisfying \eqref{eq:convradius}. Thus, for any $y$ satisfying \eqref{eq:yball}, we may invert the
relation 
\begin{equation}
y=\frac{\sqrt{\phi_{z^{2}}(\zeta,t)}}{\sqrt{2}}(z-\zeta)\sqrt{1+h(z,t)}\label{eq:yfuncz}
\end{equation}
and - after using an analytic continuation argument - obtain a function
$z=z(y)$ which is analytic in the entire ball defined by \eqref{eq:yball}. 
\newline We continue by developing an asymptotic expression for the integral on a section of the curve parametrized by $z(y)$ under suitable conditions. In essence, for a given $\zeta$, we wish to understand how the integrand $e^{-n\phi(z,t)}\psi(z)$ compares to the 'central' value $e^{-n\phi(\zeta,t)}\psi(\zeta)$ for $z$ in the ball defined by \eqref{eq:convradius}. Consider the smooth
curve $\Gamma_{\epsilon}(y)$ parameterized by $z(y)$, $-\epsilon\le y\le\epsilon$, for any small $\epsilon=\epsilon(n)$ and $t \in K \subset (0,T)$ for which
\begin{equation}
\frac{\epsilon}{\sqrt{|\phi_{z^{2}}(\zeta,t)|^{3}}t^{2}}=o(1)\label{eq:epsiloncond}
\end{equation}
uniformly on $K$ as $n\rightarrow\infty$\footnote{This requirement implies for example that $t \neq T_1$. In addition, the condition may fail to hold for certain ranges of $t$ near $0, T_1$ and $T_2$, which is the principal reason for us having to handle these cases separately in Sections \ref{sec:smallT-t}, \ref{sec:SmallT1-t} and \ref{sec:smallt}.}.  If $z\in\Gamma_{\epsilon}(y)$, then 
\[
\phi(z,t)-\phi(\zeta,t)=y^{2},
\]
and using \eqref{eq:boundhzt} and \eqref{eq:yfuncz} we get
\begin{align}
z & =\zeta+\frac{\sqrt{2}y}{\sqrt{\phi_{z^{2}}(\zeta,t)}}\left(1+\mathcal{O}\left(\frac{z-\zeta}{t^{2}\phi_{z^{2}}(\zeta,t)}\right)\right)\nonumber \\
 & =\zeta+\frac{\sqrt{2}y}{\sqrt{\phi_{z^{2}}(\zeta,t)}}\left(1+\mathcal{O}\left(\frac{\epsilon}{t^{2}\sqrt{\phi_{z^{2}}(\zeta,t)^{3}}}\right)\right).\label{eq:critcontour}
\end{align}
Rearranging \eqref{eq:critcontour} and invoking condition \eqref{eq:epsiloncond} shows that given any $\xi \ll 1$, if $n \gg 1$, then
\begin{eqnarray*}
|z-\zeta|&=&\frac{\sqrt{2}|y|}{\sqrt{|\phi_{z^2}(\zeta,t)|}}\left(1+\mathcal{O}\left(\frac{\epsilon}{t^{2}\sqrt{|\phi_{z^{2}}(\zeta,t)|^{3}}}\right)\right)\\
&\leq&\frac{\sqrt{2}\epsilon}{\sqrt{\phi_{z^2}(\zeta,t)^3}t^2}|\phi_{z^2}(\zeta,t)| t^2\left(1+\mathcal{O}\left(\frac{\epsilon}{t^{2}\sqrt{\phi_{z^{2}}(\zeta,t)^{3}}}\right)\right)\\
&\leq& \xi |\phi_{z^2}(\zeta,t)| t^2,
\end{eqnarray*}
and hence the curve $\Gamma_{\epsilon}(y)$ lies in the first open quadrant. On the ball defined by \eqref{eq:convradius}, we write $\psi(z)=e^{\Psi(z)}$,
where 
\[
\Psi(z)=-\Log z+\frac{1}{2}\Log(z_{1}-z)+\frac{1}{2}\Log(z_{2}-z)+\frac{1}{2}\Log(z_{1}+z)+\frac{1}{2}\Log(z_{2}+z).
\]
We apply arguments similar to those leading to \eqref{eq:Taylorphitail} along with equation \eqref{eq:critcontour} and the fact that $\sqrt{|\Psi_{z}(\zeta)|}|\zeta-z_{1}|\asymp1$ for small $t$, to obtain
\begin{align*}
\Psi(z)-\Psi(\zeta) & \ll\frac{|z-\zeta|}{|\zeta-z_{1}|}\\
 & \ll\frac{\epsilon}{\sqrt{|\phi_{z^{2}}(\zeta,t)|}|\zeta-z_{1}|}=o(1),
\end{align*}
where the last equality relies on the facts $\phi_{z^{2}}(\zeta,t)\asymp1/t$ and $\zeta-z_{1}\asymp t$ for small $t$, as well as condition \eqref{eq:epsiloncond}. 
In summary, for $z \in \Gamma_{\epsilon}$ and $t \in K \subset (0,T)$ for which condition \eqref{eq:epsiloncond} holds,
\begin{eqnarray*}
e^{-n\phi(z,t)}&=&e^{-n(\phi(\zeta,t)+y^2)} \\
\psi(z)&=&e^{\Psi(z)}=e^{\Psi(\zeta)+(\Psi(z)-\psi(\zeta))}=\psi(\zeta)(1+o(1)), \quad \textrm{and}\\
dz&=& \frac{dz}{dy}dy=\frac{\sqrt{2}}{\sqrt{\phi_{z^2}(\zeta,t)}}(1+o(1))dy.
\end{eqnarray*} 
Consequently,
\begin{align*}
\int_{\Gamma_{\epsilon}}e^{-n\phi(z,t)}\psi(z)dz & =\frac{\sqrt{2}\psi(\zeta)e^{-n\phi(\zeta,t)}}{\sqrt{\phi_{z^{2}}(\zeta,t)}}\int_{-\epsilon}^{\epsilon}e^{-ny^{2}}\left(1+o(1)\right)dy\\
 & =\frac{\sqrt{2}\psi(\zeta)e^{-n\phi(\zeta,t)}}{\sqrt{n\phi_{z^2}(\zeta,t)}}\int_{-\epsilon\sqrt{n}}^{\epsilon\sqrt{n}}e^{-y^{2}}(1+o(1))dy.
\end{align*}
The reader will note that for any $\epsilon,n>0$, 
\[
\sqrt{\pi}=\int_{\mathbb{R}}e^{-y^{2}}dy=2\int_{\epsilon\sqrt{n}}^{\infty}e^{-y^{2}}+\int_{-\epsilon\sqrt{n}}^{\epsilon\sqrt{n}}e^{-y^{2}}dy,
\]
and that
\[
\int_{\epsilon\sqrt{n}}^{\infty}e^{-y^{2}}dy\le\int_{\epsilon\sqrt{n}}^{\infty}\frac{y}{\epsilon\sqrt{n}}e^{-y^{2}}dy=\frac{e^{-\epsilon^{2}n}}{2\epsilon\sqrt{n}}.
\]
As a result, we obtain the following asymptotic expression for the integral over $\Gamma_{\epsilon}$: 
\begin{equation} \label{eq:Gammaepsasymptotic}
\int_{\Gamma_{\epsilon}}e^{-n\phi(z,t)}\psi(z)dz=\frac{\sqrt{2\pi}\psi(\zeta)e^{-n\phi(\zeta,t)}}{\sqrt{n\phi_{z^{2}}(\zeta,t)}}\left(1+\mathcal{O}\left(\frac{e^{-\epsilon^{2}n}}{\epsilon\sqrt{n}}\right)+o(1)\right)
\end{equation}
for $\epsilon$ and $t \in K \subset (0,T)$ satisfying \eqref{eq:epsiloncond}. 

In the next section we extend append the curve $\Gamma_{\epsilon}$ with two tails going to $\infty$ in the upper and lower half planes  respectively, for $t\ne T_{1}$ (i.e. $\phi_{z^{2}}(\zeta,t)\ne0$).  We will also demonstrate that the integrals over these two tails are dominated by the integral in \eqref{eq:Gammaepsasymptotic}.  We shall employ the convention that a complex number approaches $\infty$ in the upper half (resp. lower half) plane if its modulus approaches $\infty$ and its argument lies in $(0,\pi)$
(resp. $(-\pi,0)$). 

\subsubsection{The tails of the approximation} \label{sec:tailapprox}

To ensure the integrals over the two tails are dominated by the integral in \eqref{eq:Gammaepsasymptotic}, we will choose these tails so that $\Re\phi(z,t)>\Re\phi(\zeta(\Gamma(\pm\epsilon),t)=\Re\phi(\zeta(t),t)+\epsilon^{2}$
for all $z$ in the tails. This section is dedicated to showing that such tails exist, and to demonstrating the claimed dominance. 
\begin{lem}
\label{lem:rephiposcomp}  Let $\phi$ be as defined in \eqref{eq:phidef}, $T$ be as defined in \eqref{eq:Tdefn}. Let $r\in\mathbb{R}$ and $t\in(0,T)$. If
$\mathcal{\mathcal{R}}_{r}$ denotes the intersection of a connected component of $\text{\ensuremath{\left\{ z:\Re\phi(z,t)>r\right\} }}$
with the first open quadrant, then $\partial\mathcal{R}_r \cap (0,\infty)$ contains a non-empty open interval. 
\end{lem}

\begin{proof}
The claim is trivial if $\mathcal{R}_{r}$ is unbounded since $\left\{ z:\Re\phi(z,t)>r\right\} $
contains all $z$ with large modulus. Assume next that $\mathcal{R}_{r}$ is bounded, and
note that $\Re\phi(z,t)$ is continuous at every $z\in\overline{\mathcal{R}_{r}}$,
except perhaps at $z_{1}$ and $z_{2}$, should they lie on $\partial\mathcal{R}_{r}$.
Thus, it suffices to show that there exists a real $z_{0}\in\partial\mathcal{R}_{r}\backslash\{z_{1},z_{2}\}$ such that $\Re\phi(z_{0},t)>r$. Suppose, by way of contradiction, that there is no such $z_{0}$. By Lemma \ref{lem:rephiimaxis}, the intersection of $\partial\mathcal{R}_{r}$ and the $y$-axis is either a point
where $\Re\phi(z,t)=r$ or it is empty. We deduce that $\Re\phi(z,t)=r$ for all $z\in\partial\mathcal{R}_r$, except at $z_{1}$ or $z_{2}$ were they to lie on $\partial\mathcal{R}_{r}$.
For $ \delta >0 $ set 
\[
\mathcal{R}_{r,\delta}=\left\{ z\in\mathcal{R}_{r}:|z-z_{1}|>\delta\text{ and }|z-z_{2}|>\delta\right\},
\]
and select $z^{*}\in\mathcal{R}_{r}$. Then $z^{*}\in\mathcal{R}_{r,\delta}$
for sufficiently small $\delta$, and since $\Re\phi(z^{*},t)>r$ and $\Im\phi(z,t)\rightarrow-\infty$
as $z\rightarrow z_{1}$ or $z\rightarrow z_{2}$, wee see that the
map $\phi(\cdot,t)-\phi(z^{*},t):\partial\mathcal{R}_{r,\delta} \to \mathbb{C}$ maps $\partial\mathcal{R}_{r,\delta}$
into the complement of the closed first quadrant. Since $\phi(z,t)$
is analytic on a region containing $\overline{\mathcal{R}_{r,\delta}}$, we may apply the argument principle theorem to conclude that
\[
0=i\Delta \arg_{z\in \partial\mathcal{R}_{r,\delta}}  \phi(z,t)-\phi(z^{*},t) =\int_{\partial\mathcal{R}_{r,\delta}}\frac{\phi'(z,t)}{\phi(z,t)-\phi(z^{*},t)}dz.
\]
On the other hand,
\[
\frac{1}{2\pi i}\int_{\partial\mathcal{R}_{r,\delta}}\frac{\phi'(z,t)}{\phi(z,t)-\phi(z^{*},t)}dz\geq 1,
\]
since $z^{*}\in\mathcal{R}_{r,\delta}$ is a zero of $\phi(z,t)-\phi(z^{*},t)$ inside the curve $\partial\mathcal{R}_{r,\delta}$, and we have reached a contradiction. The result follows.
\end{proof}
\begin{rem}
\label{rem:wholeint}For each fixed $t\in(0,T)$, as a function in
$x$, $\Re\phi(x+i\delta,t)$ converges uniformly to the function
\[
u(x)=\begin{cases}
\ln x & \text{ if }0<x<z_{1}\\
\ln x-t\pi & \text{ if }z_{1}<x<z_{2}\\
\ln x-2t\pi & \text{ if }z_{2}<x
\end{cases}
\]
on any real compact subset of $(0,\infty)\backslash\{z_{1},z_{2}\}$
as $\delta\rightarrow0^+$. Since $u(x)$ is increasing,
if $(a,b)\subseteq\partial\mathcal{R}_{r}\cap(0,z_{1})$, then so
is $(a,z_{1})$ and the similar conclusions hold for the intervals $(z_{1},z_{2})$
and $(z_{2},\infty)$. 
\end{rem}

\begin{defn}
\label{def:Ckdef} Let $T$ be as defined in \eqref{eq:Tdefn}. For each $t\in(0,T)$, we set $\mathcal{C}_{1}:=\mathcal{C}_1(t)$ (resp. 
$\mathcal{C}_{2}$) to be the intersection of the first open quadrant with the connected component of $\{z:\Re\phi(z,t)-\Re\phi(\zeta,t)>0\}$ containing $\Gamma_{\epsilon}(y)$ for $-\epsilon\le y<0$ (resp. $0<y\le\epsilon$).
\end{defn}
Before we present our next result, we recall two standard theorems from complex analysis.
\begin{lem}\cite[Lemma 1.2, p.511]{palka} \label{lem:palka1}. Let $\gamma:[a,b] \to \mathbb{C}$ be a simple, closed path, and let $c$ be a point of $(a,b)$ at which $\gamma$ is differentiable with $\dot \gamma(c) \neq 0$. There exists an $\epsilon >0$ for which it is true that the sets $I^+_{\epsilon}=\{\gamma(c)+s i \dot \gamma(c): 0< s \leq \epsilon\}$ and $I^-_{\epsilon}=\{\gamma(c)+s i \dot \gamma(c): -\epsilon \leq s <0\}$ lie in different components of $\mathbb{C}\setminus | \gamma|$.
\end{lem}
 In the following theorem $n(\gamma,z)$ denotes the winding number of a simple closed curve $\gamma$ about the point $z \in \mathbb{C}$.

\begin{thm}\cite[Theorem 1.3, p.553]{palka} \label{thm:palka2} Let $\gamma:[a,b] \to \mathbb{C}$ be a simple, closed, piecewise smooth path and let $D$ be the bounded component of $\mathbb{C} \setminus |\gamma |$. Then either $n(\gamma,z)=1$ for every $z \in D$ or $n(\gamma,z)=-1$ for all such $z$.
\end{thm}
\begin{lem}
\label{lem:distinctcomp} $T$ be as defined in \eqref{eq:Tdefn}, and let $\mathcal{C}_1$
and $\mathcal{C}_2$ be as in Definition \ref{def:Ckdef}. Then $\mathcal{C}_1 \neq \mathcal{C}_2$ 
\end{lem}

\begin{proof}
If $\mathcal{C}_{1}=\mathcal{C}_{2}$, then in this component we can
extend $\Gamma_{\epsilon}(y)$ to a simple, closed, piecewise smooth curve $\gamma$ on which $\Re\phi(z,t)\geq \Re\phi(\zeta,t)$ with equality only when $z=\zeta$. This implies that if $z_0 \in \mathbb{C}$ is any point for which $\Re \phi(\zeta,t)>\Re \phi(z_0,t)$, then the winding number of $\phi(\gamma,t)-\phi(z_0,t)$ about the origin is zero, and consequently, $n(\gamma,z_0)=0$. 
For small $y >0$, consider the points 
\[
z_0^{\pm}=\zeta\pm i\frac{\sqrt{2}y}{\sqrt{\phi_{z^{2}}(\zeta,t)}}.
\]
Using the expression in \eqref{eq:phizaroundzeta} we find that
\[
\phi(z_0^{\pm},t)-\phi(\zeta,t)=\frac{\phi_{z^2}(\zeta,t)}{2}\left(\pm i\frac{\sqrt{2}y}{\sqrt{\phi_{z^{2}}(\zeta,t)}}\right)^2(1+h(z_0^{\pm},t))=-y^2+\mathcal{O}(y^3).
\]
It follows that $\Re \phi(z_0^{\pm},t)-\Re\phi(\zeta,t) <0$, and hence $z_0^{\pm} \notin \gamma$, and $n(\gamma,z_0^{\pm})=0$. By Lemma \ref{lem:palka1}, the points $z_0^{\pm}$ lie in different components of the complement of the trace of $\gamma$ in the first open quadrant, and by \cite[Theorem 1.3, page 553]{palka} these components are both unbounded. We have reached contradiction since the complement of the trace of $\gamma$ has only one unbounded component.
\end{proof}
\begin{rem}
Using similar arguments we also establish that there are two distinct connected components
of the set $\left\{ z|\Re \phi(z,t)>\Re \phi(\zeta_{2},t)\right\} $ whose boundaries
contain $\zeta_{2}$. 
\end{rem}
The next result shows that the endpoints of the curve $\Gamma_{\epsilon}$ lie in regions of the plane from which it is possible to continue the curve to the point at infinity both in the upper and in the lower half planes,  while maintaining the desired relation $\Re \phi(z,t)-\Re \phi(\zeta,t) > 0$ for all $z$ in the tails.

\begin{figure}
\begin{centering}
\includegraphics[scale=0.3]{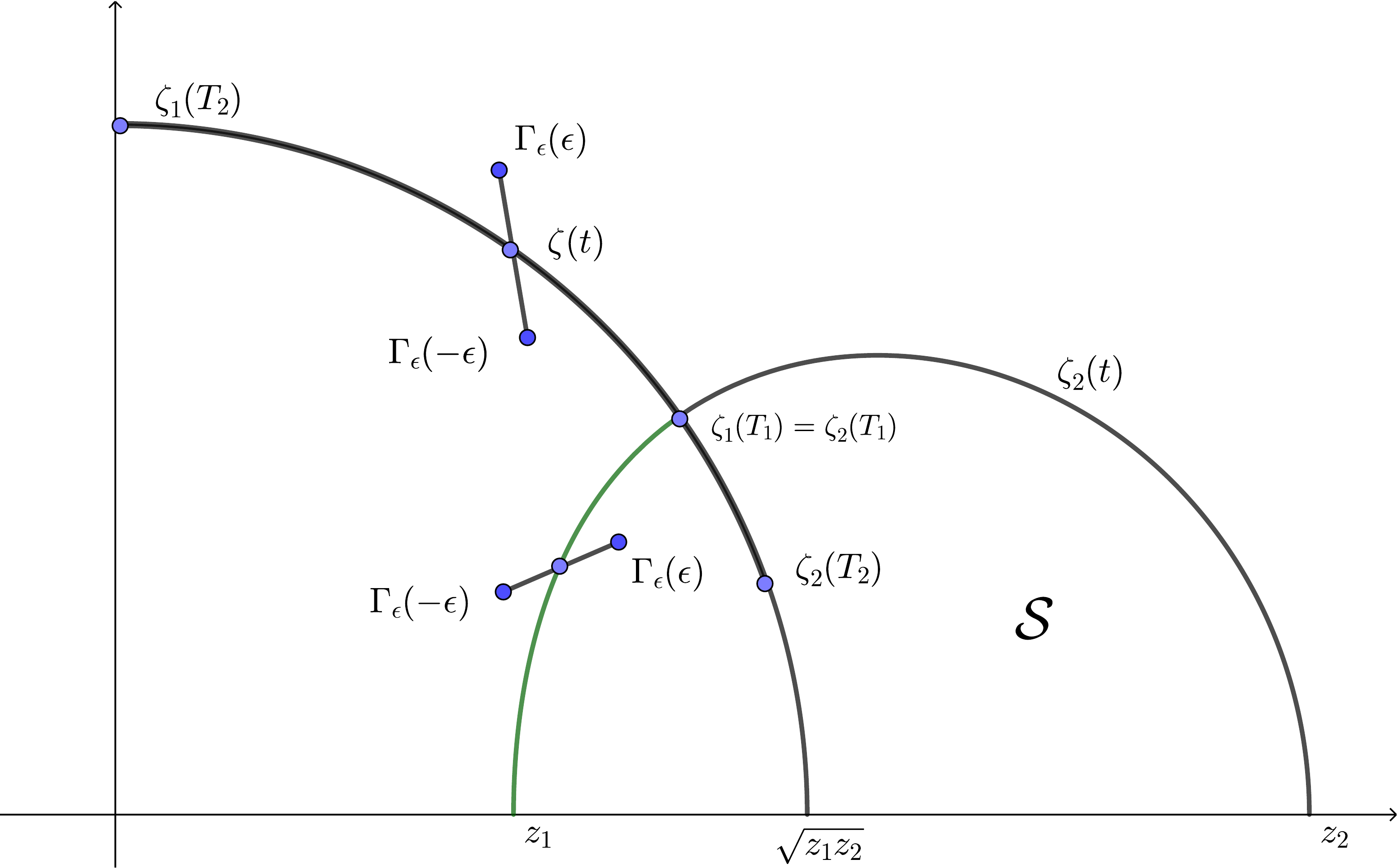}
\par\end{centering}
\caption{Possible locations of the $\Gamma_{\epsilon}$ curve}
\label{fig:intcontour}
\end{figure}

\begin{lem}
\label{lem:gammae-int}Suppose $T=T_{2}$ and \eqref{eq:epsiloncond}.
Let $\mathcal{S}$ be the region enclosed by $\zeta_{2}(\tau)$, $0\le\tau\le T_{1}$,
the circular arc radius $\sqrt{z_{1}z_{2}}$, and the segment from
$\sqrt{z_{1}z_{2}}$ to $z_{2}$ (see Figure \ref{fig:intcontour}). If $\epsilon >0$ is such that condition \eqref{eq:epsiloncond} holds, then $\Gamma_{\epsilon}(\pm\epsilon)$ lies outside $\overline{\mathcal{S}}$.
\end{lem}

\begin{proof}
Suppose first that $t\in[T_{1},T_{2})$. Then for all $y\in[-\epsilon,\epsilon]$ and for all $\tau\in(0,T)$,
\begin{align*}
\Re\phi(z(y),t) & >\Re\phi(\zeta_{1}(t),t) \qquad \textrm{(Def. of $\Gamma_{\epsilon}$)}\\
 & =\Re\phi(\zeta_{2}(t),t) \qquad \textrm{(Remark \ref{rem:rephicirc})}\\
 & \ge\Re\phi(\zeta_{2}(\tau),t).\qquad \textrm{(Lemma \ref{lem:rephizetacurve})}
\end{align*}
Were $\Gamma_{\epsilon}(\pm \epsilon) \in \overline{\mathcal{S}}$, the continuity of $\Gamma_{\epsilon}$ would necessitate its crossing of either the $\zeta_2$ curve, or the circular arc.\footnote{Recall that $\Gamma_{\epsilon}$ lies entirely in the first open quadrant, so it certainly doesn't cross the real axis.} Either of these cases would furbish a $y \in [-\epsilon,\epsilon]$ for which the above inequality would fail. \\
Suppose now that $t\in(0,T_{1})$. We claim that for any $z=\sqrt{z_{1}z_{2}}e^{i\theta}$, $0<\theta<\pi/2$, $\Re\phi(z,t)>\Re\phi(\zeta(t),t)+\epsilon^{2}$
 from which the result will follow, because $\Re \phi(\zeta,t)+\epsilon^2 >\Re \phi(z(y),t)$ for all $y \in [-\epsilon,\epsilon]$ and $t \in (0,T_1)$. In establishing the claim, it suffices to consider
$T_{1}-t=o(1)$ due to the fact that $\epsilon=o(1)$. For any $z$ on the circular arc (i.e. of the form $z=\sqrt{z_{1}z_{2}}e^{i\theta}$, $0<\theta<\pi/2$),
\begin{align*}
\Re\phi(z,t)-\Re\phi(\zeta(t),t) & =\Re\left(\phi(\zeta(T_{1}),t)-\phi(\zeta(t),t)\right).
\end{align*}
Using the Taylor series expansion of $\phi(z,t)$ as a function of two complex variables centered at $(\zeta(T_{1}),T_{1})$, and Lemma \ref{lem:asympzetaT1} (along with the constant $d$ defined therein (c.f. equation \eqref{eq:dform})) yield
\begin{align*}
 & \phi(\zeta(t),t)-\phi(\zeta(T_{1}),t)\\
= & d \cdot \phi_{z,t}(\zeta(T_{1}),T_{1})(t-T_{1})\sqrt{T_{1}-t} +\frac{1}{3!}d^{3}\cdot \phi_{z^{3}}(\zeta(T_{1}),T_{1})(T_{1}-t)^{3/2}+\mathcal{O}(T_{1}-t)^{2}.
\end{align*}
Taking real parts and invoking Remark \ref{rem:dprop} gives
\[
\Re\left(\phi(\zeta(T_{1}),t)-\phi(\zeta(t),t)\right)\asymp(T_{1}-t)^{3/2}.
\]
The inequality $\Re\left(\phi(\zeta(T_{1}),t)-\phi(\zeta(t),t)\right)>\epsilon^{2}$
follows from condition \eqref{eq:epsiloncond} and the asymptotic equivalence
\[
\phi_{z^{2}}(\zeta(t),t)\asymp(T_{1}-t)^{1/2},
\]
and the proof is complete.
\end{proof}
We now complete the argument demonstrating the existence of the two tails of $\Gamma_{\epsilon}$ needed to complete the asymptotic estimate for $\int_{\Gamma_2} f(z,t)dz$, provided that condition \eqref{eq:epsiloncond} holds. To this end, let $\mathcal{R}_{1}\subset\mathcal{C}_{1}$ (resp. 
$\mathcal{R}_{2}\subset\mathcal{C}_{2}$) be a connected component
of 
\[
\left\{ z:\Re\phi(z,t)>\Re\phi(\Gamma_{\epsilon}(\epsilon),t)=\Re\phi(\Gamma_{\epsilon}(-\epsilon),t)=\Re(\phi(\zeta,t))+\epsilon^{2}\right\} 
\]
whose boundary contains $\Gamma_{\epsilon}(-\epsilon)$ (resp. $\Gamma_{\epsilon}(\epsilon)$).
We argue that one of the two sets $\mathcal{R}_1$ and $\mathcal{R}_2$ is unbounded, while  the other contains an interval of the form $(z_1-\delta,z_1)$ in its closure for some $\delta>0$. Since for any $x_{0}\in(0,z_1)$, the function $\Re\phi(x_0-iy,t)$
is increasing in $y\in(0,\infty)$, the existence of the tail will follow because (i) in the unbounded component we can find a path to $\infty$ from one endpoint of $\Gamma_{\epsilon}$, and (ii) in the bounded component we can connect the other endpoint of $\Gamma_{\epsilon}(-\epsilon)$ to $x _0 \in (z_1,-\delta,z_0)$ and then append a ray $x_0-iy$, $y \geq 0$. 

Let $t \in (0,T)$, and suppose that condition \eqref{eq:epsiloncond} holds. Assume first that $\Re\phi(\zeta_{1}(t),t)+\epsilon^{2}\ge\Re\phi(\zeta_{2}(t),t)$.
Since $\epsilon=o(1)$, Proposition \ref{lem:dominantcrit} implies that $T=T_2$, $T_1< t<T_2$\footnote{we remind the reader that by convention, condition \eqref{eq:epsiloncond} implies that $t \neq T_1$} and $\Re \phi(\zeta_1(t),t)=\Re \phi(\zeta_2(t),t)$. By Lemma \ref{lem:rephizetacurve} we have
$\Re\phi(\zeta_{2}(t),t)\ge\Re\phi(\zeta_{2}(\tau),t)$ $\forall\tau\in(0,T)$, which means that if $z=\zeta_2(t)$ for any $t \in (0,T)$, then $z \notin \overline{\mathcal{R}}_k$, $k=1,2$. In addition, if $z \in \partial S$ with $|z|=\sqrt{z_1z_2}$, then by Remark \ref{rem:rephicirc}, $\Re \phi(z,t)=\Re \phi(\zeta(t),t)$ and once more we conclude that $z \notin \overline{\mathcal{R}}_k$, $k=1,2$. Thus, either $ \overline{\mathcal{R}}_k \subset S \cup [z_1,z_2]$, or $\overline{\mathcal{R}}_k \cap \overline{S} =\emptyset$. Since Lemma \ref{lem:gammae-int} shows that $\Gamma_{\epsilon}(\pm\epsilon)\notin\overline{\mathcal{S}}$, we conclude that  $\overline{\mathcal{R}}_k\cap\overline{\mathcal{S}}=\emptyset$ for $k=1,2$. Lemma \ref{lem:rephiposcomp} now establishes the claim, since it implies that one of the $\mathcal{R}_k$ intersects the positive real axis in a interval $(A,\infty)$ (and is hence unbounded), while the other intersects the positive real axis in an interval of the form $(z_1-\delta,z_1)$ for some $\delta >0$.

Next, we consider the case when $t \in (0,T)$ is such that $\Re\phi(\zeta_{1}(t),t)+\epsilon^{2}<\Re\phi(\zeta_2(t),t)$. In this case we must have $t\in(0,T_{1})$ (regardless off whether  $T=T_{1}$ or
$T=T_{2}$). Thus $\zeta_{2}(t)$ lies outside the closed ball
centered at the origin with radius $\sqrt{z_{1}z_{2}}$, and for all $z$ in the boundary of this ball in the first quadrant we have
\[
\Re\phi(\zeta_{2}(t),t)>\Re\phi(z,t).
\]
Consider now the two distinct, connected components of the set $\left\{ z|\Re\phi(z,t)>\Re\phi(\zeta_{2},t)\right\}$, whose boundaries contain $\zeta_2(t)$. By the argument above, we see that neither of these two components contain any points inside the ball with radius $\sqrt{z_1z_2}$. By Lemma \ref{lem:rephiposcomp} one will intersect the positive real axis in an interval of the form $(z_2-\nu,z_2)$, and the other in an interval of the form $(B, \infty)$ for some $B >z_2$. Since the connected components $\mathcal{R}_k$, $k=1,2$ are distinct,  we see that one of these has to have a common point with either $(z_2-\nu,z_2)$ for some $\nu>0$, or with $(B, \infty)$. The inequality $\Re\phi(\zeta_{1}(t),t)+\epsilon^{2}<\Re\phi(\zeta_2(t),t)$ now implies that both components of $\left\{ z|\Re\phi(z,t)>\Re\phi(\zeta_{2},t)\right\}$ with $\zeta_2(t)$ in their boundary are contained in $\mathcal{R}_{k}$ for $k=1$ or $k=2$ (w.o.l.g. we may assume that $k=2$). It follows now that $\partial\mathcal{R}_{2}$
contains two real intervals whose right endpoints are $z_{2}$ and
$\infty$, and $\partial\mathcal{R}_{1}$ contains an interval whose
right endpoint is $z_{1}$ (recall that $\mathcal{R}_{1}$ and $\mathcal{R}_{2}$
are distinct).

We summarize these discussions in the following proposition. 
\begin{prop} \label{lem:existenceGamma}
Let $\phi$ be as defined in \eqref{eq:phidef}, $T$ be as defined in \eqref{eq:Tdefn}, and $f$ be as defined in Section \ref{sec:gammadeform}. Suppose $t \in (0,T)$ and $\epsilon>0$ are such that condition \eqref{eq:epsiloncond} holds. Then there exists a piecewise smooth curve
$\Gamma$  parameterized by $z(y)$, $-\infty<y<\infty$ such that 
\begin{itemize}
\item[(i)] $\int_{\Gamma}f(z,t)dz=\pm\int_{\Gamma_{2}}f(z,t)dz$,
\item[(ii)] $\lim_{y\rightarrow\pm\infty}z(y)=\infty$ ,
\item[(iii)] $\phi(z(y),t)-\phi(\zeta,t)=y^2$ for all $-\epsilon<y<\epsilon$, 
\item[(iv)] $\Re\phi(z(y),t)\ge\Re\phi(\zeta,t)+\epsilon^{2}$ for all $|y| \geq \epsilon$, and with equality if and only if $y=\pm\epsilon$. 
\end{itemize}
\end{prop}

\begin{figure}[h]
\begin{centering}
\includegraphics[scale=0.3]{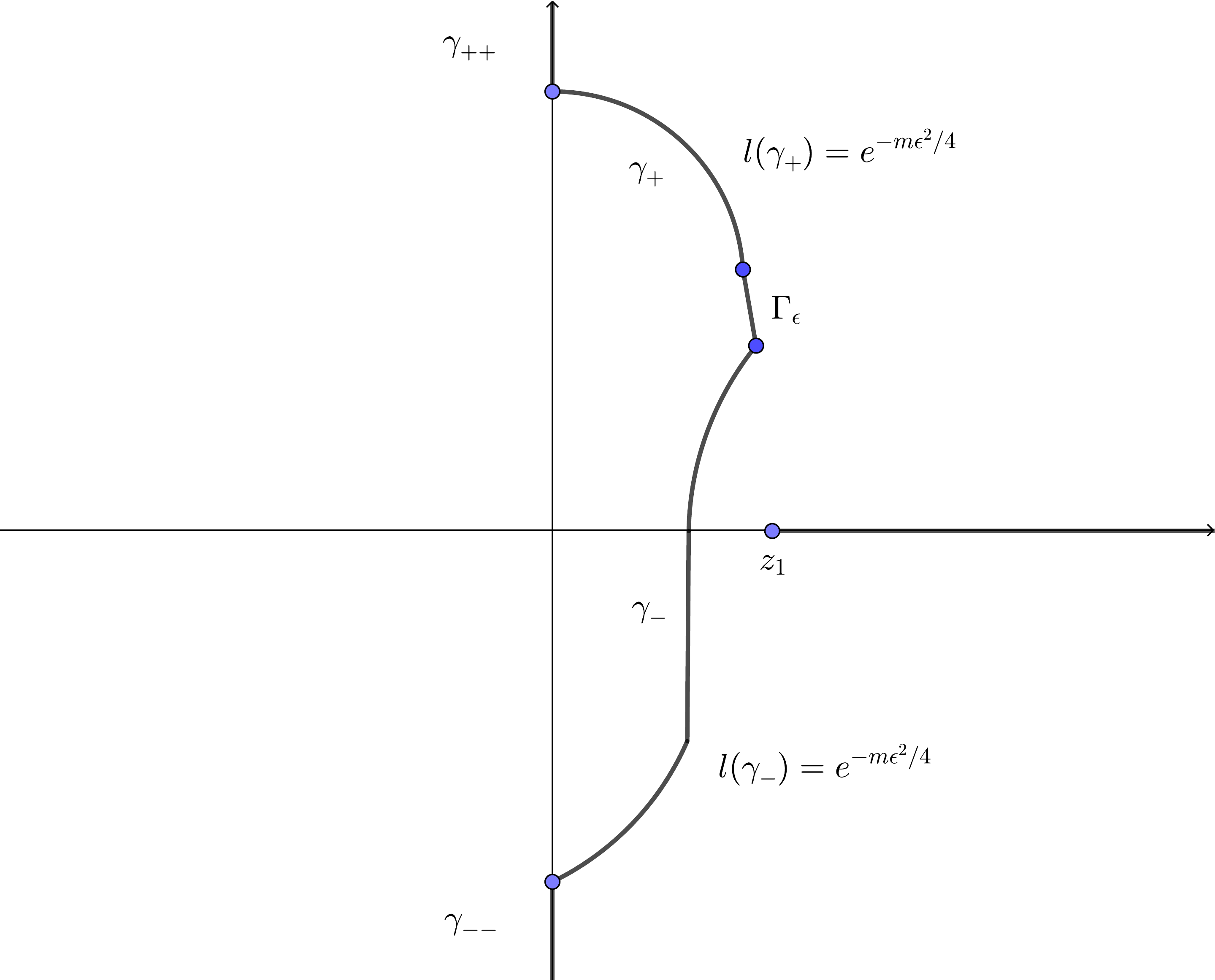}
\par\end{centering}
\caption{The curve $\Gamma$}
\label{fig:Gammacurve}
\end{figure}

\begin{rem}
\label{rem:largez}Since $\Re\phi(z,t)$ is large for large $z$,
we may assume that for large $z\in \Gamma$, the curve $\Gamma$ is given
by $z=iy$ and still have (iv) in Proposition \ref{lem:existenceGamma} hold. Figure \ref{fig:Gammacurve} illustrates the proposed path of integration.
\end{rem}
We break the path $\Gamma \setminus \Gamma_{\epsilon}$ into the pieces $\gamma_+, \gamma_-, \gamma_{++}$ and $\gamma_{--}$ as indicated in Figure \ref{fig:Gammacurve}, and give bounds on the integral of $e^{-n\phi(z,t)}\psi(z)$ over these segments under the assumptions that
\begin{align*} \label{cond:asymptbound}
(\dag) \qquad \qquad  & (i) \quad  \textrm{Condition \eqref{eq:epsiloncond} holds, and}\\
& (ii) \quad n \epsilon^2 \to \infty \quad \textrm{as} \quad  n \to \infty. \nonumber
\end{align*}
When combined, these estimates provide an asymptotic bound for the integral
\begin{equation}
\int_{\Gamma\backslash\Gamma_{\epsilon}}e^{-n\phi(z,t)}\psi(z)dz.\label{eq:tailint}
\end{equation}
Recall the definition of $\psi$ (c.f. equation \eqref{eq:psidef}) which implies that for $z\in\Gamma\backslash\Gamma_{\epsilon}$,
$\psi(z)=\mathcal{O}(z)$. On the portions $\gamma_+$ and $\gamma_-$ of $\Gamma\backslash\Gamma_{\epsilon}$ starting at the points $\Gamma_{\epsilon}(\pm \epsilon)$ with length $l(\gamma_+)=l(\gamma_-)=e^{n\epsilon^{2}/4}$ (which is large if $n\gg 1$ by assumption (ii) above), we have $|z|=\mathcal{O}(e^{n\epsilon^{2}/4})$, and consequently 
\[
\int_{\gamma_{+/-}}e^{-n\phi(z,t)}\psi(z)dz=\mathcal{O}\left(e^{-n\phi(\zeta,t)-n\epsilon^{2}}\int_{\gamma_{+/-}}|z||dz|\right)=\mathcal{O}\left(e^{-n\phi(\zeta,t)-n\epsilon^{2}/2}\right).
\]
On the other hand, by Remark \ref{rem:largez}, on the segments $\gamma_{++}$ and $\gamma_{--}$ we have
\[
\int_{\gamma_{++/--}}e^{-n\phi(z,t)}\psi(z)dz=\mathcal{O}\left(\int_{Y}^{\infty}\frac{dy}{y^{n-1}}\right)=\mathcal{O}\left(\frac{1}{(n-2)Y^{n-2}}\right),
\]
where $Y\gg1$. Since $\epsilon=o(1)$ and $\phi(\zeta,t)$ is bounded on $(0,T)$, we see that
\[
\frac{1}{Y} \asymp e^{-\phi(\zeta,t)-\epsilon^2/2},
\]
and consequently,
\[
\int_{\gamma_{++/--}}e^{-n\phi(z,t)}\psi(z)dz=\mathcal{O}\left(e^{-n\phi(\zeta,t)-n\epsilon^2/2}\right).
\]
We conclude that under the assumptions in $(\dag)$, 
\[
\int_{\Gamma_{2}}f(z,t)dz=\pm\frac{\sqrt{2\pi}\psi(\zeta)e^{-n\phi(\zeta,t)}}{\sqrt{n\phi_{z^{2}}(\zeta,t)}}\left(1+\mathcal{O}\left(\frac{e^{-\epsilon^{2}n}}{\epsilon\sqrt{n}}+\frac{e^{-\epsilon^{2}n/2}\sqrt{n\phi_{z^{2}}(\zeta,t)}}{\psi(\zeta)}\right)\right),
\]
which implies that
\begin{equation}
\left(\int_{\Gamma_{2}}f(z,t)dz\right)^{2}=\frac{2\pi\psi^{2}(\zeta)e^{-2n\phi(\zeta,t)}}{n\phi_{z^{2}}(\zeta,t)}\left(1+\mathcal{O}\left(\frac{e^{-\epsilon^{2}n}}{\epsilon\sqrt{n}}+\frac{e^{-\epsilon^{2}n/2}\sqrt{n\phi_{z^{2}}(\zeta,t)}}{\psi(\zeta)}\right)\right).\label{eq:squareasymp}
\end{equation}

We close this section by noting that we may not be able to satisfy the assumptions in $(\dag)$ when $\sqrt{|\phi_{z^{2}}(\zeta,t)|^3}t^{2}$ approaches $0$ too rapidly. Since $\sqrt{|\phi_{z^{2}}(\zeta,t)|^3}t^2$ is small when $t$ is
close to $0$, $T_{1}$, or $T_{2}$ (when $T=T_{2}$), we need separate arguments to develop asymptotic expressions for $\int_{\Gamma_2} f(z,t)dz$ for $t$ in these ranges.

\subsubsection{\label{sec:smallT-t} The asymptotics when $1/n^{2/3}\ll T-t\le\ln^{2}n/n^{2/3}$ }

We continue our work by looking at the case when $t$ is close to $T$ (as defined in equation \eqref{eq:Tdefn}). We begin with developing an asymptotic expression for $\phi(z,t)-\phi(\zeta,t)$ for $z$ close to $\zeta$. Using Lemma \ref{lem:asympzetaT} we conclude that for large $n$, 
\[
\zeta(t)-\zeta(T)=c\sqrt{T-t}+\mathcal{O}(T-t)\le2|c|\ln n/n^{1/3},
\]
where 
\begin{equation}
c=\begin{cases}
\frac{z_{1}+z_{2}}{2}\left(-\sqrt{2T_{1}}+\frac{T_{1}\sqrt{2T_{1}}}{\sqrt{2T_{1}^{2}-1}}\right) & \text{ if }T=T_{1}\\
\frac{\sqrt{32}(z_{1}z_{2})^{3/4}}{\sqrt{(z_{1}+z_{2})(-z_{1}^{2}+6z_{1}z_{2}-z_{2}^{2})}} & \text{ if }T=T_{2}
\end{cases}.\label{eq:cform}
\end{equation}
For $z$ in a small neighborhood of $\zeta$, we expand $\phi(z,t)$ (as a function of $z$) about $\zeta$:
\[
\phi(z,t)=\phi(\zeta,t)+\frac{\phi_{z^{2}}(\zeta,t)}{2!}(z-\zeta)^{2}+\frac{\phi_{z^{3}}(\zeta,t)}{3!}(z-\zeta)^{3}+\mathcal{O}\left((z-\zeta)^{4}\right).
\]
Since $\phi_{z^2}(\zeta(T),T)=0$, expanding $\phi_{z^2}(\zeta,t)$ in a bi-variate series centered at $(\zeta(T),T)$ yields
\begin{align*}
\phi_{z^{2}}(\zeta,t) & =\phi_{z^{3}}(\zeta(T),T)\left(\zeta-\zeta(T)\right)+\mathcal{O}(T-t)\\
 & =c\phi_{z^{3}}(\zeta(T),T)\sqrt{T-t}+\mathcal{O}(T-t).
\end{align*}
Simliarly, epxanding $\phi_{z^3}(\zeta,t)$ in a series centered at $(\zeta(T),T)$ gives
\begin{align*}
\phi_{z^{3}}(\zeta,t) & =\phi_{z^{3}}(\zeta(T),T)(1+\mathcal{O}(|\zeta-\zeta(T)|+|T-t|)\\
 & =\phi_{z^{3}}(\zeta(T),T)\left(1+\mathcal{O}\left(\sqrt{T-t}\right)\right).
\end{align*}
Thus 
\begin{align}
\phi(z,t)-\phi(\zeta,t) & =\frac{c^{3}\phi_{z^{3}}(\zeta(T),T)}{6}(z/c-\zeta/c)^{2}\left(3\sqrt{T-t}-(z/c-\zeta/c)\right)\nonumber \\
 & +\mathcal{O}\left((z-\zeta)^{2}(T-t)+(z-\zeta)^{3}\sqrt{T-t}+(z-\zeta)^{4}\right).\label{eq:phiapprox}
\end{align}
We also remark that using the definition of $\phi(\zeta,t)$ and a CAS, one easily verifies that if $T=T_{2}$, then
\[
c^{3}\phi_{z^{3}}(\zeta(T),T)=\frac{128i\sqrt{2}(z_{1}z_{2})^{3/4}}{(z_{1}+z_{2})^{3}\sqrt{(z_{1}+z_{2})\left(-z_{1}^{2}+6z_{1}z_{2}-z_{2}^{2}\right)}}\in i\mathbb{R}^{+},
\]
and if $T=T_{1}$, then
\[
c^{3}\phi_{z^{3}}(\zeta(T),T)=\frac{2i\sqrt{2}(z_{1}+z_{2})^{2}}{\sqrt{-z_{1}^{4}+6z_{1}^{3}z_{2}-6z_{1}z_{2}^{3}+z_{2}^{4}}}\in i\mathbb{R}^{+}.
\]

\begin{prop}
\label{lem:zycurve}Let $T$ be as defined in equation \eqref{eq:Tdefn}, $c$ be as defined in equation \eqref{eq:cform}, and let $t$ satisfy $1/n^{2/3}\ll T-t\le\ln^{2}n/n^{2/3}$. Then there exists a function $z(y)$ analytic in a neighborhood
of $\mathbb{R}$ such that 
\begin{itemize}
\item[(i)] $z(0)=\zeta$, 
\item[(ii)] $
3\sqrt{T-t}-(z(y)/c-\zeta/c)\notin(-\infty,0] \ \forall y\in\mathbb{R}$, and
\item[(iii)]
\begin{equation}
y=\frac{\sqrt{c^{3}\phi_{z^{3}}(\zeta(T),T)}}{\sqrt{6}}(z(y)/c-\zeta/c)\sqrt{3\sqrt{T-t}-(z(y)/c-\zeta/c)}.\label{eq:maincurve}
\end{equation}
\end{itemize}
\end{prop}

\begin{proof} Let $c,t$ and $T$ be as in the statement. Let $y\in\mathbb{R}$ and set
\begin{align}
A&=-\frac{\sqrt{6}}{\sqrt{c^{3}\phi_{z^{3}}(\zeta(T),T)}\sqrt{3\sqrt{T-t}}^{3}}, \nonumber\\
B&=\frac{1}{2}c^{3}\phi_{z^{3}}(\zeta(T),T)\sqrt{T-t}, \qquad \textrm{and} \nonumber \\
Z&=\frac{1}{3A\sqrt{B}}\left(\frac{r(y)}{\sqrt[3]{2}}+\frac{\sqrt[3]{2}}{r(y)}-1\right),\label{eq:yfunc}
\end{align}
where
\begin{equation}
r(y)^{3}=27A^{2}y^{2}-2+\sqrt{27}y\sqrt{27A^{4}y^{2}-4A^{2}}.\label{eq:aydef}
\end{equation}
The last remarks immediately preceding the proposition imply that $A^4 \in \mathbb{R}$, and that $A^2 \in i\mathbb{R}^-$. Consequently, $27A^{4}y^{2}-4A^{2}\notin\mathbb{R}$,
 which implies that $r(y)^{3}$ is analytic in a neighborhood of $\mathbb{R}$. The existence of an analytic cube root of $r(y)^3$ will follow once we show that  
 \[
27A^{2}y^{2}-2+\sqrt{27}y\sqrt{27A^{4}y^{2}-4A^{2}}\notin[0,\infty),
\]
since in this case we can choose $[0,\infty)$ as the cut to define $r(y)$. Suppose to the contrary, that 
 \[
27A^{2}y^{2}-2+\sqrt{27}y\sqrt{27A^4y^2-4A^2}\in [0,\infty).
\]
The identity 
\[
(\star) \quad \left(27A^{2}y^{2}-2+\sqrt{27}y\sqrt{27A^{4}y^{2}-4A^{2}}\right)\left(27A^{2}y^{2}-2-\sqrt{27}y\sqrt{27A^{4}y^{2}-4A^{2}}\right)=4 \qquad  (y \in \mathbb{R})
\]
implies that $r^{3}(y)\ne0 \forall y\in\mathbb{R}$, and $ 27A^{2}y^{2}-2-\sqrt{27}y\sqrt{27A^{4}y^{2}-4A^{2}}\in[0,\infty)$ for all $y \in \mathbb{R}$. Consequently, the sum of the two factors in $(\star)$ is real, and non-negative. But this sum is equal to $27A^2y^2-2$, which belongs to $\mathbb{C} \setminus [0,\infty)$, as $A^2 \in i\mathbb{R}^-$. We have reached a contradiction, and conclude that $r^3(y)$ has an analytic cube root $r(y)$ in a neighborhood of $\mathbb{R}$.

Define $z(y)$ by the relation
\[
Z(y)=\sqrt{B}(z(y)/c-\zeta/c).
\]
Since $Z(y)$ is analytic on a neighborhood or $\mathbb{R}$, so is $z(y)$. We now verify the claims (i)-(iii) in the statement of the proposition. For (i), we compute
\[
Z(0)=\frac{1}{3A\sqrt{B}}\left(\frac{r(0)}{\sqrt[3]{2}}+\frac{\sqrt[3]{2}}{r(0)}-1\right)=\frac{1}{3A\sqrt{B}}\left(\frac{\sqrt[3]{2}e^{\pi i/3}}{\sqrt[3]{2}}+\frac{\sqrt[3]{2}}{\sqrt[3]{2}e^{\pi i/3}}-1\right)=0,
\]
from which we readily deduce that $z(0)=\zeta$. For (ii), note that the definitions of $A,B,Z$ and $r^3(y)$ imply that
\[
y^{2}=AZ^{3}+Z^{2},
\]
or equivalently,
\[
y^{2}=\frac{c^{3}\phi_{z^{3}}(\zeta(T),T)}{6}(z/c-\zeta/c)^{2}\left(3\sqrt{T-t}-(z/c-\zeta/c)\right).
\]
It is immediate then that if $3\sqrt{T-t}-(z/c-\zeta/c) \in (-\infty,0]$, then $(z/c-\zeta/c) \in \mathbb{R}$, and hence we must have 
$\frac{c^{3}\phi_{z^{3}}(\zeta(T),T)}{6} \in\mathbb{R}$, which we know to be false. We conclude that $3\sqrt{T-t}-(z/c-\zeta/c) \notin (-\infty,0]$. To establish (iii), we note that given (ii), we may deduce from the equation above that either 
\[
y=\frac{\sqrt{c^{3}\phi_{z^{3}}(\zeta(T),T)}}{\sqrt{6}}(z(y)/c-\zeta/c)\sqrt{3\sqrt{T-t}-(z(y)/c-\zeta/c)}
\]
or 
\[
-y=\frac{\sqrt{c^{3}\phi_{z^{3}}(\zeta(T),T)}}{\sqrt{6}}(z(y)/c-\zeta/c)\sqrt{3\sqrt{T-t}-(z(y)/c-\zeta/c)}.
\]
In the first case the result follows. In the second case chose the the analytic function $z(-y)$ in place of $z(y)$. The proof is complete.
\end{proof}
We continue our discussion by developing a path of integration akin to $\Gamma$ from the previous section, with a central segment and two tails, and appropriate asymptotic expressions for the integral of $f(z,t)$ on each. We begin with the central segment. The function $z(y)$ as defined in Proposition \ref{lem:zycurve} is unbounded, as are the sets $\left\{ z(y)|y\in(0,\infty)\right\} $ and $\{z(y)|y\in(-\infty,0)\}$. Thus, given $n \in \mathbb{N}$ and $c$ as in equation \eqref{eq:cform}, there exists $a<0<b$ (depending on $n$, c.f. Lemma \ref{lem:ablowerbound}) such that 
\begin{enumerate}
\item $\left|z(y)/c-\zeta/c\right|<7(\ln n)/n^{1/3}$ , $\forall y\in(a,b)$, 
\item $|z(a)/c-\zeta/c|=|z(b)/c-\zeta/c|=7\ln n/n^{1/3}$. 
\end{enumerate}
This means that the curve $z(y)$ cuts the boundary of the ball centered at $\zeta(t)$ with radius $c \cdot 7 (\ln n)/n^{1/3}$ at $z(a)$ and $z(b)$ and for no other $y \in (a,b)$. While this ball may not, in its entirety, be contained in the first open quadrant, the curve $z(y)$ for $a<y<b$ is contained in the first open quadrant. 
\begin{lem}
\label{lem:firstquadrant} Let $n \gg 1$, and let $a<0<b$ be such that conditions (1) and (2) in the above discussion hold. Then for any $y \in [a,b],$ $z(y)$ lies in the
first quadrant.
\end{lem}

\begin{proof}
Recall that $1/n^{2/3}\ll T-t\le\ln^{2}n/n^{2/3}$, and hence $\zeta\rightarrow\zeta(T)$ as $n \to \infty$. In addition, the definition of $a$ and $b$ imply that $|z(y)-\zeta|\rightarrow0$ as $n \to \infty$. Since $\Im\zeta(T)>0$, we deduce that $\Im z(y)>0$ for all large $n$.
Note that 
\begin{align*}
\frac{\pi}{2}-\Im\phi(\zeta,t) & =\Im\phi(\zeta(T),T)-\Im\phi(\zeta,t).
\end{align*}
With $\zeta(T)-\zeta=c\sqrt{T-t}+\mathcal{O}(T-t)$, $\phi_{z}(\zeta(T),T)=0$,
$\phi_{z^{2}}(\zeta(T),T)=0$, and $\Im\phi_{t}(\zeta(T),T)=0$, we
apply the Taylor series expansion to the bivariate function $\phi(z,t)$
in a neighborhood of $(\zeta(T),T)$ to arrive at
\[
\Im\phi(\zeta(T),T)-\Im\phi(\zeta,t)=\Im\left((\zeta(T)-\zeta)(T-t)\phi_{z,t}(\zeta(T),T)+\mathcal{O}((T-t)^{2})\right)\asymp(T-t)^{3/2}\gg\frac{1}{n}.
\]
If there were a $y\in[a,b]$ such that $\Re z(y)=0$, then by the definition
of $\phi(z,t)$, we would have $\Im\phi(z(y),t)=\pi/2$, and consequently,
\[
\Im\phi(z(y),t)-\Im\phi(\zeta,t)=\frac{\pi}{2}-\Im\phi(\zeta,t) \gg\frac{1}{n}.
\]
This, however, contradicts the estimate (obtained from \eqref{eq:phiapprox} and \eqref{eq:maincurve})
\[
\Im\phi(z(y),t)-\Im\phi(\zeta,t)=\Im \left(y^{2}+\phi(z(y),t)-\phi(\zeta,t)\right)=\mathcal{O}\left(\frac{\ln^{4}n}{n^{4/3}}\right).
\]
Since $\Re z(0)=\Re \zeta(t) >0$, we conclude that $\Re z(y)>0$ for all $a<y<b$, and the result follows. 
\end{proof}
\begin{lem}
\label{lem:ablowerbound} Let $a$, $b$ be as in given in Lemma \ref{lem:firstquadrant}. Then $a\asymp\ln^{3/2}n/n^{1/2}$ and $b\asymp\ln^{3/2}n/n^{1/2}$. 
\end{lem}

\begin{proof}
We compute
\begin{align*}
\left|3\sqrt{T-t}-(z(b)/c-\zeta/c)\right|^{2} & =9(T-t)+|z(b)/c-\zeta/c|^{2}-6\sqrt{T-t}\Re(z(b)/c-\zeta/c)\\
 & \ge|z(b)/c-\zeta/c|(|z(b)/c-\zeta/c|-6\sqrt{T-t})\\
 & \geq \frac{7 \ln n}{n^{1/3}}\left(\frac{7 \ln n}{n^{1/3}}-6\frac{\ln n}{n^{1/3}} \right)\\
 & =\frac{7\ln^{2}n}{n^{2/3}},
\end{align*}
and
\[
|3\sqrt{T-t}-(z(b)/c-\zeta/c)|\le3\sqrt{T-t}+|z(b)/c-\zeta/c|\le\frac{10\ln n}{n^{1/3}}.
\]
Combining these estimates yields
\begin{align*}
b & =\frac{\sqrt{|c^{3}\phi_{z^{3}}(\zeta(T),T)|}}{\sqrt{6}}|z(b)/c-\zeta/c|\sqrt{|3\sqrt{T-t}-(z(b)/c-\zeta/c)|}\\
 & \asymp\frac{\ln^{3/2}n}{n^{1/2}}.
\end{align*}
Similar computations establish the claim for $a$ as well. 
\end{proof}
We now develop an estimate for the central integral (analogous to the estimate \eqref{eq:Gammaepsasymptotic} of section \ref{sec:maintermapprox}). Let $a,b$ and $z(y)$ be as above, and let $\Gamma(a,b)$ be the curve parameterized by $z(y)$, $a\le y\le b$. Since $\psi$ is analytic near $\zeta(T)$, we see that 
\begin{eqnarray} 
\int_{\Gamma(a,b)}e^{-n\phi(z,t)}\psi(z)dz&=&\psi(\zeta(T))\int_{\Gamma(a,b)}e^{-n\phi(z,t)}dz\left(1+o(1)\right) \nonumber \\
&=&\psi(\zeta(T))\int_{\Gamma(a,b)}e^{-n(\phi(z,t)-\phi(\zeta,t)+\phi(\zeta,t))}dz\left(1+o(1)\right). \label{eq:Gammaabest}
\end{eqnarray}
We use equation \eqref{eq:phiapprox} to estimate $\phi(z,t)-\phi(\zeta,t)$, along with the estimate 
\[
\exp\left(n\mathcal{O}\left((z-\zeta)^{2}(T-t)+(z-\zeta)^{3}\sqrt{T-t}+(z-\zeta)^{4}\right)\right)=e^{\mathcal{O}(\ln^{4}n/n^{1/3})}=1+\mathcal{O}(\ln^{4}n/n^{1/3})
\]
to conclude that the expression in \eqref{eq:Gammaabest} is equal to 
\begin{align*}
 & \psi(\zeta(T))e^{-n\phi(\zeta,t)}\int_{\Gamma(a,b)}\exp\left(-n\frac{\phi_{z^{3}}(\zeta(T),T)}{6}(z-\zeta)^{2}\left(3c\sqrt{T-t}-(z-\zeta)\right)\right)dz(1+o(1))\\
= & \psi(\zeta(T))e^{-n\phi(\zeta,t)}\int_{a}^{b}e^{-ny^{2}}z'(y)dy(1+o(1)).
\end{align*}

In order to find an expression for $z'(y)$, we square both sides of equation \eqref{eq:maincurve} and compute the derivatives with respective to $y$ to obtain 
\[
2y=\frac{c^{2}\phi_{z^{3}}(\zeta(T),T)}{2}(z/c-\zeta/c)\left(2\sqrt{T-t}-(z/c-\zeta/c)\right)z'(y),
\]
which we solve for $z'(y)$:
\[
z'(y)=\frac{4y}{c^{2}\phi_{z^{3}}(\zeta(T),T)(z/c-\zeta/c)\left(2\sqrt{T-t}-(z/c-\zeta/c)\right)}.
\]
Since
\begin{align*}
(z/c-\zeta/c)\left(2\sqrt{T-t}-(z/c-\zeta/c)\right) & =(z/c-\zeta/c)\left(3\sqrt{T-t}-(z/c-\zeta/c)\right)-\sqrt{T-t}(z/c-\zeta/c)\\
 & =\frac{6y^{2}}{c^{3}\phi_{z^{3}}(\zeta(T),T)(z/c-\zeta/c)}-\sqrt{T-t}(z/c-\zeta/c),
\end{align*}
we conclude that
\begin{align}
z'(y) & =\frac{4c}{\sqrt{c^{3}\phi_{z^{3}}(\zeta(T),T)}}\left(\frac{6y}{\sqrt{c^{3}\phi_{z^{3}}(\zeta(T),T)}(z/c-\zeta/c)}-\frac{\sqrt{T-t}(z/c-\zeta/c)\sqrt{c^{3}\phi_{z^{3}}(\zeta(T),T)}}{y}\right)^{-1} \label{eq:zprimeexp}\\
 & =\frac{4c}{\sqrt{6}\sqrt{c^{3}\phi_{z^{3}}(\zeta(T),T)}}\left(\sqrt{3\sqrt{T-t}-(z/c-\zeta/c)}-\frac{\sqrt{T-t}}{\sqrt{3\sqrt{T-t}-(z/c-\zeta/c)}}\right)^{-1}, \nonumber
\end{align}
and consequently 
\begin{align}
 & \int_{\Gamma(a,b)}e^{-n\phi(z,t)}\psi(z)dz\nonumber \\
= & \frac{4c\psi(\zeta(T))e^{-n\phi(\zeta,t)}}{\sqrt{6}\sqrt{c^{3}\phi_{z^{3}}(\zeta(T),T)}}\int_{a}^{b}e^{-ny^{2}}\left(\sqrt{3\sqrt{T-t}-(z/c-\zeta/c)}-\frac{\sqrt{T-t}}{\sqrt{3\sqrt{T-t}-(z/c-\zeta/c)}}\right)^{-1}dy(1+o(1)).\label{eq:mainint}
\end{align}
In an effort to develop an asymptotic lower bound for $\int_{\Gamma(a,b)}e^{-n\phi(z,t)}\psi(z)dz$, we start with the following result.
\begin{lem}
\label{lem:posintegrand} Let $a,b$ be as in Lemma \ref{lem:firstquadrant}. Then for any $y\in(a,b)$, 
\[
\Re\left(\sqrt{3\sqrt{T-t}-(z/c-\zeta/c)}-\frac{\sqrt{T-t}}{\sqrt{3\sqrt{T-t}-(z/c-\zeta/c)}}\right)\ge0.
\]
\end{lem}
\begin{rem} The conclusion of Lemma \ref{lem:posintegrand} is equivalent to the non-negativity of the real part of the integrand of the right hand side of \eqref{eq:mainint}, but is easier to establish.
\end{rem}
\begin{proof}
Since $\sqrt{T-t}\in\mathbb{R}$, 
\[
\Arg\sqrt{3\sqrt{T-t}-(z/c-\zeta/c)}=-\Arg\frac{\sqrt{T-t}}{\sqrt{3\sqrt{T-t}-(z/c-\zeta/c)}},
\]
and this common argument lies between $-\pi/2$ and $\pi/2$. Consequently,
the sign of 
\[
\Re\left(\sqrt{3\sqrt{T-t}-(z/c-\zeta/c)}\right)-\Re\left(\frac{\sqrt{T-t}}{\sqrt{3\sqrt{T-t}-(z/c-\zeta/c)}}\right)
\]
is the same as the sign of 
\begin{eqnarray}
 &  & \left|\sqrt{3\sqrt{T-t}-(z/c-\zeta/c)}\right|-\left|\frac{\sqrt{T-t}}{\sqrt{3\sqrt{T-t}-(z/c-\zeta/c)}}\right|\nonumber \\
 & = & \left|\frac{1}{\sqrt{3\sqrt{T-t}-(z/c-\zeta/c)}}\right|E(y),\label{eq:modulusdiff}
\end{eqnarray}
where
\[
E(y):=\left|3\sqrt{T-t}-(z/c-\zeta/c)\right|-\sqrt{T-t}.
\]
We now argue that $E(y) \geq 0$ for all $y \in [a,b]$. We first check the endpoints. If $y=a$ or $y=b$, then
$|z(y)/c-\zeta/c|=\frac{7\ln n}{n^{1/3}}$, and hence 
\[
\left|3\sqrt{T-t}-(z(y)/c-\zeta/c)\right|-\sqrt{T-t}\geq|z(y)/c-\zeta/c|-4\sqrt{T-t}\ge\frac{3\ln n}{n^{1/3}}>0.
\]
Next we demonstrate that the condition $E(y)=0$ holds at no more than one point in $[a,b]$. To this end, assume that $E(y)=0$. Since 
\[
\Re\left(3\sqrt{T-t}-(z/c-\zeta/c)\right)=3\sqrt{T-t}-\Re(z/c-\zeta/c),
\]
the assumption
\[
\left|3\sqrt{T-t}-(z/c-\zeta/c)\right|=\sqrt{T-t}
\]
implies that $\Re(z/c-\zeta/c)>0$, or equivalently $-\pi/2<\Arg(z/c-\zeta/c)<\pi/2$.
In order to obtain more information on the quantity $\Arg(z/c-\zeta/c)$, consider
the identity 
\begin{equation}
-(z/c-\zeta/c)e^{-i\Arg(z/c-\zeta/c)}=(3\sqrt{T-t}-(z/c-\zeta/c))e^{-i\Arg(z/c-\zeta/c)}-3\sqrt{T-t}e^{-i\Arg(z/c-\zeta/c)}.\label{eq:trianglerel}
\end{equation}
Taking the imaginary parts of both sides of \eqref{eq:trianglerel} and using the assumption $E(y)=0$ we conclude that
\begin{align}
0=\sin\left(\Arg(3\sqrt{T-t}-(z/c-\zeta/c))-\Arg(z/c-\zeta/c)\right) & +3\sin\left(\Arg(z/c-\zeta/c)\right).\label{eqn:auxeq}
\end{align}
On the other hand, taking the arguments of both sides of \eqref{eq:maincurve}
yields 
\begin{equation}
\frac{\pi}{4}+\Arg(z/c-\zeta/c)+\frac{1}{2}\Arg(3\sqrt{T-t}-(z/c-\zeta/c))=\begin{cases}
0 & \text{ if }y>0\\
\pm\pi & \text{ if }y<0
\end{cases},\label{eq:angleeq}
\end{equation}
which implies that 
\[
\Arg(3\sqrt{T-t}-(z/c-\zeta/c))=-\frac{\pi}{2}-2\Arg(z/c-\zeta/c)\pmod{2\pi}.
\]
Substituting into equation \eqref{eqn:auxeq} and rearranging gives
the equation 
\[
0=\cos(3\Arg(z/c-\zeta/c))+3\sin(\Arg(z/c-\zeta/c)),
\]
which has exactly one solution $\Arg(z/c-\zeta/c)=-0.247872...$ on
$(-\pi/2,\pi/2)$. Now we take the real parts of both sides of \eqref{eq:trianglerel}
and employ the assumption $E(y)=0$ to solve for $|z/c-\zeta/c|$: 
\begin{align*}
-\frac{|z/c-\zeta/c|}{\sqrt{T-t}} & =\cos\left(\Arg(3\sqrt{T-t}-(z/c-\zeta/c))-\Arg(z/c-\zeta/c)\right)-3\cos\left(\Arg(z/c-\zeta/c)\right)\\
 & =-\sin(3\Arg(z/c-\zeta/c))-3\cos\left(\Arg(z/c-\zeta/c)\right)
\end{align*}
from which we conclude that $E(y)=0$ on $(a,b)$ can only occur at $y$ corresponding to $\Arg(z/c-\zeta/c)=-0.247872\ldots$.
The result now follows from the continuity of $E(y)$. 
\end{proof}
\begin{rem}
Using \eqref{eq:angleeq} and the fact that a possible zero of $E(y)$ on $[a,b]$
satisfies $\Arg(z/c-\zeta/c)=-0.247872\ldots$, we conclude that $E(y)$ has no zero in $y$ on $(a,0)$. 
\end{rem}
The following lemma gives the range of $\Arg (z/c-\zeta/c)$ for $y \in (a,0)$.
\begin{lem}
\label{lem:anglez-zeta}For any $y\in(a,0)$, $\pi>\Arg(z/c-\zeta/c)\ge\pi/4$
or $-\pi<\Arg(z/c-\zeta/c)\le-3\pi/4$. 
\end{lem}

\begin{proof} For ease of notation, set
\[
\bigboxvoid=\sqrt{3\sqrt{T-t}-(z/c-\zeta/c)}-\frac{\sqrt{T-t}}{\sqrt{3\sqrt{T-t}-(z/c-\zeta/c)}}.
\]
Equation \eqref{eq:zprimeexp} implies that
\[
z'(y)\frac{\sqrt{c^{3}\phi_{z^{3}}(\zeta(T),T)}}{c}=\frac{4}{\sqrt{6}}\left(\bigboxvoid \right)^{-1}.
\]
By Lemma \ref{lem:posintegrand} $\Re \bigboxvoid \geq 0$, and hence 

\[
\Re \frac{4}{\sqrt{6}}\left( \bigboxvoid \right)^{-1}=\Re\left(\frac{\sqrt{c^{3}\phi_{z^{3}}(\zeta(T),T)}}{c}z'(y)\right)\ge0.
\]
Combining this inequality with 
\[
\Re\left(\sqrt{c^{3}\phi_{z^{3}}(\zeta(T),T)}(z(0)/c-\zeta/c)\right)=0
\]
gives 
\[
\Re\left(\sqrt{c^{3}\phi_{z^{3}}(\zeta(T),T)}(z(y)/c-\zeta/c)\right)\le0,\qquad\forall y\in(a,0).
\]
The lemma now follows from the fact that 
\[
\Arg\left(\sqrt{c^{3}\phi_{z^{3}}(\zeta(T),T)}\right)=\pi/4.
\]
\end{proof}
\begin{lem}
\label{lem:maintermlowerbound} Let $z(y)$ be the function defined on $[a,b]$ by Proposition \ref{lem:zycurve}, where $a,b$ are as in the discussion preceding Lemma \ref{lem:firstquadrant}. Then 
\[
\left|\int_{a}^{b}\left(\sqrt{3\sqrt{T-t}-(z/c-\zeta/c)}-\frac{\sqrt{T-t}}{\sqrt{3\sqrt{T-t}-(z/c-\zeta/c)}}\right)^{-1}dy\right|\gg\frac{1}{n^{1/3}\sqrt{\ln n}}, \qquad \textrm{as} \quad n \to \infty.
\]
\end{lem}

\begin{proof}
By Lemma \ref{lem:posintegrand}, we have 
\begin{align*}
 & \left|\int_{a}^{b}\left(\sqrt{3\sqrt{T-t}-(z/c-\zeta/c)}-\frac{\sqrt{T-t}}{\sqrt{3\sqrt{T-t}-(z/c-\zeta/c)}}\right)^{-1}dy\right|\\
\ge & \int_{a}^{0}\Re\left(\sqrt{3\sqrt{T-t}-(z/c-\zeta/c)}-\frac{\sqrt{T-t}}{\sqrt{3\sqrt{T-t}-(z/c-\zeta/c)}}\right)^{-1}dy.
\end{align*}
The integrand in the last integral is the quotient of 
\[
\left(\left|\sqrt{3\sqrt{T-t}-(z/c-\zeta/c)}\right|-\frac{\sqrt{T-t}}{|\sqrt{3\sqrt{T-t}-(z/c-\zeta/c)}|}\right)\cos\left(\Arg\sqrt{3\sqrt{T-t}-(z/c-\zeta/c)}\right)
\]
and 
\[
\left|\sqrt{3\sqrt{T-t}-(z/c-\zeta/c)}-\frac{\sqrt{T-t}}{\sqrt{3\sqrt{T-t}-(z/c-\zeta/c)}}\right|^{2}.
\]
Rearranging this quotient gives
\begin{equation}
\frac{\left(\sqrt{|3\sqrt{T-t}-(z/c-\zeta/c)}|^{3}-\sqrt{T-t}\left|\sqrt{3\sqrt{T-t}-(z/c-\zeta/c)}\right|\right)\cos\left(\Arg\sqrt{3\sqrt{T-t}-(z/c-\zeta/c)}\right)}{|2\sqrt{T-t}-(z/c-\zeta/c)|^{2}}.\label{eq:reintegrand}
\end{equation}
To bound the denominator, we compute
\[
|2\sqrt{T-t}-(z/c-\zeta/c)|^{2} \leq \left(\frac{2\ln n}{n^{1/3}}+\frac{7\ln n}{n^{1/3}}\right)^{2}=\frac{81\ln^{2}n}{n^{2/3}}.
\]
Turning our attention to the numerator, we recall Lemma \ref{lem:anglez-zeta}, which -- since $y \in (a,0)$ -- implies that
\[
-\frac{3\pi}{4}\le\Arg\left(3\sqrt{T-t}-(z/c-\zeta/c)\right)\le\frac{\pi}{4},
\]
and consequently 
\begin{equation}
\cos\left(\Arg\sqrt{3\sqrt{T-t}-(z/c-\zeta/c)}\right)\ge\cos\frac{3\pi}{8}.\label{eq:cosineq}
\end{equation}
We now find a lower bound for the remaining part of the numerator. To this end, note that 
\begin{align*}
\left|3\sqrt{T-t}-(z/c-\zeta/c)\right|^{2} & =9(T-t)+|z/c-\zeta/c|^{2}-6\sqrt{T-t}\Re(z/c-\zeta/c)\\
 & =9(T-t)+|z/c-\zeta/c|^{2}-6\sqrt{T-t}|z/c-\zeta/c|\cos\Arg(z/c-\zeta/c).
\end{align*}
By Lemma \ref{lem:anglez-zeta}, $\cos\left(z/c-\zeta/c \right)\leq \sqrt{2}/2$. Therefore
\begin{align}
& 9(T-t)+|z/c-\zeta/c|^{2}-6\sqrt{T-t}|z/c-\zeta/c|\cos\Arg(z/c-\zeta/c) \nonumber \\ & \geq \underbrace{9(T-t)-3\sqrt{2}\sqrt{T-t}|z/c-\zeta/c|+|z/c-\zeta/c|^{2}}_{\textrm{quadratic in}\  |z/c-\zeta/c|}
 \ge\frac{9}{2}(T-t).\label{eq:quadineq}
\end{align}
On the other hand,
\[
\left|\sqrt{3\sqrt{T-t}-(z/c-\zeta/c)}\right|^{3}-\sqrt{T-t}\left|\sqrt{3\sqrt{T-t}-(z/c-\zeta/c)}\right|
\]
is an increasing function in terms of $\left|\sqrt{3\sqrt{T-t}-(z/c-\zeta/c)}\right|$
when this quantity is at least $\sqrt{\sqrt{T-t}/3}$. Since $\sqrt{9/2(T-t)} >\sqrt{\sqrt{T-t}/3}$, we conclude that
\begin{align*}
 & \left|\sqrt{3\sqrt{T-t}-(z/c-\zeta/c)}\right|^{3}-\sqrt{T-t}\left|\sqrt{3\sqrt{T-t}-(z/c-\zeta/c)}\right|\\
\ge & \left(\frac{9}{2}(T-t)\right)^{1/4}\left(\frac{3}{\sqrt{2}}\sqrt{T-t}-\sqrt{T-t}\right)\asymp(T-t)^{3/4}\gg1/\sqrt{n}.
\end{align*}
Putting all this together shows that the expression in \eqref{eq:reintegrand} is at least a constant multiple of $\displaystyle{\frac{n^{1/6}}{\ln^{2}n}}$. The result now follows from Lemma \ref{lem:ablowerbound}, as 
\[
\int_0^a \frac{n^{1/6}}{\ln^2 n}dy=a\frac{n^{1/6}}{\ln^2 n} \asymp \frac{1}{m^{1/3}\sqrt{\ln m} }.
\] 
\end{proof}
In order to extend $\Gamma_{(a,b)}$ to the point at infinity we use arguments similar to those in Section 2.3 using $a$ and $b$ instead
of $\pm\epsilon$ and extend the curve $z(y)$, $a\le y\le b$,
by adding two tails $\Gamma_{a}$and $\Gamma_{b}$ starting from
$z(a)$ and $z(b)$ to $\infty$ in the lower and upper half planes respectively, 
such that $\Re\phi(z)\ge\Re\phi(z(a))+a^{2}$, $\forall z\in\Gamma_{a}$
and $\Re\phi(z)\ge\Re\phi(z(b))+b^{2},\forall z\in\Gamma_{b}$ with
equality only at $z=z(a)$ or $z=z(b)$ respectively. Repeating the arguments
which provided the asymptotic bound for \eqref{eq:tailint} we obtain  
\[
\int_{\Gamma_{a}}e^{-n\phi(z,t)}\psi(z)dz=\mathcal{O}\left(e^{-n\phi(\zeta,t)-na^{2}/2}\right)
\]
and 
\[
\int_{\Gamma_{b}}e^{-n\phi(z,t)}\psi(z)dz=\mathcal{O}\left(e^{-n\phi(\zeta,t)-nb^{2}/2}\right),
\]
which are valid since 
\begin{align*}
na^2 & \asymp \sqrt{n} \ln ^{3/2} n \qquad \textrm{as} \quad n \to \infty, \qquad \textrm{and} \\
nb^2 &\asymp \sqrt{n} \ln ^{3/2} n \qquad \textrm{as} \quad n \to \infty.
\end{align*}
Finally, from equation \eqref{eq:mainint} and Lemmas \ref{lem:ablowerbound}
and \ref{lem:maintermlowerbound} we obtain the estimate
\[
\int_{\Gamma_{2}}f(z,t)dz=K_c(n)\int_{a}^{b}e^{-ny^{2}}\left(\sqrt{3\sqrt{T-t}-(z/c-\zeta/c)}-\frac{\sqrt{T-t}}{\sqrt{3\sqrt{T-t}-(z/c-\zeta/c)}}\right)^{-1}dy(1+o(1)),
\]
where
\[
K_c(n):=\frac{4c\psi(\zeta(T))e^{-n\phi(\zeta,t)}}{\sqrt{6}\sqrt{c^{3}\phi_{z^{3}}(\zeta(T),T)}}.
\]

\subsubsection{\label{sec:SmallT1-t}The asymptotics when $T=T_2$ and $\frac{1}{n}\ll|T_{1}-t|\le\ln^{2}n/n^{2/3}$}

We focus on the case $t<T_1$ -- the same arguments will apply
for the case $T_{1}<t$. Recall from Lemma \ref{lem:asympzetaT1} that 
\begin{align*}
\zeta(t)-\zeta(T_{1})= d\sqrt{T_{1}-t}+\mathcal{O}(T_{1}-t),
\end{align*}
where $d^{3}\cdot \phi_{z^{3}}(\zeta(T_{1}),T_{1})\in\mathbb{R}^{+}$.
Following arguments analogous to those in the previous section, replacing $c$ with $d$, and $T$ with $T_1$, we conclude that 
\[
\int_{\Gamma_{2}}f(z,t)dz=K_d(n)\int_{a}^{b}e^{-ny^{2}}\left(\sqrt{3\sqrt{T_{1}-t}-(z/d-\zeta/d)}-\frac{\sqrt{T_{1}-t}}{\sqrt{3\sqrt{T_{1}-t}-(z/d-\zeta/d)}}\right)^{-1}dy(1+o(1)),
\]
where
\[
K_d(n)=\frac{4d\psi(\zeta(T_{1}))e^{-n\phi(\zeta,t)}}{\sqrt{6}\sqrt{d^{3}\phi_{z^{3}}(\zeta(T_{1}),T_{1})}}.
\]
 For the sake of brevity, instead of reproducing the argument in its entirety, we content ourselves with highlighting the differences. The curve $z(y)$ (c.f. Proposition \ref{lem:zycurve}) is now only piecewise smooth, and is defined
by replacing the function in \eqref{eq:aydef} by 
\[
r(y)^{3}=\begin{cases}
27A^{2}y^{2}-2+\sqrt{27}yi\sqrt{4A^{2}-27A^{4}y^{2}} & \text{ if }27A^{4}y^{2}-4A^{2}<0\\
27A^{2}y^{2}-2+\sqrt{27}y\sqrt{27A^{4}y^{2}-4A^{2}} & \text{ if }27A^{4}y^{2}-4A^{2}\ge0
\end{cases}.
\]
When $27A^{4}y^{2}-4A^{2}<0$ we use the cut $[0,\infty)$ to
define $r(y)$ and when $27A^{4}y^{2}-4A^{2}\ge0$, we define $r(y)$
by 
\[
r(y)=\begin{cases}
\sqrt[3]{27A^{2}y^{2}-2+\sqrt{27}y\sqrt{27A^{4}y^{2}-4A^{2}}} & \text{ if }y>0\\
e^{-2\pi i/3}\sqrt[3]{27A^{2}y^{2}-2+\sqrt{27}y\sqrt{27A^{4}y^{2}-4A^{2}}} & \text{ if }y<0
\end{cases}.
\]
Lemma \ref{lem:firstquadrant} holds trivially as $\zeta(T_{1})$
lies in the open first quadrant. In the proof of Lemma \ref{lem:posintegrand}, equation \eqref{eq:angleeq}
becomes 
\[
\Arg(z/d-\zeta/d)+\frac{1}{2}\Arg(3\sqrt{T_{1}-t}-(z/d-\zeta/d))=\begin{cases}
0 & \text{ if }y>0\\
\pm\pi & \text{ if }y<0
\end{cases},
\]
which implies that 
\[
\Arg(3\sqrt{T_{1}-t}-(z/d-\zeta/d))=-2\Arg(z/d-\zeta/d)\pmod{2\pi}.
\]
The analogue of equation \eqref{eqn:auxeq} reads
\[
0=\sin(3\Arg(z/d-\zeta/d))+3\sin(\Arg(z/d-\zeta/d)),
\]
which has exactly one solution, namely $\Arg(z/d-\zeta/d)=0$, on $(-\pi/2,\pi/2)$. Using this solution we also find $|z/d-\zeta/d|$. 

Since $d^{3} \cdot \phi_{z^{3}}(\zeta(T_{1}),T_{1})\in\mathbb{R^{+}}$, the
result corresponding to Lemma \ref{lem:anglez-zeta} is that for $y\in(a,0)$
\[
\Re(z/d-\zeta/d)\le0,
\]
while the inequality corresponding to \eqref{eq:cosineq} is 
\[
\cos \left(\Arg\sqrt{3\sqrt{T_{1}-t}-(z/d-\zeta/d)}\right) \ge\frac{\sqrt{2}}{2}.
\]
Using that $\cos\Arg(z/d-\zeta/d)\le0$, inequality \eqref{eq:quadineq} is replaced
by 
\[
9(T_{1}-t)+|z/d-\zeta/d|^{2}\ge9(T_{1}-t),
\]
leading to the asymptotic lower bound
\begin{align*}
 & \left|\sqrt{3\sqrt{T_{1}-t}-(z/d-\zeta/d)}\right|^{3}-\sqrt{T_{1}-t}\left|\sqrt{3\sqrt{T_{1}-t}-(z/d-\zeta/d)}\right|\\
\ge & \left(9(T_{1}-t)\right)^{1/4}\left(3\sqrt{T_{1}-t}-\sqrt{T_{1}-t}\right)\asymp(T_{1}-t)^{3/4}\gg1/n^{3/4}.
\end{align*}
With these differences, the estimate in the statement of Lemma \ref{lem:maintermlowerbound} becomes
\[
\left|\int_{a}^{b}\left(\sqrt{3\sqrt{T_{1}-t}-(z/d-\zeta/d)}-\frac{\sqrt{T_{1}-t}}{\sqrt{3\sqrt{T_{1}-t}-(z/d-\zeta/d)}}\right)^{-1}dy\right|\gg\frac{1}{n^{7/12}\sqrt{\ln n}}.
\]

For the case $t>T_{1}$, we note that as $t\rightarrow T_{1}$

\begin{align*}
\sqrt{1-2t^{2}-2t\sqrt{t^{2}-T_{1}^{2}}} & =\sqrt{(1-2T_{1}^{2})\left(1-2T_{1}\frac{\sqrt{t^{2}-T_{1}^{2}}}{1-2T_{1}^{2}}+\mathcal{O}(t-T_{1})\right)}\\
 & =\sqrt{1-2T_{1}^{2}}\left(1-T_{1}\frac{\sqrt{t^{2}-T_{1}^{2}}}{1-2T_{1}^{2}}+\mathcal{O}(t-T_{1})\right),
\end{align*}
and consequently equation \eqref{eq:zeta1defnT1} yields 
\begin{align*}
\zeta(t)-\zeta(T_{1})= & \frac{z_{1}+z_{2}}{2}\left(i\sqrt{t^{2}-T_{1}^{2}}-\frac{T_{1}\sqrt{t^{2}-T_{1}^{2}}}{\sqrt{1-2T_{1}^{2}}}+\mathcal{O}(t-T_{1})\right)\\
= & \widehat{d} \cdot \sqrt{t-T_{1}}+\mathcal{O}(t-T_{1})
\end{align*}
where 
\begin{equation}
\widehat{d}=\frac{(z_{1}+z_{2})\sqrt{T_{1}}}{\sqrt{2}}\left(i-\frac{T_{1}}{\sqrt{1-2T_{1}^{2}}}\right).\label{eq:eform}
\end{equation}

In order to be able to extend the curve $z(y)$, $a<y<b$ with two tails going to infinity we replace Lemma \ref{lem:gammae-int} with the following result.
\begin{lem}
\label{lem:zabsmallT1-t} Let $\mathcal{S}$ be defines as in Lemma \ref{lem:gammae-int}, and let $a,b$ be defined as in the discussion preceding Lemma \ref{lem:firstquadrant}. If $\frac{1}{n}\ll T_{1}-t\le\ln^{2}n/n^{2/3}$,
then $z(a)$ and $z(b)$ lie outside the region $\mathcal{S}$. 
\end{lem}

\begin{proof}
We recall that 
\[
y=\frac{\sqrt{d^{3}\phi_{z^{3}}(\zeta(T_{1}),T_{1})}}{\sqrt{6}}(z(y)/d-\zeta/d)\sqrt{3\sqrt{T_1-t}-(z(y)/d-\zeta/d)}
\]
and 
\begin{align*}
A & =\frac{\sqrt{6}}{\sqrt{d^{3}\phi_{z^{3}}(\zeta(T_1),T_1)}\sqrt{3\sqrt{T_1-t}^{3}}}\in\mathbb{R}^{+},\\
B & =\frac{1}{2}d^{3}\phi_{z^{3}}(\zeta(T_1),T_1)\sqrt{T_1-t}\in\mathbb{R}^{+}.
\end{align*}
Using the conditions
\[
|z(a)/d-\zeta/d|=|z(b)/d-\zeta/d|=7\ln n/n^{1/3},
\]
we conclude that when $y=a$, 
\[
A^{2}a^{2}=|A^{2}a^{2}|\ge\frac{7^{2}\cdot4}{3},
\]
which in turn implies that $27A^{4}a^{2}-4A^{2}>0$. Combining this with the similar
inequality $27A^{4}b^{2}-4A^{2}>0$, we deduce that 
\begin{equation}
r(y)=\begin{cases}
\sqrt[3]{27A^{2}b^{2}-2+\sqrt{27}b\sqrt{27A^{4}b^{2}-4A^{2}}} & \text{ if }y=b\\
e^{-2\pi i/3}\sqrt[3]{27A^{2}a^{2}-2+\sqrt{27}a\sqrt{27A^{4}a^{2}-4A^{2}}} & \text{ if }y=a
\end{cases}.\label{eq:rab}
\end{equation}

Lemma \ref{lem:asympzetaT1} provides for any $|\tau-T_{1}|=o(1)$ the estimate
\[
\zeta_{2}(\tau)-\zeta_{1}(t)\sim\begin{cases}
-d\sqrt{T_{1}-\tau}-d\sqrt{T_{1}-t} & \text{ if }\tau<T_{1}\\
id\sqrt{\tau-T_{1}}-d\sqrt{T_{1}-t} & \text{ if }\tau>T_{1}
\end{cases},
\]
from which we deduce that if $z\in\mathcal{S}$ and $z-\zeta=o(1)$,
then $(z-\zeta)/d$ lies in the second quadrant. Consequently, the expression \eqref{eq:rab} when $y=b$ together with
\[
\frac{z(b)-\zeta}{d}=\frac{1}{3A\sqrt{B}}\left(\frac{r(b)}{\sqrt[3]{2}}+\frac{\sqrt[3]{2}}{r(b)}-1\right)>0
\]
implies that $z(b)\notin\mathcal{S}.$ Using an analogous argument we also conclude $(z(a)-\zeta)/d$
lies in the third quadrant (see discussion after equation \eqref{eq:aydef} for why $\sqrt[3]{2}/r(a)$ is greater than 1 in modulus), and hence $z(a) \notin \mathcal{S}$.  The proof is complete.
\end{proof}
Repeating the arguments provided for the case $t<T_{1}$ \textit{mutatis mutandis}, we conclude that when $t>T_{1}$, 
\[
\int_{\Gamma_{2}}f(z,t)dz=K_{\widehat{d}}(n)\int_{a}^{b}e^{-ny^{2}}\left(\sqrt{3\sqrt{t-T_{1}}-(z/\widehat{d}-\zeta/\widehat{d})}-\frac{\sqrt{t-T_{1}}}{\sqrt{3\sqrt{t-T_{1}}-(z/\widehat{d}-\zeta/\widehat{d})}}\right)^{-1}dy(1+o(1)).
\]
We hasten to note that Lemma \ref{lem:zabsmallT1-t} is not necessary
in this case by the first paragraph in the proof of Lemma \ref{lem:gammae-int}. This completes the asymptotic analysis of the key integral away from the point $t=0$. We deal with this range in the next section.

\subsubsection{\label{sec:smallt} The asymptotics when $t\ll\ln^{4}n/n$}

\begin{figure}
\begin{centering}
\includegraphics[scale=0.2]{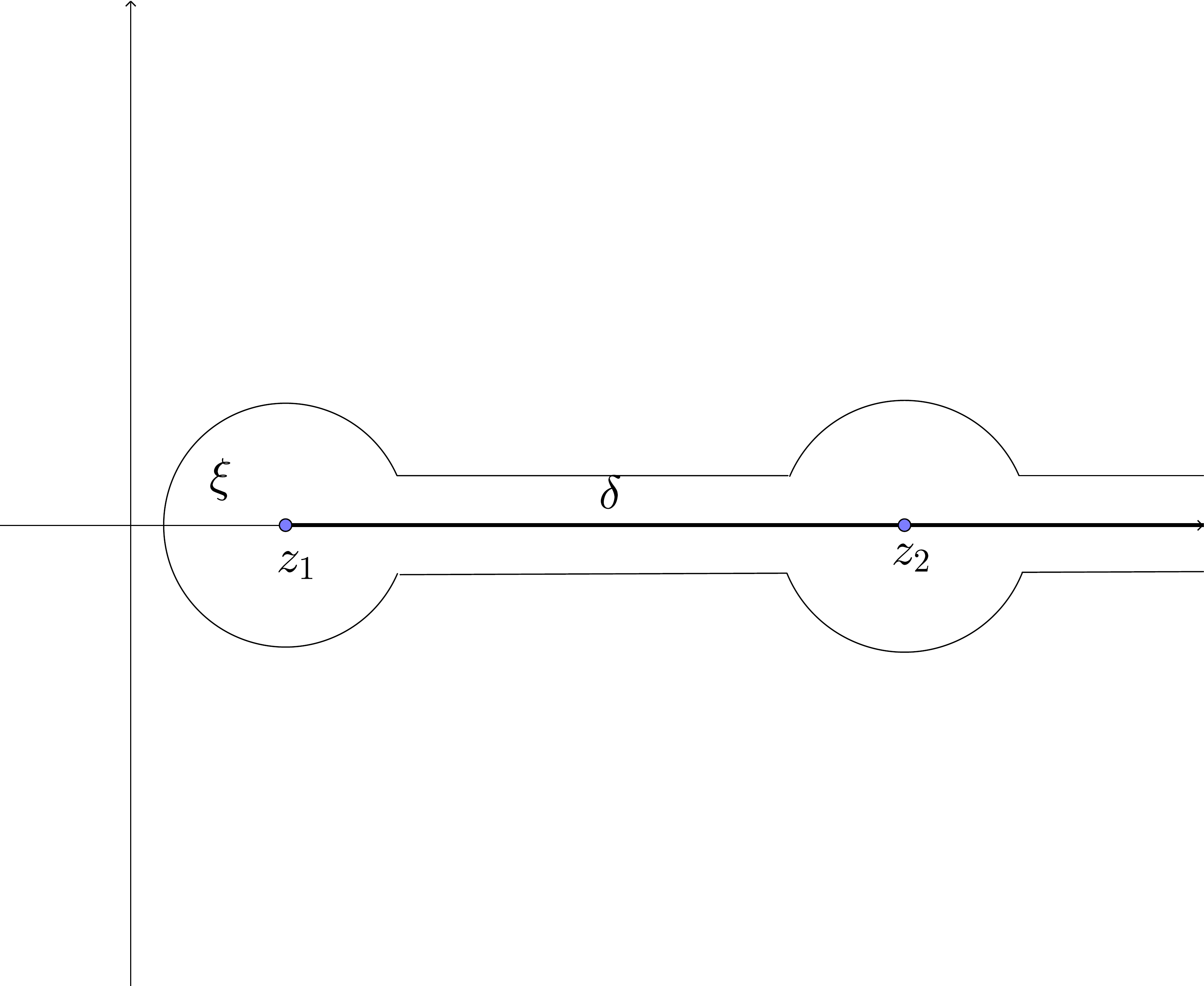} 
\par\end{centering}
\caption{\label{fig:ContourCurve}Contour curve}
\end{figure}

With a slight abuse of notation we let $\Gamma_{2}$ be the curve in Figure \ref{fig:ContourCurve}
where each circle around $z_{1}$ and $z_{2}$ has small radius
$\xi$, and each horizontal line segment has distance $\delta$ from
the $x$-axis.

On each arc $\gamma$ of the circles with radius $\xi$ centered at $z_1$ and $z_2$ we employ the basic estimates
\[
\left|e^{-n\phi(z,t)}\right|\le|z|^{-n}e^{4\pi t}\qquad\text{and}\qquad|\psi(z)|=\frac{|Q(z)Q(-z)|^{1/2}}{|z|}
\]
to conclude that 
\[
\int_{\gamma}f(z,t)dz\rightarrow0
\]
as $\xi\rightarrow0$. By letting $\delta\rightarrow0$ and $\xi\rightarrow0$
we may then rewrite the integral
\[
\int_{\Gamma_{2}}\frac{Q(z)^{1/2+int}Q(-z)^{1/2-int}}{z^{n+1}}dz
\]
as 
\begin{align*}
 & \left(e^{-i\pi(1/2+int)}-e^{i\pi(1/2+int)}\right)\int_{z_{1}}^{z_{2}}(z-z_{1})^{1/2+int}(z_{2}-z)^{1/2+int}(z_{1}+z)^{1/2-int}(z_{2}+z)^{1/2-int}\frac{dz}{z^{n+1}}\\
+ & (e^{-2i\pi(1/2+int)}-e^{2i\pi(1/2+int)}\int_{z_{2}}^{\infty}(z-z_{1})^{1/2+int}(z-z_{2})^{1/2+int}(z_{1}+z)^{1/2-int}(z_{2}+z)^{1/2-int}\frac{dz}{z^{n+1}}.
\end{align*}
We make the substitutions $w=z/z_{1}$ and $w=z/z_{2}$ in the first
and the second integral respectively to arrive at the expression 
\begin{align}
 & \frac{-i\left(e^{n\pi t}+e^{-n\pi t}\right)}{z_{1}^{n}}\int_{1}^{z_{2}/z_{1}}((wz_{1}-z_{1})(z_{2}-wz_{1}))^{1/2+int}((z_1+wz_1)(z_2+wz_1))^{1/2-int}\frac{dw}{w^{n+1}}\nonumber \\
+ & \frac{(-e^{2\pi nt}+e^{-2\pi nt})}{z_{2}^{n}}\int_{1}^{\infty}((wz_{2}-z_{1})(wz_{2}-z_{2}))^{1/2+int}((z_{1}+wz_{2})(z_{2}+wz_{2}))^{1/2-int}\frac{dw}{w^{n+1}}.\label{eq:sumtwoint}
\end{align}
We claim that the first summand in \eqref{eq:sumtwoint} is asymptotically equivalent to 
\begin{equation}
\frac{-i(e^{n\pi t}+e^{-n\pi t})}{z_{1}^{n-1}}2^{1/2-int}(z_{2}-z_{1})^{1/2+int}(z_{2}+z_{1})^{1/2-int}\Gamma(3/2+int)e^{-(3/2+int)\ln n}.\label{eq:firstint}
\end{equation}
In order to demonstrate the claim, we first apply the substitution $e^{z}=w$ to transform the first summand in \eqref{eq:sumtwoint} into
\[
\frac{-i\left(e^{n\pi t}+e^{-n\pi t}\right)}{z_{1}^{n}}\int_{0}^{\ln(z_{2}/z_{1})}\left(z\left(\frac{z_{1}e^{z}-z_{1}}{z}\right)(z_{2}-e^{z}z_{1})\right)^{1/2+int}((z_{1}+e^{z}z_{1})(z_{2}+e^{z}z_{1}))^{1/2-int}\frac{dz}{e^{nz}}.
\]
Then we split the range of integration into the intervals $(0,\ln^{5}n/n)$ and $(\ln^{5}n/n,\ln(z_{2}/z_{1}))$.
The contribution corresponding to the second interval is 
\[
\frac{-i\left(e^{n\pi t}+e^{-n\pi t}\right)}{z_{1}^{n}}\mathcal{O}\left(e^{-\ln^{5}n}\right),
\]
since 
\[
\left| \left(z\left(\frac{z_{1}e^{z}-z_{1}}{z}\right)(z_{2}-e^{z}z_{1})\right)^{1/2+int}((z_{1}+e^{z}z_{1})(z_{2}+e^{z}z_{1}))^{1/2-int}\right| =\mathcal{O}(1).
\]
Next we find the contribution corresponding to the first interval:
\[
\frac{-i\left(e^{n\pi t}+e^{-n\pi t}\right)}{z_{1}^{n}}z_{1}^{1/2+int}(z_{2}-z_{1})^{1/2+int}(2z_{1})^{1/2-int}(z_{2}+z_{1})^{1/2-int}\left(\int_{0}^{\ln^{5}n/n}z^{1/2+int}e^{-nz}dz\right)\left(1+\mathcal{O}(\ln^{5}n/n)\right).
\]
We write the integral in the above expression as 
\begin{align*}
\int_{0}^{\ln^{5}n/n}z^{1/2+int}e^{-nz}dz & =\frac{1}{n^{3/2+int}}\int_{0}^{\ln^{5}n}z^{1/2+int}e^{-z}dz\\
 & =\frac{1}{n^{3/2+int}}(\Gamma(3/2+int)-\int_{\ln^{5}n}^{\infty}z^{1/2+int}e^{-z}dz)\\
 & =\frac{1}{n^{3/2+int}}(\Gamma(3/2+int)+\mathcal{O}(e^{-\ln^{5}n/2})),
\end{align*}
where $\Gamma(\cdot)$ denotes the (complete) gamma function. Using Sterling's approximation to the gamma function 
\[
\Gamma(z)=\exp\left((z+1/2)\Log z-z+\frac{1}{2}\ln2\pi+\mathcal{O}(z^{-1})\right) \qquad (|z| \gg 1)
\]
 and the condition $t\ll\ln^{4}n/n$ we find that
\begin{align*}
|\Gamma(3/2+int)| & =\exp\left(2\ln|3/2+int|-nt\Arg(3/2+int)-\frac{3}{2}-\frac{1}{2}\ln2\pi+\mathcal{O}(1/nt)\right)\\
 & \gg\exp(-\frac{\pi\ln^{4}n}{4}).
\end{align*}
Consequently,
\[
\Gamma(3/2+int)+\mathcal{O}(e^{-\ln^{5}n/2})=\Gamma(3/2+int)(1+o(1)),
\]
and hence
\[
\int_{0}^{\ln^{5}n/n}z^{1/2+int}e^{-nz}dz=\frac{1}{n^{3/2+int}}\Gamma(3/2+int)(1+o(1))=\Gamma(3/2+int)e^{-(3/2+int) \ln n}(1+o(1)).
\]
Assembling the estimates results in the claimed equivalence. Similar arguments show that the second summand in \eqref{eq:sumtwoint} is asymptotic to 
\[
\frac{(-e^{2\pi nt}+e^{-2\pi nt})}{z_{2}^{n}}2^{1/2-int}(z_{2}-z_{1})^{1/2+int}(z_{2}+z_{1})^{1/2-int}\Gamma(3/2+int)e^{-(3/2+int)\ln n}.
\]
Since $t\ll\ln^{4}n/n$ and $z_2>z_1$,
\[
\frac{(-e^{2\pi nt}+e^{-2\pi nt})}{z_{2}^{n}}=\left(\frac{e^{\pi nt}(1+e^{-2\pi nt})}{z_1\left(\frac{z_2}{z_1}\right)^n} \right)\frac{(-e^{\pi nt}+e^{-\pi nt})}{z_{1}^{n-1}}=\frac{(-e^{\pi nt}+e^{-\pi nt})}{z_{1}^{n-1}}\cdot o(1),
\]
 and we conclude that the entire expression \eqref{eq:sumtwoint} is
asymptotic to \eqref{eq:firstint}. We thus obtain the estimate 
\[
\int_{\Gamma_2} f(z,t)dz \sim \frac{-i(e^{n\pi t}+e^{-n\pi t})}{z_{1}^{n-1}}2^{1/2-int}(z_{2}-z_{1})^{1/2+int}(z_{2}+z_{1})^{1/2-int}\Gamma(3/2+int)e^{-(3/2+int)\ln n},
\]
completing the section. We are now in a position to account for the zeros of the polynomials in the Sheffer sequence $(H_n)_{n \gg1}$.

\subsection{\label{sec:zerodist} The zeros of the polynomials $H_n$.}

Recall (c.f. the discussion preceding equation\eqref{eq:intloopcut}) that given $n \in \mathbb{N}$, the polynomial of interest, namely $\pi H_{n}(1/2+int)$, is the imaginary part or $-i$ times the real
part of 
\[
\int_{\Gamma_{2}}f(z,t)dz
\]
when $n$ is even or odd respectively. For large $n$, we now find
a lower bound for the number of real zeros of $H_{n}(1/2+int$) on
$t\in(0,T)$ and compare this number to the degree of $H_n$.

Recall that for each $t\in(0,T)$, $\zeta$ is a solution of \eqref{eq:critpointseq}.
Using this relation we express $t$ as a function of $\zeta$:
\begin{equation}
it=\frac{(\zeta^{2}-z_{1}^{2})(\zeta^{2}-z_{2}^{2})}{2\zeta(z_{1}+z_{2})(\zeta^{2}-z_{1}z_{2})},\label{eq:tfunczeta}
\end{equation}
and rewrite $\int_{\Gamma_{2}}f(z,t)dz$ as a function
of $\zeta$:
\[
h(\zeta):=\int_{\Gamma_{2}}f\left(z,\frac{(\zeta^{2}-z_{1}^{2})(\zeta^{2}-z_{2}^{2})}{2i\zeta(z_{1}+z_{2})(\zeta^{2}-z_{1}z_{2})}\right)dz.
\]
Since zeros of $H_{n}(1/2+int)$ correspond to the intersections of $h(\zeta)$ with the imaginary or the real axis depending on the parity of $n$, we focus on the change of argument of $h(\zeta)$, provided that it does not pass through the origin. Our next results describes ranges of $t$ on which this condition on $h(\zeta)$ is met. In order to ease the exposition, we introduce the following notation.
\begin{notation} Given two functions $f_1, f_2: \mathbb{N} \to \mathbb{R}$, we will denote the set of $t \in(0,T)$ satisfying $f_1(n) \ll t < f_2(n)$ as $n \to \infty$ by
$\langle\langle f_1(n), f_2(n))$. Similarly, we will denote the set of $t \in(0,T)$ satisfying $f_1(n) < t \ll f_2(n)$ as $m \to \infty$ by $( f_1(n), f_2(n) \langle\langle$.
\end{notation}
\begin{lem}
\label{lem:changeargsaddle} Let $\phi, \psi$ and $T$ be defined as in equations \eqref{eq:phidef}, \eqref{eq:psidef} and \eqref{eq:Tdefn}. Let
\begin{align*}
I_1&:= \langle \langle \ln^4 n/n, \  T-\ln^2n/n^{2/3}) \\
I_2&:= \langle \langle \ln^4 n/n, \ T_1-\ln^2n/n^{2/3}) \\
I_3&:=  [T_1+\ln^2 n/n^{2/3}, \ T_2-\ln^2n/n^{2/3}],
\end{align*}
and set
\begin{equation}
g(\zeta):=\frac{2\pi\psi^{2}(\zeta)e^{-2n\phi(\zeta,t)}}{n\phi_{z^{2}}(\zeta,t)}.\label{eq:gzetadef}
\end{equation}
\begin{itemize}
\item[(i)] If $T=T_{1}$, then $h(\zeta(t))\ne0$ for $t \in I_1$.
\item[(ii)] If $T=T_{2}$, then $h(\zeta(t))\ne0$ for $t \in I_2 \cup I_3$.
\end{itemize}
In addition, for $j=1,2,3$, 
\[
2\Delta_{I_j}\arg h(\zeta) =\Delta_{I_j}\arg g(\zeta)+o(1).
\]
\end{lem}

\begin{proof}
We establish the claims in the case $T=T_{2}$ and $t \in I_2$. The remaining cases are obtained using analogous arguments. For $t \in I_2$, equation
\eqref{eq:squareasymp} provides
\[
h^{2}(\zeta)=\frac{2\pi\psi^{2}(\zeta)e^{-2n\phi(\zeta,t)}}{n\phi_{z^{2}}(\zeta,t)}\left(1+\mathcal{O}\left(\epsilon^{2/3}+\frac{e^{-\epsilon^{2}n}}{\epsilon\sqrt{n}}+e^{-\epsilon^{2}n/2}\frac{\sqrt{n\phi_{z^{2}}(\zeta,t)}}{\psi(\zeta)}\right)\right)
\]
for any $\epsilon(n)$ satisfying the conditions $n\epsilon^{2}\gg1$ and 
\[
\frac{\epsilon}{\sqrt{|\phi_{z^{2}}(\zeta,t)|^{3}}t^{2}}=o(1).
\]
Let $\epsilon=\ln n/\sqrt{n}$. Then $n\epsilon^2=n(\ln n/\sqrt{n})^2=\ln^2 n \gg 1$. In addition, since 
\[
\phi_{z^{2}}(\zeta,t)\asymp\begin{cases}
\frac{1}{t} & \text{ for small }t\\
\sqrt{|T_{1}-t|} & \text{ for small }T_{1}-t
\end{cases},
\]
$\psi(\zeta)\asymp\sqrt{t}$ for small $t$, and $t\gg\ln^{4}n/n$, we also get
\[
\frac{\epsilon}{\sqrt{|\phi_{z^{2}}(\zeta,t)|^{3}}t^{2}}=o(1),
\]
as well as 
\[
e^{-\epsilon^{2}n/2}\frac{\sqrt{n\phi_{z^{2}}(\zeta,t)}}{\psi(\zeta)}=o(1).
\]
Consequently, 
\begin{equation}
h^{2}(\zeta)=\frac{2\pi\psi^{2}(\zeta)e^{-2n\phi(\zeta,t)}}{n\phi_{z^{2}}(\zeta,t)}\left(1+o(1)\right)=g(\zeta)(1+o(1))\label{eq:asympsaddle}
\end{equation}
from which the claims $h(\zeta)\ne0$ and $ 2\Delta_{I_2}\arg h(\zeta)=\Delta_{I_2}\arg g(\zeta)+o(1)$ follow.
\end{proof}
\begin{lem}
\label{lem:changeargsmallT-t} Let $T$ be as defined in equation \eqref{eq:Tdefn} and $g(\zeta)$ as defined in \eqref{eq:gzetadef}. \\ Set $I=T-(\ln^2n/n^{2/3}, \ 1/n^{2/3} \langle\langle$. Then $h(\zeta(t))\ne0$ on $I$, and $\Delta_{I}\arg h(\zeta)=\frac{1}{2}\Delta_{I}\arg g(\zeta)+C$, where $|C|<\pi/2+o(1)$. 
\end{lem}

\begin{proof}
The relevant estimate for $h(\zeta)$ in this case (c.f. Section \ref{sec:smallT-t}) is 
\[
\frac{4c\psi(\zeta(T))e^{-n\phi(\zeta,t)}}{\sqrt{6}\sqrt{c^{3}\phi_{z^{3}}(\zeta(T),T)}}\int_{a}^{b}e^{-ny^{2}}\left(\sqrt{3\sqrt{T-t}-(z/c-\zeta/c)}-\frac{\sqrt{T-t}}{\sqrt{3\sqrt{T-t}-(z/c-\zeta/c)}}\right)^{-1}dy(1+o(1)).
\]
Since the real part of the integrand is positive, we immediately obtain
that $h(\zeta)\ne0$. To compute the change of argument of the above expression
over $I$ we first note that 
\[
\Delta_{I}\arg\frac{4c\psi(\zeta(T))e^{-n\phi(\zeta,t)}}{\sqrt{6}\sqrt{c^{3}\phi_{z^{3}}(\zeta(T),T)}}=\Delta_{I}\arg e^{-n\phi(\zeta,t)}.
\]
Next we employ the estimate $\zeta(t)-\zeta(T)=c\sqrt{T-t}+\mathcal{O}(T-t)$ to deduce that 
\begin{equation}
g(\zeta)=\frac{2\pi\psi^{2}(\zeta)e^{-2n\phi(\zeta,t)}}{n\phi_{z^{2}}(\zeta,t)}\sim\frac{2\pi\psi^{2}(\zeta(T))e^{-2n\phi(\zeta,t)}}{n\phi_{z^{3}}(\zeta(T),T)c\sqrt{T-t}},\label{eq:gzetatcloseT}
\end{equation}
and hence 
\[
\Delta_{I}\arg e^{-n\phi(\zeta,t)}=\frac{1}{2}\Delta_{I}\arg g(\zeta)+o(1).
\]

Since 
\begin{equation}
\int_{a}^{b}e^{-ny^{2}}\left(\sqrt{3\sqrt{T-t}-(z/c-\zeta/c)}-\frac{\sqrt{T-t}}{\sqrt{3\sqrt{T-t}-(z/c-\zeta/c)}}\right)^{-1}dy\label{eq:intfactor}
\end{equation}
is in the right half plane for $t \in I$,
the change in argument of this integral over $I$ is no more than the
difference in arguments corresponding to $t=T-\ln^2n/n^{2/3}$ and $t=T-1/n^{2/3}$. When $t=T-\ln^{2}n/n^{2/3}$, the square of the
integral in \eqref{eq:intfactor} is 
\[
h^{2}(\zeta)\left(\frac{16\psi^{2}(\zeta(T))e^{-2n\phi(\zeta,t)}}{6c\phi_{z^{3}}(\zeta(T),T)}\right)^{-1}(1+o(1)).
\]
This expressionsm by \eqref{eq:asympsaddle} is equal to
\[
g(\zeta)\left(\frac{16\psi^{2}(\zeta(T))e^{-2n\phi(\zeta,t)}}{6c\phi_{z^{3}}(\zeta(T),T)}\right)^{-1}(1+o(1)),
\]
whose argument is $o(1)$ by equation \eqref{eq:gzetatcloseT}. Thus, regardless of what the argument is at $t=T-1/n^{2/3}$, the absolute
value of the change in argument of \eqref{eq:intfactor} for $t \in I$ is at most $\pi/2+o(1)$ and the result follows. 
\end{proof}
\begin{lem}
\label{lem:changeargsmallT1-t} Suppose that $T=T_2$ and $h(\zeta(T_{1}))\ne0$. Let 
\begin{align*}
I_1&=(T_{1}-\ln^{2}n/n^{2/3}, \ T_{1}+\ln^{2}n/n^{2/3}), \\
I_2 &=T_1-(\ln^2 n/n^{2/3}, \ 1/n\langle\langle, \quad \textrm{and}\\ 
I_3 &= T_1+\langle \langle 1/n, \ \ln^2n /n^{2/3}).
\end{align*}
 Then
the $h(\zeta(t))\ne0$ on $I_1$, and 
\begin{align*}
 & 2\Delta_{I_1}\arg h(\zeta)\\
= & \Delta_{I_2}\arg g(\zeta)+\Delta_{I_3}\arg g(\zeta)+C
\end{align*}
where $|C|\le\pi+o(1)$. 
\end{lem}

\begin{proof} On the interval $I_1$, we either have $|T_1-t|=o(\frac{1}{n})$ or $\frac{1}{n} \ll |T_1-t|<\ln^2 n/n^{2/3}$.
If $|T_{1}-t|=o\left(\frac{1}{n}\right)$, then the boundedness of 
\[
\Log(z_{1}-z)+\Log(z_{2}-z)-\Log(z_{1}+z)-\Log(z_{2}+z)
\]
on $\Gamma_{2}$ implies that 
\begin{align*}
\phi(z,t) & =\phi(z,T_{1})+i(t-T_{1})\left(\Log(z_{1}-z)+\Log(z_{2}-z)-\Log(z_{1}+z)-\Log(z_{2}+z)\right)\\
 & =\phi(z,T_{1})+o(1/n) \qquad (z \in \Gamma_2).
\end{align*}
Thus 
\begin{align*}
h(\zeta(t)) & =\int_{\Gamma_{2}}e^{-n\phi(z,t)}\psi(z)dz\\
 & =h(\zeta(T_{1}))(1+o(1)),
\end{align*}
from which we deduce that $h(\zeta(t))\ne0$ for $|T_{1}-t|=o(1/n)$
and the change of argument of $h(\zeta(t))$ for $t$ in this range
is $o(1)$. We now consider the case $\frac{1}{n} \ll |T_1-t|<\ln^2 n/n^{2/3}$.
Using an argument completely analogous to the proof of Lemma \ref{lem:changeargsmallT-t},
we conclude that $h(\zeta(t))\ne0$ on $I_2 \cup I_3$, and that
\begin{align*}
 &\Delta_{I_2}\arg h(\zeta)+\Delta_{I_3}\arg h(\zeta)
\le \frac{1}{2} \Delta_{I_2}\arg g(\zeta)+ \frac{1}{2}\Delta_{I_3}\arg g(\zeta)+C,
\end{align*}
where $|C|\le\pi+o(1)$. Combining these estimates establishes the claim.
\end{proof}
\begin{lem}
\label{lem:changeargsmallt} If $\tau$ is a large constant multiple
of $\ln^{4}n/n$, then $h(\zeta(t))\ne0$ for $0<t\le\tau$ and 
\[
\Delta_{0<t\le\tau}\arg h(\zeta(t))=\frac{1}{2}\lim_{\xi\rightarrow0}\Delta_{\xi<t\le\tau}\arg g(\zeta(t))+o(1).
\]
\end{lem}

\begin{figure}
\begin{centering}
\includegraphics[scale=0.3]{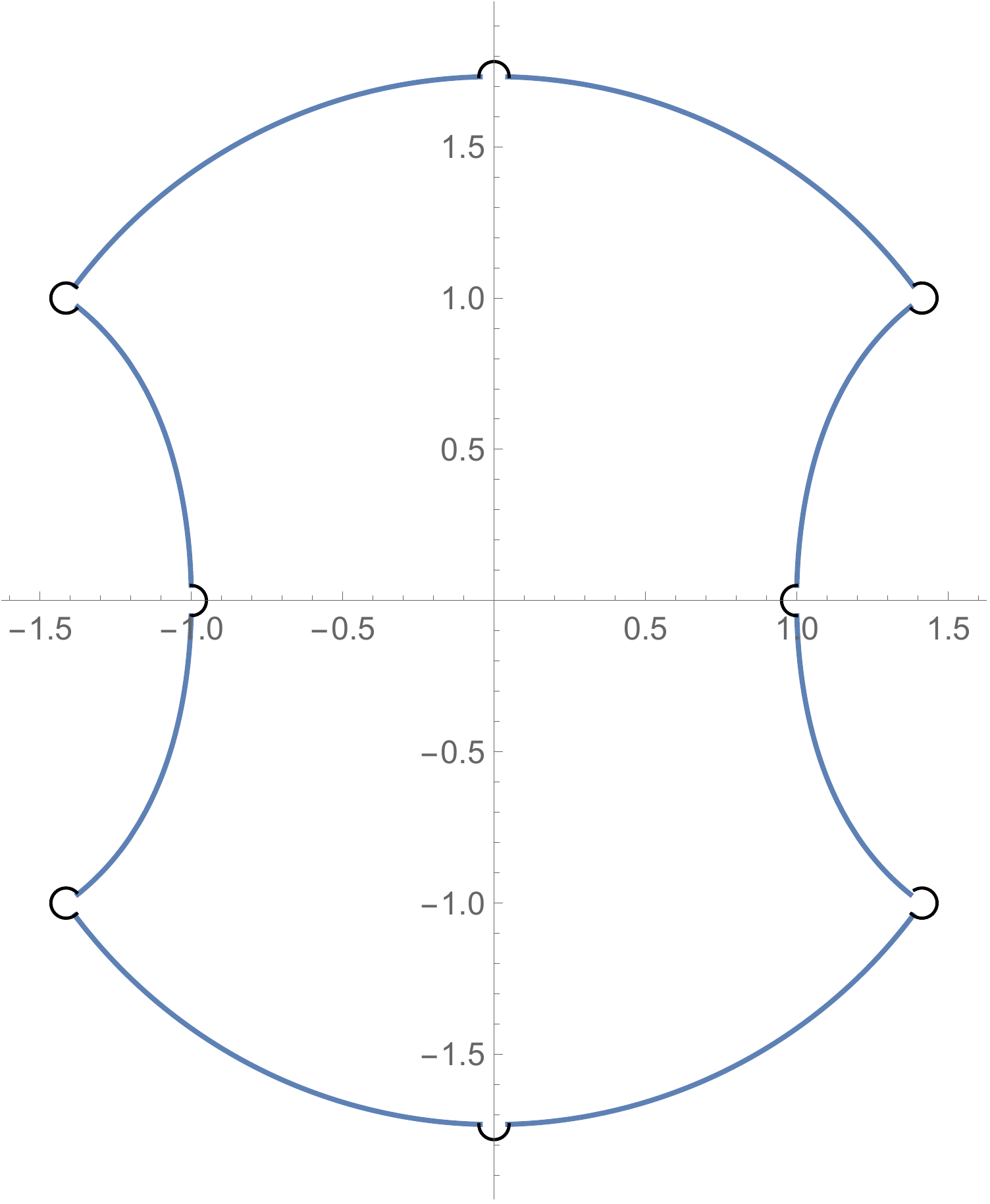} $\qquad$\includegraphics[scale=0.3]{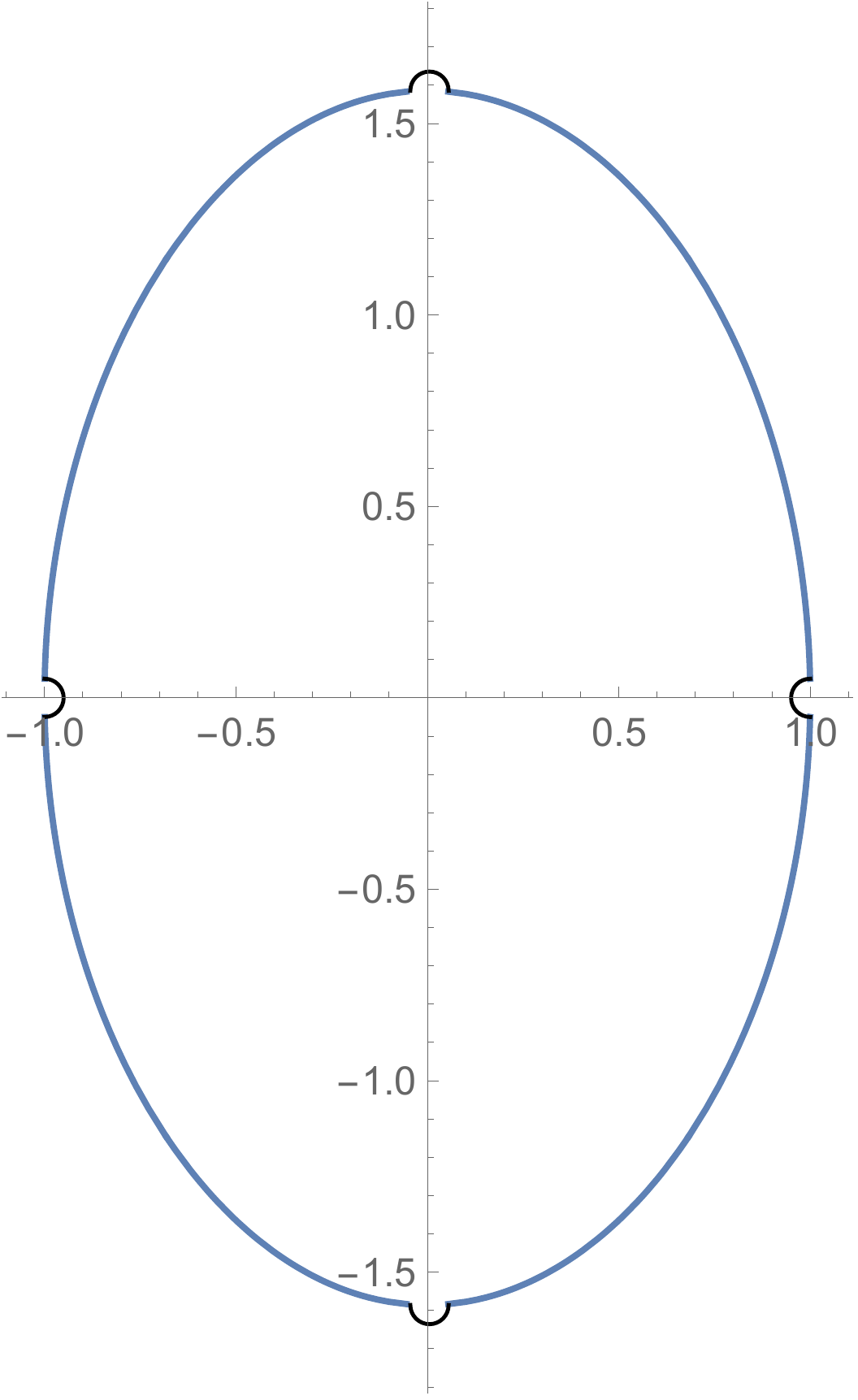} 
\par\end{centering}
\caption{\label{fig:gammacurve}The curve $\gamma$ for for $(z_{1},z_{2})=(1,3)$
(left) and $(1,7)$ (right)}
\end{figure}

\begin{proof}
Recall from Section \ref{sec:smallt} that under the assumptions of the Lemma, 
\begin{equation*}
h(\zeta(t)) \asymp \frac{-i(e^{n\pi t}+e^{-n\pi t})}{z_{1}^{n-1}}2^{1/2-int}(z_{2}-z_{1})^{1/2+int}(z_{2}+z_{1})^{1/2-int}\Gamma(3/2+int)e^{-(3/2+int)\ln n}
\end{equation*}
 for $0<t\le\tau$. It follows that $h(\zeta(t))\ne0$
on this interval. Using Stirling's formula once more we find that the change of the argument $\Gamma(s)$ along the line $L:s=3/2+int$, $0<t\le\tau$
is 
\begin{align*}
 \Delta_{0<t\le\tau} \arg \Gamma(3/2+int)& =\Im\left((1/2+in\tau)\Log(3/2+in\tau)\right)-n\tau+\mathcal{O}\left(\frac{1}{\ln^{4}n}\right)\\
 &= n\tau\ln(n\tau)+\frac{\pi}{4}-n\tau+\mathcal{O}\left(\frac{1}{\ln^{4}n}\right).
\end{align*}
The change of arguments of the factors $2^{1/2-int}$, $(z_{2}-z_{1})^{1/2+int}$,
$(z_{2}+z_{1})^{1/2-int}$, and $e^{-(3/2+int)\ln n}$ are given by $-n\tau\ln2$,
$n\tau\ln(z_{2}-z_{1})$, $-n\tau\ln(z_{2}+z_{1})$, and $-n\tau\ln n$
respectively. Thus, the change of argument of $h(\zeta)$ over the range in question is 
\[
\Delta_{0<t\le\tau} \arg h(\zeta(t))=n\tau\ln\frac{\tau(z_{2}-z_{1})}{2(z_{2}+z_{1})}-n\tau+\frac{\pi}{4}+o(1).
\]
We next consider the change in argument of 
\[
g(\zeta)=\frac{2\pi\psi^{2}(\zeta)e^{-2n\phi(\zeta)}}{n\phi_{z^{2}}(\zeta,t)}=\frac{2\pi(z_{1}-\zeta)^{1+2int}(z_{2}-\zeta)^{1+2int}(z_{1}+\zeta)^{1-2int}(z_{2}+\zeta)^{1-2int}}{n\phi_{z^{2}}(\zeta,t)\zeta^{2n+2}}.
\]
With $z_{1}-\zeta=-iz_{1}t+\mathcal{O}(t^{2})$ (c.f. equation \eqref{eq:z1-zeta}), the change of argument
of the factor $(z_{1}-\zeta(t))^{1+2imt}$ is 
\begin{align*}
 & \lim_{\xi\rightarrow0}\Im\left((1+2imt)(\ln|z_{1}t|-i\frac{\pi}{2}+\mathcal{O}(t))\right)|_{\xi}^{\tau}\\
= & 2m\tau\ln(z_{1}\tau)+o(1),
\end{align*}
while the change of argument of $\zeta(t)^{2n+2}$ is 
\[
\Im\left((2n+2)\Log\zeta(\tau)\right)=(2n+2)(\tau+\mathcal{O}(\tau^{2}))=2n\tau+o(1).
\]
Similarly, the change of argument of the factors $(z_{2}-\zeta(t))^{1+2int}$,
$(z_{1}+\zeta(t))^{1-2int}$, and $(z_{2}+\zeta(t))^{1-2int}$ are given by
$2n\tau\ln(z_{2}-z_{1})+o(1)$, $-2n\tau\ln(2z_{1})+o(1)$, and $-2n\tau\ln(z_{1}+z_{2})$
respectively. Since 
\[
\phi_{z^{2}}(\zeta,t)=-\frac{1}{\zeta^{2}}-it\left(\frac{1}{(z_{1}-\zeta)^{2}}+\frac{1}{(z_{2}-\zeta)^{2}}+\frac{1}{(z_{1}+\zeta)^{2}}+\frac{1}{(z_{2}+\zeta)^{2}}\right),
\]
the corresponding change in argument of $\phi_{z^{2}}(\zeta,t)$ for
$0<t\le\tau$ is $\pi/2+o(1)$. We conclude that the change in argument
of $g(\zeta(t))$ is 
\[
2n\tau\ln\frac{\tau(z_{2}-z_{1})}{2(z_{1}+z_{2})}-2n\tau+\frac{\pi}{2}+o(1)
\]
and the result follows. 
\end{proof}
Lemmas \ref{lem:changeargsaddle}, \ref{lem:changeargsmallT-t}, \ref{lem:changeargsmallT1-t},
and \ref{lem:changeargsmallt} show that in order to understand the change in the argument of $h(\zeta(t))$, it suffices to study of the change
of argument of $g(\zeta)$ defined in \eqref{eq:gzetadef}, where we (re)write $\phi_{z^2}$ as
\[
\phi{}_{z^{2}}(\zeta,t)=\frac{(\zeta^{2}+z_{1}z_{2})(\zeta^{4}+(z_{1}^{2}-4z_{1}z_{2}+z_{2}^{2})\zeta^{2}+z_{1}^{2}z_{2}^{2})}{\zeta^{2}(\zeta^{2}-z_{1}^{2})(\zeta^{2}-z_{2}^{2})(\zeta^{2}-z_{1}z_{2})}.
\]
Let $\gamma$ be the simple closed curve with counter clockwise orientation
formed by the traces of $\zeta(t),\overline{\zeta(t)},-\zeta(t),$
and $-\overline{\zeta(t)}$ for $0\le t\le T$ and small deformations
around 
\begin{equation}
\begin{cases}
\pm i\sqrt{z_{1}z_{2}},\pm\zeta(T_{1}),\pm\overline{\zeta(T_{1})} & \text{ if }z_{1}^{2}-6z_{1}z_{2}+z_{2}^{2}<0\\
\pm i\zeta(T) & \text{ if }z_{1}^{2}-6z_{1}z_{2}+z_{2}^{2}\ge0
\end{cases}\label{eq:zerosphi''}
\end{equation}
such that the region enclosed by $\gamma$ contains the points defined in \eqref{eq:zerosphi''}.
We also deform $\gamma$ around $\pm z_{1}$ so that the cuts $(-\infty,-z_{1}]$
and $[z_{1},\infty)$ lie outside this region (see Figure \ref{fig:gammacurve}).

By computing the logarithmic derivative of $g(\zeta)$, we find that
\begin{align*}
\Res\left(\frac{g'(\zeta)}{g(\zeta)},0\right) & =\Res\left(-2n\frac{d\phi(\zeta,t)}{d\zeta},0\right),
\end{align*}
since $\frac{\psi^2(\zeta)}{\phi_{z^2}(\zeta,t)}$ is analytic in a neighborhood of the origin. Furthermore,
\[
\frac{d\phi(\zeta,t)}{d\zeta}=\frac{1}{\zeta}-\frac{d}{d\zeta}\left(\frac{(\zeta^{2}-z_{1}^{2})(\zeta^{2}-z_{2}^{2})}{2\zeta(z_{1}+z_{2})(\zeta^{2}-z_{1}z_{2})}.\left(\frac{1}{\zeta-z_{1}}+\frac{1}{\zeta-z_{2}}-\frac{1}{\zeta+z_{1}}-\frac{1}{\zeta+z_{2}}\right)\right),
\]
and hence $d\phi(\zeta,t)/d\zeta-1/\zeta$ is the derivative of
a meromorphic function around the origin.Thus, its residue at $0$
is $0$ by its Laurent series expansion. We conclude that 
\[
\Res\left(\frac{g'(\zeta)}{g(\zeta)},0\right)=-2n\Res\left(\frac{d\phi(\zeta,t)}{d\zeta},0\right)=-2n.
\]
Since the points in \eqref{eq:zerosphi''} are all simple poles of $g(\zeta)$ with residue $1$, we conclude from
the residue theorem that 
\begin{equation}
\frac{1}{2\pi i}\int_{\gamma}\frac{g'(\zeta)}{g(\zeta)}d\zeta=\begin{cases}
-2n-6 & \text{ if }z_{1}^{2}-6z_{1}z_{2}+z_{2}^{2}<0\\
-2n-2 & \text{ if }z_{1}^{2}-6z_{1}z_{2}+z_{2}^{2}\ge0
\end{cases}.\label{eq:argchangegzeta}
\end{equation}

Suppose that $\zeta$ not a singularity of $g(\zeta)$. Then $g(\overline{\zeta})=\overline{g(\zeta)}$. In addition, since 
\[
e^{-2n\Log\zeta}=\frac{1}{\zeta^{2n}},
\]
we also have
\begin{equation}
g(-\zeta)=g(\zeta).\label{eq:gsymyaxis}
\end{equation}

If $\zeta\in i\mathbb{R}$, then $\zeta=-\overline{\zeta}$. In this
case we apply \eqref{eq:gsymyaxis} to conclude that 
\[
g(\zeta)=g(-\overline{\zeta})=g(\overline{\zeta})=\overline{g(\zeta)},
\]
and consequently $g(\zeta)\in\mathbb{R}$.

If we let $\gamma_{1}$, $\gamma_{2}$, $\gamma_{3}$, and $\gamma_{4}$
be the portion of $\gamma$ on the first, second, third, and fourth
quadrant respectively, then -- using the above relations -- we see that
\[
\Delta_{\gamma_{4}}\arg g(\zeta)=-\Delta_{\gamma_{4}}\arg\overline{g(\zeta)}=-\Delta_{\gamma_{4}}\arg g(\overline{\zeta})=\Delta_{\gamma_{1}}\arg g(\zeta),
\]
and 
\[
\Delta_{\gamma_{3}}\arg g(\zeta)=\Delta_{\gamma_{2}}\arg g(\zeta).
\]
Similarly 
\begin{align*}
\Delta_{\gamma_{3}}\arg g(\zeta) & =\Delta_{\gamma_{3}}\arg g(-\zeta)\\
 & =\Delta_{\gamma_{1}}\arg g(\zeta).
\end{align*}
We conclude that 
\begin{equation}
\Delta_{\gamma_{1}}\arg g(\zeta)=\frac{\Delta_{\gamma}\arg g(\zeta)}{4}=\begin{cases}
-(n+3)\pi & \text{ if }z_{1}^{2}-6z_{1}z_{2}+z_{2}^{2}<0\\
-(n+1)\pi & \text{ if }z_{1}^{2}-6z_{1}z_{2}+z_{2}^{2}\ge0
\end{cases}.\label{eq:arggzeta}
\end{equation}
The following three lemmas demonstrate that the change in the argument of $g(\zeta(t))$ near the points $T_1$ and $T_2$ are small. 
\begin{lem}
\label{lem:taud} Let $g(\zeta)$ be as defined in equation \eqref{eq:gzetadef}, and let $T$ be as defined in \eqref{eq:Tdefn}. If $\tau<T$ satisfies $T-\tau\ll1/n^{2/3}$, then
\[
\lim_{\xi\rightarrow0}\Delta_{\tau\le t<T-\xi}\arg g(\zeta(t))\ll1.
\]
\end{lem}
\begin{proof}
Using the definition of $g(\zeta)$ and the estimates $\phi_{z^{2}}(\zeta,t)\asymp\zeta(t)-\zeta(T)\asymp\sqrt{T-t}$,
it suffices to show that
\[
\Im\phi(\zeta(\tau),\tau)-\Im\phi(\zeta(T),T)\ll\frac{1}{n}.
\]
With $\zeta(T)-\zeta(\tau)=c\sqrt{T-\tau}+\mathcal{O}(T-\tau)$, $\phi_{z}(\zeta(T),T)=0$,
$\phi_{z^{2}}(\zeta(T),T)=0$, and $\Im\phi_{t}(\zeta(T),T)=0$, we
expand the bivariate function $\phi(z,t)$ in a Taylor series 
in a neighborhood of $(\zeta(T),T)$ to arrive at
\[
\Im\phi(\zeta(\tau),\tau)-\Im\phi(\zeta(T),T)=\Im\left((\zeta(T)-\zeta(\tau))(T-\tau)\phi_{z,t}(\zeta(T),T)+\mathcal{O}((T-\tau))\right)\asymp(T-\tau)^{3/2},
\]
and since $(T-\tau)^{3/2}\ll\frac{1}{n}$, the result follows.
\end{proof}
 Since the proofs of the next two lemmas are essentially identical to the one we just gave, we omit them, and state only the results.

\begin{lem}
\label{lem:changeargT1+} Let $T=T_2$, and suppose that $\tau>T_{1}$ satisfies
$\tau-T_{1}\ll1/n$. Then 
\[
\lim_{\xi \to 0}\Delta_{T_1+\xi <t<\tau }\arg g(\zeta(t))\ll1.
\]
\end{lem}
\begin{lem}
\label{lem:changeargT1-} Let $T=T_2$ and suppose tha t$\tau<T_{1}$ satisfies
$T_{1}-\tau\ll1/n$. Then 
\[
\lim_{\xi\rightarrow0}\Delta_{\tau\le t< T_{1}-\xi}\arg g(\zeta(t))\ll1.
\]
\end{lem}
Next, we address the  change in the argument of $g(\zeta)$ on the small deformations near the points $z_1, \zeta(T_1)$ and $\zeta(T_2)$.
We begin with the small arc of $\gamma_{1}$ around $z_{1}$. Note that for
any fixed $n$, as $\zeta\rightarrow z_{1}$, 
\begin{align*}
\psi^{2}(\zeta) &\sim C_1( z_{1}-\zeta),\\
\phi_{z^2}(\zeta) &\sim C_2/(z_{1}-\zeta), \qquad \textrm{and}\\
\exp\left(2n(\zeta-z_{1})(z_{1}-z_{2})\Log(z_{1}-\zeta)/z_{1}\right)&\rightarrow1
\end{align*}
for some constants $C_1,C_2$. We conclude that as $\zeta\rightarrow z_{1}$,
\[
g(\zeta)=\frac{2\pi\psi^{2}(\zeta)e^{-2n\phi(\zeta)}}{n\phi_{z^2}(\zeta)} \sim C_3 (z_{1}-\zeta)^{2}, \qquad (C_3 \in \mathbb{C})
\]
 Using the estimate $z_{1}-\zeta=-iz_{1}t+\mathcal{O}(t^{2})$ developed for small $t$ in equation \eqref{eq:z1-zeta}, we find the change of argument of $g(\zeta)$
on the small arc of $\gamma_{1}$ around $z_{1}$ to be $-\pi+o(1)$.

We continue by considering the small arc of $\gamma_1$ around $\zeta(T)$. Using Lemma \ref{lem:asympzetaT} we deduce that the change of argument
of $\zeta(t)$ on the arc of $\gamma_{1}$ around $\zeta(T)$ approaches
$\pi/2$ when the arc is small. Since $g(\zeta)$ has a simple pole
at $\zeta(T)$, the change of argument of $g(\zeta)$ on the piece of
$\gamma_{1}$ around $\zeta(T)$ approaches $-\pi/2$. 

Finally, we look at the change of argument in $g(\zeta)$ on the small arc near $T_1$ when $T=T_{2}$. If $T=T_2$, then equations \eqref{eq:dform} and \eqref{eq:eform} imply that $d/\widehat{d}=-i$. Consequently, the change of argument of $g(\zeta)$ on the arc
of $\gamma_{1}$ around $\zeta(T_{1})$ approaches $-3\pi/2$.

We are now in position to count the number of zeros of our polynomials $H_n$ -- thereby completing the proof of the main result -- via bounding the change in the argument of $h(\zeta(t))$ from below on the interval $(0,T)$. There are two cases to consider, corresponding to whether $T=T_1$, or $T=T_2$. 

If $T=T_{1}$, then the observations in the preceding paragraphs along with
Lemmas \ref{lem:changeargsmallT-t} and \ref{lem:changeargsmallt}
imply that for some $|C|<\pi/2+o(1)$, 
\begin{align*}
\Delta_{0<t<T}\arg h(\zeta) & =\frac{1}{2}\left(\Delta_{\gamma_{1}}\arg g(\zeta)+\frac{3\pi}{2}\right)+C\\
 & \stackrel{\eqref{eq:arggzeta}}{=}-\frac{n\pi}{2}+\frac{\pi}{4}+C.
\end{align*}
Since $\pi H_{n}(1/2+int)$ is either the imaginary part or $-i$ times the real
part of $h(\zeta)$, the number of zeros of $H_{n}(1/2+int)$ on $(0,T)$
is at least 
\[
\left\lfloor \frac{\Delta_{0<t<T}\arg h(\zeta)}{\pi}\right\rfloor \ge\left\lfloor \frac{n}{2}-\frac{3}{4}+o(1)\right\rfloor .
\]
It follows that $H_{n}(x)$ has at least 
\[
2\left\lfloor \frac{n}{2}-\frac{3}{4}+o(1)\right\rfloor =\begin{cases}
n-3 & \text{ if }2\nmid n\\
n-2 & \text{ if }2\mid n
\end{cases}
\]
nonreal zeros on the line $\Re x=1/2$. It remains to account for the missing $2$ or $3$ zeros depending on the parity of $n$. Once this is accomplished (see Lemma \ref{lem:trivialzeros}), and we establish a bound on the degree of $H_n$ (see Lemma \ref{lem:degree}), Theorem \ref{thm:maintheorem} will follow from the fundamental theorem of algebra. 
\begin{lem}
\label{lem:trivialzeros}If $n>2$ is even, then 
\begin{enumerate}
\item $x=0$ and $x=1$ are zeros of $H_{n}(x)$ 
\item $H'_{n}(0)<0$ and $H_{n}'(1)>0$. 
\end{enumerate}
If $n>1$ is odd, then 
\begin{enumerate}
\item $x=0$ , $x=1/2$, and $x=1$ are zeros of $H_{n}(x)$ 
\item $H'_{n}(0)<0$ and $H'_{n}(1)<0$. 
\end{enumerate}
\end{lem}

\begin{proof}
With the substitution $x$ by $1-x$ and $z$ by $-z$ in \eqref{eq:genfunc}
we conclude 
\[
H_{n}(x)=(-1)^{n}H_{n}(1-x).
\]
Consequently $H_{n}(1/2)=0$ when $n$ is odd and it suffices to consider
the case $x=0$. Plugging $x=0$  into the right side of equation \eqref{eq:genfunc}
gives $H_{n}(0)=0$ $\forall n > 2$. Similarly, evaluating the derivative
of this expression as a function of $x$ at $x=0$ gives 
\begin{align*}
\sum_{n=0}^{\infty}H_{n}'(0)\frac{z^{n}}{n!} & =(z_{1}+z)(z_{2}+z)\left(\Log(z_{1}-z)+\Log(z_{2}-z)-\Log(z_{1}+z)-\Log(z_{2}+z)\right)\\
 & =(z_{1}+z)(z_{2}+z)\left(\Log\left(1-\frac{z}{z_{1}}\right)-\Log\left(1+\frac{z}{z_{1}}\right)+\Log\left(1-\frac{z}{z_{2}}\right)-\Log\left(1+\frac{z}{z_{2}}\right)\right).
\end{align*}
Since all the coefficients in the Taylor series of the right side
is positive, we conclude $H'_{n}(0)>0$. 
\end{proof}
\begin{lem}
\label{lem:degree} Suppose that $(H_n(x))_{n \geq 0}$ is as in the statement of Theorem \ref{thm:maintheorem}. Then for each $n$, the degree of $H_{n}(x)$ is at most $n$. Furthermore,
$\lim_{x\rightarrow-\infty}H_{n}(x)=+\infty$ and 
\[
\lim_{s\rightarrow+\infty}H_{n}(x)=\begin{cases}
+\infty & \text{ if }2|n\\
-\infty & \text{ if }2\nmid n.
\end{cases}
\]
\end{lem}

\begin{proof}
We apply the binomial expansion to each factor of 
\[
Q(z)^{x}Q(-z)^{1-x}=z_{1}z_{2}(1-z/z_{1})^{x}(1-z/z_{2})^{x}(1+z/z_{1})^{1-x}(1+z/z_{2})^{1-x}
\]
and collect the $z^{n}$-coefficients to conclude that the degree
of $H_{n}(x)$ is at most $n$. Also from this binomial expansion,
we see that all the coefficients in the power series in $z$ of each
factor are positive as $x\rightarrow-\infty$. Thus $\lim_{x\rightarrow-\infty}H_{n}(x)=+\infty$.
We complete the proof by applying the identity $H_{n}(x)=(-1)^{n}H_{n}(1-x)$. 
\end{proof}
\noindent With these results, the proof of Theorem \ref{thm:maintheorem} is complete in the case when $T=T_1$.

If $T=T_{2}$ and $h(\zeta(T_{1}))\ne0$, Lemmas \ref{lem:changeargsmallT-t},
\ref{lem:changeargsmallT1-t}, and \ref{lem:changeargsmallt} imply
that for some $|C|<3\pi/2+o(1)$ 
\begin{align*}
\Delta_{0<t<T}\arg h(\zeta) & =\frac{1}{2}\left(\Delta_{\gamma_{1}}\arg g(\zeta)+3\pi\right)+C\\
 & \stackrel{\eqref{eq:arggzeta}}{=}-\frac{n\pi}{2}+C.
\end{align*}
Consequently, the number of zeros of $H_{n}(1/2+int)$ on $(0,T)$
is at least 
\[
\left\lfloor \frac{\Delta_{0<t<T}\arg h(\zeta)}{\pi}\right\rfloor \ge\left\lfloor \frac{n}{2}-\frac{3}{2}+o(1)\right\rfloor .
\]
and $H_{n}(x)$ has at least 
\[
2\left\lfloor \frac{n}{2}-\frac{3}{2}+o(1)\right\rfloor =\begin{cases}
n-5 & \text{ if }2\nmid n\\
n-4 & \text{ if }2\mid n
\end{cases}
\]
non-real zeros on the line $\Re x=1/2$. Since $H_{n}(x)$ has zeros
at $x=0$ and $x=1$ and a zero at $x=1/2$ when $n$ is odd, there
are at most two possible zeros of $H_{n}(x)$ distinct from $0,1$
and not on the line $\Re x=1/2$. By the symmetry of zeros of $H_{n}(x)$
along the line $\Re x=1/2$ and the real axis, these two zeros must
be real and symmetric about the line $\Re x=1/2$. 

If $n$ is even,
then Lemmas \ref{lem:trivialzeros} and \ref{lem:degree} imply that
$H_{n}(+\infty)H_{n}(0^{-})>0$, $H_{n}(0^{+})H_{n}(1^{-})>0$,
$H_{n}(1^{+})H_{n}(+\infty)>0$ and these two possible exceptional
zeros do not exist. A similar argument applies to the case $n$ is
odd.

If $T=T_{2}$ and $h(\zeta(T_{1}))=0$, then the number of real zeros
of $H_{n}(1/2+int)$ on $(0,T)\backslash\{T_{1}\}$ is at least 
\begin{align*}
 & \lim_{\xi\rightarrow0}\left\lfloor \frac{\Delta_{0<t<T_{1}-\xi}\arg h(\zeta)}{\pi}\right\rfloor +\left\lfloor \frac{\Delta_{T_{1}+\xi<t<T}\arg h(\zeta)}{\pi}\right\rfloor \\
\ge & \lim_{\xi\rightarrow0}\left\lfloor \frac{\Delta_{0<t<T_{1}-\xi}\arg h(\zeta)}{\pi}+\frac{\Delta_{T_{1}+\xi<t<T}\arg h(\zeta)}{\pi}\right\rfloor -1.
\end{align*}
If we count $T_{1}$ as another zero of $H_{n}(1/2+int)$ on $(0,T)$,
we obtain the same number of zeros of this polynomial as in the case
$h(\zeta(T_{1}))\ne0$. This completes the proof of Theorem \ref{thm:maintheorem}.

\end{document}